\documentclass[11pt]{article}

\usepackage[scale = 0.703]{geometry}
\usepackage[affil-it]{authblk}

\usepackage[dvipsnames]{xcolor}
\usepackage{tikz}

\usepackage{amsthm}
\usepackage{mathtools}

\newtheorem{assumption}{Assumption}
\newtheorem{theorem}{Theorem}
\newtheorem{lemma}{Lemma}
\newtheorem{corollary}{Corollary}

\usepackage{microtype}
\usepackage{amssymb}
\usepackage{bm}
\usepackage{bbm}

\usepackage{hyperref}
\usepackage{enumitem}

\usepackage[backend = biber, natbib = true, style = authoryear-comp]{biblatex}

\newcommand{\CustomQED}{\hfill \qedsymbol}
\newenvironment{CustomProof}{\par\noindent{\textbf{Proof}\ }}{\CustomQED \medskip}

\newcommand*{\E}{\mathbb{E}}
\newcommand*{\R}{\mathbb{R}}
\newcommand*{\bbX}{\mathbb{X}}

\newcommand*{\indic}{\mathbbm{1}}

\newcommand*{\bfone}{\mathbf{1}}
\newcommand*{\bfzero}{\mathbf{0}}

\newcommand*{\bfa}{\mathbf{a}}
\newcommand*{\bfe}{\mathbf{e}}
\newcommand*{\bfu}{\mathbf{u}}
\newcommand*{\bfv}{\mathbf{v}}
\newcommand*{\bfw}{\mathbf{w}}
\newcommand*{\bfx}{\mathbf{x}}

\newcommand*{\bfbeta}{\bm{\beta}}
\newcommand*{\bfxi}{\bm{\xi}}
\newcommand*{\bfeps}{\bm{\varepsilon}}

\newcommand*{\bfX}{\mathbf{X}}
\newcommand*{\bfY}{\mathbf{Y}}

\newcommand*{\calN}{\mathcal{N}}

\newcommand*{\inv}{^{-1}}
\newcommand*{\lin}{_{\mathrm{lin}}}
\newcommand*{\tran}{^{\top}}

\newcommand*{\eps}{\varepsilon}
\newcommand*{\simiid}{\overset{\scriptscriptstyle i.i.d.}{\sim}}
\newcommand*{\lambdamin}{\lambda_{\mathrm{min}}}

\newcommand*{\Ber}{\mathrm{Ber}}
\newcommand*{\Cov}{\mathrm{Cov}}
\newcommand*{\Diag}{\mathrm{Diag}}
\newcommand*{\Tr}{\mathrm{Tr}}

\newcommand*{\id}{\mathrm{id}}
\newcommand*{\op}{\mathrm{op}}

\newcommand*{\ol}[1]{\overline{#1}}
\newcommand*{\wh}[1]{\widehat{#1}}
\newcommand*{\wt}[1]{\widetilde{#1}}

\newcommand*{\rpA}{A^{\mathrm{rp}}}
\newcommand*{\whbeta}{\wh{\bfbeta}}
\newcommand*{\wtbeta}{\wt{\bfbeta}}
\newcommand*{\wtbetaaff}{\wt{\bfbeta}_{\mathrm{aff}}}
\newcommand*{\rpbeta}{\wtbeta^{\mathrm{rp}}}
\newcommand*{\starbeta}{\bfbeta_{\star}}

\DeclarePairedDelimiter{\abs}{\lvert}{\rvert}
\DeclarePairedDelimiter{\norm}{\lVert}{\rVert}

\DeclareMathOperator*{\argmin}{arg\,min}

\title{Dropout Regularization Versus $\ell_2$-Penalization\newline in the Linear
Model}

\author{Gabriel Clara
  \footnote{Corresponding author: \texttt{g.clara@utwente.nl}}
  \footnote{The authors thank an anonymous editor and three anonymous referees
    for their valuable time and effort; their comments improved the article
    tremendously.
    
    \vspace{6pt}

    \noindent \begin{minipage}{0.05\textwidth}
      \includegraphics[height = 2.5\baselineskip]{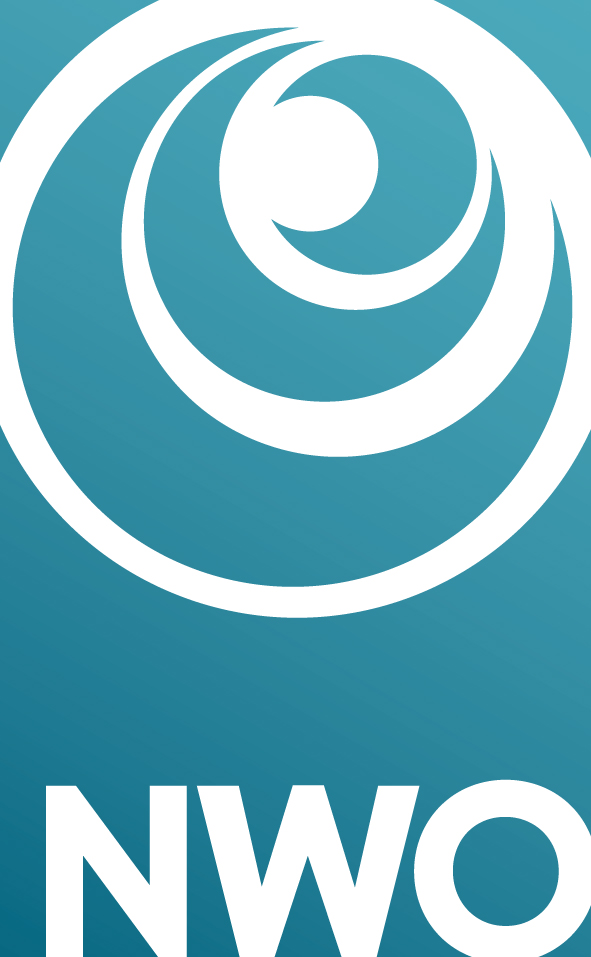}
    \end{minipage} \hfill \begin{minipage}{0.935\textwidth}
      This publication is part of the project \textit{Statistical foundation for
      multilayer neural networks} (project number VI.Vidi.192.021 of the Vidi
      ENW programme) financed by the Dutch Research Council (NWO).
    \end{minipage}}
}
\author{Sophie Langer \protect\footnotemark[2] }
\author{Johannes Schmidt-Hieber \protect\footnotemark[2] }

\affil{Faculty of Electrical Engineering, Mathematics, and Computer Science \\
  University of Twente}

\begin{document}

\maketitle

\begin{abstract}
  We investigate the statistical behavior of gradient descent iterates with
  dropout in the linear regression model. In particular, non-asymptotic bounds
  for the convergence of expectations and covariance matrices of the iterates
  are derived. The results shed more light on the widely cited connection
  between dropout and $\ell_2$-regularization in the linear model. We indicate a
  more subtle relationship, owing to interactions between the gradient descent
  dynamics and the additional randomness induced by dropout. Further, we study a
  simplified variant of dropout which does not have a regularizing effect and
  converges to the least squares estimator.
\end{abstract}


\section{Introduction}

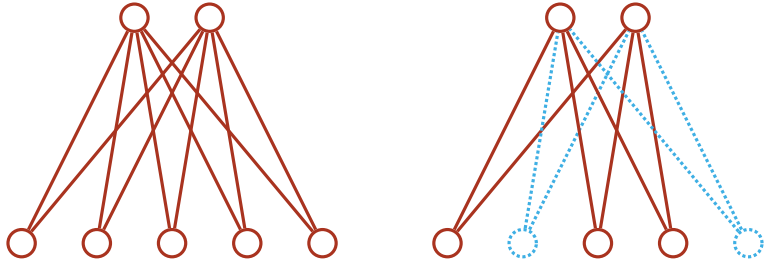
\begin{figure}
  \centering
  \begin{tikzpicture}[very thick]
    \node[shape=circle, draw=Mahogany] (N1) at (0.5, 1) {};
    \node[shape=circle, draw=Mahogany] (N2) at (-0.5, 1) {};

    \node[shape=circle, draw=Mahogany] (I1) at (-2, -2) {};
    \node[shape=circle, draw=Mahogany] (I2) at (-1, -2) {};
    \node[shape=circle, draw=Mahogany] (I3) at (0, -2) {};
    \node[shape=circle, draw=Mahogany] (I4) at (1, -2) {};
    \node[shape=circle, draw=Mahogany] (I5) at (2, -2) {};

    \path[-, draw=Mahogany] (I1) edge node {} (N1);
    \path[-, draw=Mahogany] (I2) edge node {} (N1);
    \path[-, draw=Mahogany] (I3) edge node {} (N1);
    \path[-, draw=Mahogany] (I4) edge node {} (N1);
    \path[-, draw=Mahogany] (I5) edge node {} (N1);

    \path[-, draw=Mahogany] (I1) edge node {} (N2);
    \path[-, draw=Mahogany] (I2) edge node {} (N2);
    \path[-, draw=Mahogany] (I3) edge node {} (N2);
    \path[-, draw=Mahogany] (I4) edge node {} (N2);
    \path[-, draw=Mahogany] (I5) edge node {} (N2);
  \end{tikzpicture}
  \hspace{1cm}
  \begin{tikzpicture}[very thick]
    \node[shape=circle, draw=Mahogany] (N1) at (0.5, 1) {};
    \node[shape=circle, draw=Mahogany] (N2) at (-0.5, 1) {};

    \node[shape=circle, draw=Mahogany] (I1) at (-2, -2) {};
    \node[shape=circle, densely dotted, draw=CornflowerBlue] (I2) at (-1, -2) {};
    \node[shape=circle, draw=Mahogany] (I3) at (0, -2) {};
    \node[shape=circle, draw=Mahogany] (I4) at (1, -2) {};
    \node[shape=circle, densely dotted, draw=CornflowerBlue] (I5) at (2, -2) {};

    \path[-, draw=Mahogany] (I1) edge node {} (N1);
    \path[-, densely dotted, draw=CornflowerBlue] (I2) edge node {} (N1);
    \path[-, draw=Mahogany] (I3) edge node {} (N1);
    \path[-, draw=Mahogany] (I4) edge node {} (N1);
    \path[-, densely dotted, draw=CornflowerBlue] (I5) edge node {} (N1);

    \path[-, draw=Mahogany] (I1) edge node {} (N2);
    \path[-, densely dotted, draw=CornflowerBlue] (I2) edge node {} (N2);
    \path[-, draw=Mahogany] (I3) edge node {} (N2);
    \path[-, draw=Mahogany] (I4) edge node {} (N2);
    \path[-, densely dotted, draw=CornflowerBlue] (I5) edge node {} (N2);
  \end{tikzpicture}
  \caption{Regular neurons (left) with all connections active. Sample of the
  same neurons with dropout (right). The dashed connections are ignored during
  the current iteration of training.}
  \label{fig::drop_neuron}
\end{figure}

Dropout is a simple, yet effective, algorithmic regularization technique,
intended to prevent neural networks from overfitting. First introduced in
\cite{srivastava_et_al_2014}, the method is implemented via random masking of
neurons at training time. Specifically, during every gradient descent iteration,
the output of each neuron is replaced by zero based on the outcome of an
independently sampled $\Ber(p)$-distributed variable. This temporarily removes
each neuron with a probability of $1 - p$, see Figure \ref{fig::drop_neuron} for
an illustration. The method has demonstrated effectiveness in various
applications, see for example \cite{krizhevsky_et_al_2012,
srivastava_et_al_2014}. On the theoretical side, dropout is often studied by
exhibiting connections with explicit regularizers \citep{arora_et_al_2021,
baldi_sadowski_2013, cavazza_et_al_2018, mcallester_2013, mianjy_arora_2019,
mianjy_et_al_2018, senen-cerda_sanders_2020, srivastava_et_al_2014,
wager_et_al_2013}. Rather than analyzing the gradient descent iterates with
dropout noise, these results consider the marginalized training loss with
marginalization over the dropout noise. Within this framework,
\cite{srivastava_et_al_2014} established a connection between dropout and
weighted $\ell_2$-penalization in the linear regression model. This connection
is now cited in popular textbooks \citep{efron_hastie_2016,
goodfellow_et_al_2016}.

However, \cite{wei_et_al_2020} show empirically that injecting dropout noise in
the gradient descent iterates also induces an implicit regularization effect
that is not captured by the link between the marginalized loss and explicit
regularization. This motivates our approach to directly derive the statistical
properties of gradient descent iterates with dropout. We study the linear
regression model due to mathematical tractability and because the minimizer of
the explicit regularizer is unique and admits a closed-form expression. In line
with the implicit regularization observed in \cite{wei_et_al_2020}, our main
result provides a theoretical bound quantifying the amount of randomness in the
gradient descent scheme that is ignored by the previously considered minimizers
of the marginalized training loss. More specifically, Theorem \ref{thm::var_lim}
shows that for a fixed learning rate there is additional randomness which fails
to vanish in the limit, while Theorem \ref{thm::wtbeta_var_sub} characterizes
the gap between dropout and $\ell_2$-penalization with respect to the learning
rate, the dropout parameter $p$, the design matrix, and the distribution of the
initial iterate. Theorem \ref{thm::rp_l2_conv} shows that this gap disappears
for the Ruppert-Polyak averages of the iterates.

To provide a clearer understanding of the interplay between gradient descent and
variance, we also investigate a simplified variant of dropout featuring more
straightforward interactions between the two. Applying the same analytical
techniques to this simplified variant, Theorem \ref{thm::wh_l2_conv} establishes
convergence in quadratic mean to the conventional linear least-squares
estimator. This analysis illustrates the sensitivity of gradient descent to
small changes in the way noise is injected during training.

Many randomized optimization methods can be formulated as noisy gradient descent
schemes. The developed strategy to treat gradient descent with dropout may be
generalized to other settings. An example is the recent analysis of forward
gradient descent in \cite{bos_schmidt-hieber_2023}.

The article is organized as follows. After discussing related results below,
Section \ref{sec::drop_variations} contains preliminaries and introduces two
different variants of dropout. Section \ref{sec::avg_drop} discusses some
extensions of previous results for averaged dropout obtained by marginalizing
over the dropout distribution in the linear model considered in
\cite{srivastava_et_al_2014}. Section \ref{sec::gd_results} illustrates the main
results on gradient descent with dropout in the linear model, and examines its
statistical optimality. Section \ref{sec::disc} contains further discussion and
mentions a number of natural follow-up problems. All proofs are deferred to the
Appendix.

\subsection{Other Related Work}

Considering linear regression and the marginalized training loss with
marginalization over the dropout noise, the initial dropout article
\citep{srivastava_et_al_2014} already connects dropout with
$\ell_2$-regularization. This connection was also noted by
\cite{baldi_sadowski_2013} and by \cite{mcallester_2013}. As this argument is
crucial in our own analysis, we will discuss it in more detail in Section
\ref{sec::avg_drop}.

\cite{wager_et_al_2013} extends the reasoning to generalized linear models and
more general forms of injected noise. Employing a quadratic approximation to the
loss function after marginalization over the injected noise, the authors exhibit
an explicit regularizer. In case of dropout noise, this regularizer induces in
first-order an $\ell_2$-penalty after rescaling of the data by the estimated
inverse of the diagonal Fisher information. 

For two-layer models, marginalizing the dropout noise leads to a nuclear norm
penalty on the product matrix, both in matrix factorization
\citep{cavazza_et_al_2018} and linear neural networks \citep{mianjy_et_al_2018}.
The latter may be seen as a special case of a particular ``$\ell_2$-path
regularizer'', which appears in deep linear networks \citep{mianjy_arora_2019}
and shallow ReLU-activated networks \citep{arora_et_al_2021}. Further,
\cite{arora_et_al_2021} exhibit a data distribution-dependent regularizer in
two-layer matrix sensing/completion problems. This regularizer collapses to a
nuclear norm penalty for specific distributions.

\cite{gal_ghahramani_2016a} show that empirical risk minimization in deep neural
networks with dropout may be recast as performing Bayesian variational inference
to approximate the intractable posterior resulting from a deep Gaussian process
prior. The Bayesian viewpoint also allows for the quantification of uncertainty.
\cite{gal_ghahramani_2016b} further generalizes this technique to recurrent and
long-short-term-memory (LSTM) networks. \cite{wu_gu_2015} analyze dropout
applied to the max-pooling layers in convolutional neural networks.
\cite{wang_manning_2013} present a Gaussian approximation to the gradient noise
induced by dropout.

Generalization results for dropout training exist in various settings. Given
bounds on the norms of weight vectors, \cite{wan_et_al_2013},
\cite{gao_zhou_2016}, and \cite{zhai_wang_2018} prove decreasing Rademacher
complexity bounds as the dropout rate increases. \cite{arora_et_al_2021} bound
the Rademacher complexity of shallow ReLU-activated networks with dropout.
\cite{mcallester_2013} obtains a PAC-Bayes bound for dropout and illustrates a
trade-off between large and small dropout probabilities for different terms in
the bound.

Recently, \cite{manita_et_al_2022} demonstrated an universal approximation
result in the vein of classic results \citep{cybenko_1998, hornik_1991,
leshno_et_al_1993}, stating that any function in some generic semi-normed space
that can be $\eps$-approximated by a deterministic neural network may also be
stochastically approximated in $L^q$-norm by a sufficiently large network with
dropout.

Less is known about gradient descent training with dropout.
\cite{senen-cerda_sanders_2020} study the gradient flow associated with the
explicit regularizer obtained by marginalizing the dropout noise in a shallow
linear network. In particular, the flow converges exponentially fast within a
neighborhood of a parameter vector satisfying a balancing condition.
\cite{mianjy_arora_2020} study gradient descent with dropout on the logistic
loss of a shallow ReLU-activated network in a binary classification task. Their
main result includes an explicit rate for the misclassification error, assuming
an overparametrized network operating in the so-called lazy regime, where the
trained weights stay relatively close to their initializations, and two
well-separated classes.

\subsection{Notation}

Column vectors $\bfx = (x_1, \dots, x_d)\tran$ are denoted by bold letters. We
define $\bfzero := (0, \ldots, 0)\tran$, $\bfone := (1, \ldots, 1)\tran$, and
the Euclidean norm $\norm*{\bfx}_2 := \sqrt{\bfx\tran \bfx}$. The $d \times d$
identity matrix is symbolized by $I_d$, or simply $I$, when the dimension $d$ is
clear from context. For matrices $A,B$ of the same dimension, $A \odot B$ denotes the
Hadamard/entry-wise product $(A \odot B)_{ij} = A_{ij} B_{ij}$. We write
$\Diag(A) := I \odot A$ for the diagonal matrix with the same main diagonal as
$A$. Given $p \in (0, 1)$, we define the matrices \begin{align*}
  \ol{A} &:= A - \Diag(A)\\
  A_p &:= p A + (1 - p) \Diag(A).
\end{align*} In particular, $A_p = p \ol{A} + \Diag(A)$, so $A_p$ results from
rescaling the off-diagonal entries of $A$ by $p$.

The smallest eigenvalue of a symmetric matrix $A$ is denoted by $\lambdamin(A)$.
The operator norm of a linear operator $T : V \to W$ between normed linear
spaces is given by $\norm*{T}_{\op} := \sup_{v \in V : \norm*{v}_V \leq 1}
\norm*{T v}_W$. We write $\norm{\ \cdot\ }$ for the spectral norm of matrices,
which is the operator norm induced by $\norm{\ \cdot\ }_2$. For symmetric
matrices, the relation $A \geq B$ signifies $\bfx\tran (A - B) \bfx \geq 0$ for
all non-zero vectors $\bfx$. The strict operator inequality $A > B$ is defined
analogously.


\section{Gradient Descent and Dropout}
\label{sec::drop_variations}

\begin{figure}
  \centering
  \begin{tikzpicture}[very thick]
    \node[shape=circle, draw=Mahogany, label=below:$y$] (N) at (0,0) {};
    \node[shape=circle, draw=Mahogany, label=below:$x_1$] (I1) at (-3.5, -3)
    {};
    \node[shape=circle, draw=Mahogany, label=below:$x_2$] (I2) at (-2.5, -3)
    {};
    \node[shape=circle, draw=Mahogany, label=below:$x_{d - 1}$] (I3) at (2.5, -3)
    {};
    \node[shape=circle, draw=Mahogany, label=below:$x_d$] (I4) at (3.5, -3)
    {};

    \node (D1) at (1, -3) {};
    \node (D2) at (-1, -3) {};

    \path[dotted, draw=Mahogany] (D1) edge node {} (D2);
    \path[-, draw=Mahogany] (I1) edge node[label=left:$\beta_1$] {} (N);
    \path[-, draw=Mahogany] (I2) edge node[label=right:$\beta_2\ \cdots$] {} (N);
    \path[-, draw=Mahogany] (I3) edge node[label=left:$\beta_{d - 1}$] {} (N);
    \path[-, draw=Mahogany] (I4) edge node[label=right:$\beta_d$] {} (N);
  \end{tikzpicture}
  \caption{The linear regression model $y = \sum_{i = 1}^d \beta_i x_i$, viewed
  as a neural network without hidden layers. \label{fig::lin_network}}
\end{figure}
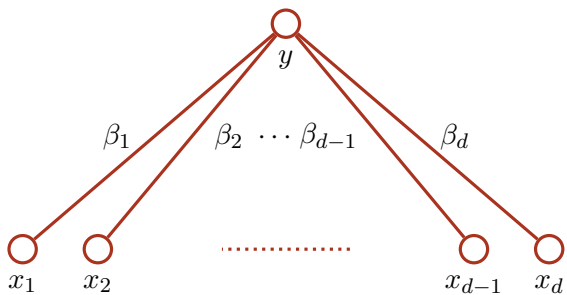

We consider a linear regression model with fixed $n \times d$ design matrix $X$
and $n$ outcomes $\bfY$, so that \begin{align}
  \label{eq::lin_reg}
  \bfY = X \starbeta + \bfeps,
\end{align} with unknown parameter $\starbeta$. We assume $\E\big[\bfeps\big] =
\bfzero$ and $\Cov\big(\bfeps\big) = I_n$. The task is to estimate $\starbeta$
from the observed data $(X, \bfY)$. As the Gram matrix $X\tran X$ appears
throughout our analysis, we introduce the shorthand \begin{equation*}
  \bbX := X\tran X.
\end{equation*} Recovery of $\starbeta$ in the linear regression model
\eqref{eq::lin_reg} may be interpreted as training a neural network without
intermediate hidden layers, see Figure \ref{fig::lin_network}. If $X$ were to
have a zero column, the corresponding regression coefficient would not affect
the response vector $\bfY$. Consequently, both the zero column and the
regression coefficient may be eliminated from the linear regression model.
Without zero columns, the model is said to be in \textit{reduced form}.

The least squares criterion for the estimation of $\starbeta$ refers to the
objective function $\bfbeta \mapsto \tfrac{1}{2} \norm*{\bfY - X \bfbeta}^2_2$.
Given a fixed learning rate $\alpha > 0$, performing gradient descent on the
least squares objective leads to the iterative scheme \begin{align}
  \wtbeta^{\mathrm{gd}}_{k + 1} = \wtbeta^{\mathrm{gd}}_k - \alpha
  \nabla_{\wtbeta^{\mathrm{gd}}_k} \dfrac{1}{2} \norm*{\bfY - X
  \wtbeta^{\mathrm{gd}}_k}^2_2 = \wtbeta^{\mathrm{gd}}_k + \alpha X\tran
  \big(\bfY -X \wtbeta^{\mathrm{gd}}_k\big)
  \label{eq::gd_no_drop}
\end{align} with $k = 0, 1, 2, \ldots$ and (possibly random) initialization
$\wtbeta^{\mathrm{gd}}_0$.

For standard gradient descent as defined in \eqref{eq::gd_no_drop}, the estimate
is updated with the gradient of the full model. Dropout, as introduced in
\cite{srivastava_et_al_2014}, replaces the gradient of the full model with the
gradient of a randomly reduced model during each iteration of training. To make
this notion more precise, we call a random diagonal $d \times d$ matrix $D$ a
$p$-\textit{dropout matrix}, or simply a \textit{dropout matrix}, if its
diagonal entries satisfy $D_{ii} \overset{\scriptscriptstyle i.i.d.}{\sim}
\Ber(p)$ for some $p \in (0, 1)$. We note that the Bernoulli distribution may
alternatively be parametrized with the failure probability $q := 1 - p$, but
following \cite{srivastava_et_al_2014} we choose the success probability $p$.

On average, $D$ has $pd$ diagonal entries equal to $1$ and $(1 - p)d$ diagonal
entries equal to $0$. Given any vector $\bfbeta$, the coordinates of $D \bfbeta$
are randomly set to $0$ with probability $1 - p$. For simplicity, the dependence
of $D$ on $p$ will be omitted.

Now, let $D_k$, $k = 1, 2, \ldots$ be a sequence of i.i.d.\ dropout matrices,
where $D_k$ refers to the dropout matrix applied in the $k$\textsuperscript{th}
iteration. Gradient descent with dropout takes the form \begin{align}
  \wtbeta_{k + 1} = \wtbeta_k - \alpha \nabla_{\wtbeta_k} \dfrac{1}{2}
  \norm*{\bfY - X D_{k + 1} \wtbeta_k}^2_2 = \wtbeta_k + \alpha D_{k + 1} X\tran
  \big(\bfY - X D_{k + 1} \wtbeta_k\big)
  \label{eq::wt_drop_def}
\end{align} with $k = 0, 1, 2, \ldots$ and (possibly random) initialization
$\wtbeta_0$. In contrast with  \eqref{eq::gd_no_drop}, the gradient in
\eqref{eq::wt_drop_def} is taken on the model reduced by the action of
multiplying $\wtbeta_k$ with $D_{k + 1}$. Alternatively, \eqref{eq::wt_drop_def}
may be interpreted as replacing the design matrix $X$ with the reduced matrix $X
D_{k + 1}$ during the $(k + 1)$\textsuperscript{th} iteration. The columns of
the reduced matrix are randomly deleted with a probability of $1 - p$. Observe
that the dropout matrix appears inside the squared norm, making the gradient
quadratic in $D_{k + 1}$.

Dropout, as defined in \eqref{eq::wt_drop_def}, considers the full gradient of
the reduced model, whereas another variant is obtained through reduction of the
full gradient. The resulting iterative scheme takes the form \begin{align}
  \whbeta_{k + 1} = \whbeta_k - \alpha D_{k + 1} \nabla_{\whbeta_k} \dfrac{1}{2}
  \norm*{\bfY - X \whbeta_k}^2_2 = \whbeta_k + \alpha D_{k + 1} X\tran \big(\bfY
  - X \whbeta_k\big)
  \label{eq::wh_drop_def}
\end{align} with $k = 0, 1, 2, \ldots$ and (possibly random) initialization
$\whbeta_0$. As opposed to $\wtbeta_k$ defined above, the dropout matrix only
occurs once in the updates, so we shall call this method \textit{simplified
dropout} from here on. As we will illustrate, the quadratic dependence of
$\wtbeta_k$ on $D_{k + 1}$ creates various challenges, whereas the analysis of
$\whbeta_k$ is more straightforward.

Both versions \eqref{eq::wt_drop_def} and \eqref{eq::wh_drop_def} coincide when
the Gram matrix $\bbX = X\tran X$ is diagonal, meaning when the columns of $X$
are orthogonal. To see this, note that diagonal matrices commute, so $D_{k +
1}^2 = D_{k + 1}$ and hence $D_{k + 1} \bbX D_{k + 1} = D_{k + 1} \bbX$.

We note that dropout need not require the complete removal of neurons. Each
neuron may be multiplied by an arbitrary (not necessarily Bernoulli distributed)
random variable. For instance, \cite{srivastava_et_al_2014} also report good
performance for $\calN(1, 1)$-distributed diagonal entries of the dropout
matrix. However, the Bernoulli variant seems well-motivated from a model
averaging perspective. \cite{srivastava_et_al_2014} propose dropout with the
explicit aim of approximating a Bayesian model averaging procedure over all
possible combinations of connections in the network. The random removal of nodes
during training is thought to prevent the neurons from co-adapting, recreating
the model averaging effect. This is the main variant implemented in popular
software libraries, such as \textit{Caffe} \citep{jia_et_al_2014},
\textit{TensorFlow} \citep{abadi_et_al_2016}, \textit{Keras}
\citep{chollet_2015}, and \textit{PyTorch} \citep{paszke_et_al_2019}.

Numerous variations and extensions of dropout exist. \cite{wan_et_al_2013} show
state-of-the-art results for networks with \textit{DropConnect}, a
generalization of dropout where connections are randomly dropped, instead of
neurons. In the linear model, this coincides with standard dropout.
\cite{ba_frey_2013} analyze the case of varying dropout probabilities, where the
dropout probability for each neuron is computed using binary belief networks
that share parameters with the underlying fully connected network. An adaptive
procedure for the choice of dropout probabilities is presented in
\cite{kingma_et_al_2015}, while also giving a Bayesian justification for
dropout. 

For a comprehensive overview of established methods and cutting-edge variants,
see \cite{moradi_et_al_2020} and \cite{santos_papa_2022}.


\section{Analysis of Averaged Dropout}
\label{sec::avg_drop}

Before presenting our main results on iterative dropout schemes, we further
discuss some properties of the marginalized loss minimizer that was first
analyzed by \cite{srivastava_et_al_2014}. For the linear regression model
\eqref{eq::lin_reg}, marginalizing the dropout noise leads to \begin{align}
  \wtbeta := \argmin_{\bfbeta} \E\Big[\norm[\big]{\bfY - X D \bfbeta}^2_2 \mid
  \bfY\Big].
  \label{eq::averaged_loss}
\end{align} One may hope that the dropout gradient descent recursion for
$\wtbeta_k$ in \eqref{eq::wt_drop_def} leads to a minimizer of
\eqref{eq::averaged_loss}, so that the marginalized loss minimizer may be
studied as a surrogate for the behaviour of $\wtbeta_k$ in the long run.

Intuitively, the gradient descent iterates with dropout represent a Monte-Carlo
estimate of some deterministic algorithm \citep{wang_manning_2013}. This can be
motivated by separating the gradient descent update into a part without
algorithmic randomness and a centered noise term, meaning \begin{equation}
  \label{eq::gd_noise_split}
  \begin{split}
    \wtbeta_{k + 1} = \wtbeta_k 
    &- \dfrac{\alpha}{2}\E\bigg[\nabla_{\wtbeta_k}
    \norm*{\bfY - X D_{k + 1} \wtbeta_k}^2_2\ \big\vert\ \bfY, \wtbeta_k\bigg]\\
    &- \dfrac{\alpha}{2} \Bigg(\nabla_{\wtbeta_k} \norm*{\bfY - X D_{k + 1}
    \wtbeta_k}^2_2 - \E\bigg[\nabla_{\wtbeta_k} \norm*{\bfY - X D_{k + 1}
    \wtbeta_k}^2_2\ \big\vert\ \bfY, \wtbeta_k\bigg]\Bigg).
  \end{split}
\end{equation} Notably, the stochastic terms form a martingale difference
sequence with respect to $(\bfY, \wtbeta_k)$. It seems conceivable that the
noise in \eqref{eq::gd_noise_split} averages out; despite the random variables
being neither independent, nor identically distributed, one may hope that a law
of large numbers still holds, see \cite{andrews_1988}. In this case, after a
sufficient number of gradient steps, \begin{equation}
  \label{eq::heuristic_av_GD}
  \begin{split}
    \wtbeta_{k + 1} &= \wtbeta_0 - \dfrac{\alpha}{2} \sum_{\ell = 1}^k
    \nabla_{\wtbeta_\ell} \norm*{\bfY - X D_{\ell + 1} \wtbeta_\ell}^2_2\\
    &\approx \wtbeta_0 - \dfrac{\alpha}{2} \sum_{\ell = 1}^k
    \E\bigg[\nabla_{\wtbeta_\ell} \norm*{\bfY - X D_{\ell + 1}
    \wtbeta_\ell}^2_2\ \big\vert\ \bfY, \wtbeta_k\bigg].
  \end{split}
\end{equation} The latter sequence could plausibly converge to the marginalized
loss minimizer $\wtbeta$. While this motivates studying $\wtbeta$, the main
conclusion of our work is that this heuristic is not entirely correct and
additional noise terms occur in the limit $k \to \infty$.

As the marginalized loss minimizer $\wtbeta$ still plays a  pivotal role in our
analysis, we briefly recount and expand on some of the properties derived in
\cite{srivastava_et_al_2014}.  Recall that $\bbX = X\tran X$, so we have
\begin{align*}
  \norm[\big]{\bfY - X D \bfbeta}^2_2 = \norm*{\bfY}^2_2 - 2 \bfY\tran X D
  \bfbeta + \bfbeta\tran D \bbX D \bfbeta.
\end{align*} Since $D$ is diagonal, $\E[D] = p I_d$, and by Lemma
\ref{lem::drop_moments_A}, $\E[D \bbX D] = p^2 \bbX + p (1 - p)\Diag(\bbX)$,
\begin{align}
  \label{eq::tikh_func_exp}
  \E\Big[\norm[\big]{\bfY - X D \bfbeta}^2_2\ \big\vert\ \bfY\Big] &=
  \norm{\bfY}_2^2-2p \bfY\tran X \bfbeta+p^2\bfbeta\tran \bbX
  \bfbeta + p (1 - p) \bfbeta\tran \Diag(\bbX) \bfbeta \notag\\
  &= \norm[\big]{\bfY - p X \bfbeta}^2_2 + p (1 - p)\bfbeta\tran \Diag(\bbX)
  \bfbeta.
\end{align} The right-hand side may be identified with a Tikhonov functional, or
an $\ell_2$-penalized least squares objective. Its gradient with respect to
$\bfbeta$ is given by \begin{align*}
  \nabla_{\bfbeta} \E\big[\norm{\bfY - X D \bfbeta}^2_2 \mid \bfY\big] = - 2 p
  X\tran \bfY + 2 \big(p^2 \bbX + p (1 - p) \Diag(\bbX)\big) \bfbeta. 
\end{align*} Recall from the discussion following Equation \eqref{eq::lin_reg}
that the model is assumed to be in reduced form, meaning $\min_i \bbX_{ii} > 0$.
In turn, \begin{align*}
  p^2 \bbX + p (1 - p) \Diag(\bbX) \geq p (1 - p) \Diag(\bbX) \geq p (1 - p)
  \min_i \bbX_{ii} \cdot I_d    
\end{align*} is bounded away from $0$, making $p^2 \bbX + p (1 - p) \Diag(\bbX)$
invertible. Solving the gradient for the minimizer $\wtbeta$ now leads to
\begin{align}
  \wtbeta = \argmin_{\bfbeta \in \R^d} \E\Big[\norm[\big]{\bfY - X D
  \bfbeta}^2_2\ \big\vert\ \bfY\Big] = p \Big(p^2 \bbX + p (1 - p)
  \Diag(\bbX)\Big)\inv X\tran \bfY = \bbX_p\inv X\tran \bfY,
  \label{eq::mlls_drop}
\end{align} where $\bbX_p := p \bbX + (1 - p) \Diag(\bbX)$. If the columns of
$X$ are orthogonal, then $\bbX$ is diagonal and hence $\bbX_p = \bbX$. In this
case, $\wtbeta$ matches the usual linear least squares estimator $\bbX\inv
X\tran \bfY$. Alternatively, $\wtbeta$ minimizing the marginalized loss can also
be deduced from the identity \begin{align}
  \label{eq::avg_obj_split}
  \E\Big[\norm[\big]{\bfY - X D \whbeta}_2^2\ \big\vert\ \bfY\Big] =
  \E\Big[\norm[\big]{\bfY - X D \wtbeta}_2^2\ \big\vert\ \bfY\Big] +
  \E\Big[\norm[\big]{X D \big(\wtbeta - \whbeta\big)}_2^2\ \big\vert\ \bfY\Big],
\end{align} which holds for all estimators $\whbeta$. See Appendix
\ref{sec::avg_drop_proof} for a proof of \eqref{eq::avg_obj_split}. We now
mention several other relevant properties of $\wtbeta$.

\textbf{Calibration:} \cite{srivastava_et_al_2014} recommend multiplying
$\wtbeta$ by $p$, which may be motivated as follows: Since $\bfY = X \starbeta +
\bfeps$, a small squared error $\norm[\big]{\bfY - p X \wt \bfbeta}^2_2$ in
\eqref{eq::tikh_func_exp} implies $\starbeta \approx p\wt \bfbeta$. Moreover,
multiplying $\wtbeta$ by $p$ leads to $p\wtbeta = \big(\bbX + (1/p - 1)
\Diag(\bbX)\big)\inv X\tran \bfY$ which may be identified with the minimizer of
the objective function \begin{align*}
  \bfbeta \mapsto \norm[\big]{\bfY - X\bfbeta}^2_2 + (p\inv - 1) \bfbeta\tran
  \Diag(\bbX) \bfbeta = \E\big[\norm{\bfY - X p\inv D \bfbeta}^2_2 \mid
  \bfY\big].
\end{align*} This recasts $p \wtbeta$ as resulting from a weighted form of ridge
regression. Comparing the objective function to the original marginalized loss
$\E\big[\norm{\bfY - X D \bfbeta}^2_2 \mid \bfY\big]$, the rescaling replaces
$D$ with the normalized dropout matrix $p\inv D$, which has the identity matrix
as its expected value. In popular machine learning software, the sampled dropout
matrices are usually rescaled by $p\inv$ \citep{abadi_et_al_2016, chollet_2015,
jia_et_al_2014, paszke_et_al_2019}.

In some settings multiplication by $p$ may worsen $\wtbeta$ as a statistical
estimator. As an example, consider the case $n = d$ with $X = n I_n$ a multiple
of the identity matrix. Now, \eqref{eq::mlls_drop} turns into $\wtbeta = n\inv
\bfY = \starbeta + n\inv \bfeps$. If the noise vector $\bfeps$ consists of
independent standard normal random variables, then $\wtbeta$ has mean squared
error $\E\big[\norm{\wtbeta - \starbeta}_2^2\big] = n\inv$. In contrast,
$\E\big[\norm{p \wtbeta - \starbeta}_2^2\big] = (1 - p) \norm{\starbeta}_2^2 +
p^2 n\inv$, so $\wtbeta$ converges to $\starbeta$ at the rate $n\inv$ while $p
\wtbeta$ cannot be consistent as $n\to \infty$, unless $\starbeta = \bfzero$.

The correct rescaling may also depend on the parameter dimension $d$ and the
spectrum of $\bbX$. Suppose that all columns of $X$ have the same Euclidean
norm, so that $\Diag(\bbX) = \bbX_{11} \cdot I_d$. Let $X = \sum_{\ell =
1}^{\mathrm{rank}(X)} \sigma_\ell \bfv_\ell \bfw_\ell\tran$ denote a singular
value decomposition of $X$, with singular values $\sigma_1 , \ldots,
\sigma_{\mathrm{rank}(X)}$. Now, $\wtbeta$ satisfies \begin{equation}
  \label{eq::wtbeta_ones}
  \begin{split}
    \wtbeta &= \sum_{\ell = 1}^{\mathrm{rank}(X)} \dfrac{1}{p \sigma_\ell^2 + (1
    - p) \bbX_{11}} \big(\sigma_\ell \bfw_\ell \bfv_\ell\tran\big) \bfY\\
    \E\big[X \wtbeta\big] &= \sum_{\ell = 1}^{\mathrm{rank}(X)} \dfrac{1}{p + (1
    - p) \bbX_{11}/\sigma_\ell^2} \big(\sigma_\ell \bfv_\ell \bfw_\ell\tran\big)
    \starbeta.
  \end{split}
\end{equation} For a proof of these identities, see Appendix
\ref{sec::avg_drop_proof}. To get an unbiased estimator for $X \starbeta =
\sum_{\ell = 1}^{\mathrm{rank}(X)} \sigma_\ell \bfv_\ell \bfw_\ell\tran
\starbeta$, we must undo the effect of the spectral multipliers $1/(p + (1 - p)
\bbX_{11}/\sigma_\ell^2)$ which take values in the interval $[0,1/p]$.
Consequently, the proper rescaling depends on the eigenspace. Multiplication of
the estimator by $p$ addresses the case where the singular values $\sigma_\ell$
are large. In particular, if $X = \sigma \bfv \bfw\tran$ with $\sigma =
\sqrt{nd}$, $\bfv = (n^{-1/2}, \ldots, n^{-1/2})\tran$, and $\bfw = (d^{-1/2},
\ldots, d^{-1/2})\tran$, then $X$ is the $n \times d$ matrix with all entries
equaling $1$. Now $\bbX_{11} = n$ and so $\bbX_{11} / \sigma^2 = d\inv$, meaning
the correct scaling factor depends explicitly on the parameter dimension $d$ and
converges to the dropout probability $p$ as $d \to \infty$.

\textbf{Invariance properties:} The minimizer $\wtbeta= \bbX_p\inv X\tran \bfY$
is scale invariant in the sense that $\bfY$ and $X$ may be replaced with $\gamma
\bfY$ and $\gamma X$ for some arbitrary $\gamma \neq 0$, without changing
$\wtbeta$. This does not hold for the gradient descent iterates
\eqref{eq::wt_drop_def} and \eqref{eq::wh_drop_def}, since rescaling by $\gamma$
changes the learning rate from $\alpha$ to $\alpha \gamma^2$.  Moreover,
$\wtbeta$ as well as the gradient descent iterates $\wtbeta_k$ in
\eqref{eq::wt_drop_def} and $\whbeta_k$ in \eqref{eq::wh_drop_def} are invariant
under replacement of $\bfY$ and $X$ by $Q \bfY$ and $Q X$ for any orthogonal $n
\times n$ matrix $Q$. See \cite{helmbold_long_2017} for further results on
scale-invariance of dropout in deep networks.

\textbf{Overparametrization:} Dropout has been successfully applied in the
overparametrized regime, see for example \cite{krizhevsky_et_al_2012}. For the
overparametrized linear regression model, the data-misfit term in
\eqref{eq::tikh_func_exp} suggests that $p X \wtbeta = X \big(\bbX + (p\inv - 1)
\Diag(\bbX)\big)\inv X\tran \bfY$ should be close to the data vector $\bfY$.
However, \begin{equation}
  \label{eq::data_miss_norm}
  \norm[\Big]{X \big(\bbX + (p\inv - 1) \Diag(\bbX)\big)\inv X\tran} < 1.
\end{equation} See Appendix \ref{sec::avg_drop_proof} for a proof. Hence, $p X
\wtbeta$ also shrinks $\bfY$ towards zero in the overparametrized regime and
does not interpolate the data. The variational formulation $\wtbeta \in
\argmin_{\bfbeta \in \R^d} \norm[\big]{\bfY - p X \bfbeta}^2_2 + p (1 - p)
\bfbeta\tran \Diag(\bbX) \bfbeta$ reveals that $\wtbeta$ is a minimum-norm
solution in the sense \begin{align*}
  \wtbeta \in \argmin_{\bfbeta : X \bfbeta = X\wtbeta}\ \bfbeta\tran \Diag(\bbX)
  \bfbeta,
\end{align*} which explains the induced shrinkage.


\section{Analysis of Iterative Dropout Schemes}
\label{sec::gd_results}

In the linear model, gradient descent with a small but fixed learning rate, as
in \eqref{eq::gd_no_drop}, leads to exponential convergence in the number of
iterations. Accordingly, we analyze the iterative dropout schemes
\eqref{eq::wt_drop_def} and \eqref{eq::wh_drop_def} for fixed learning rate
$\alpha$ and only briefly discuss the algebraically less tractable case of
decaying learning rates.

\subsection{Convergence of Dropout}
\label{sec::conv_wtbeta}

We proceed by assessing convergence of the iterative dropout scheme
\eqref{eq::wt_drop_def}, as well as some of its statistical properties. Recall
that gradient descent with dropout takes the form \begin{align}
  \label{eq::wt_drop_def_recall}
  \wtbeta_k = \wtbeta_{k - 1} + \alpha D_k X\tran \big(\bfY - X D_k \wtbeta_{k -
  1}\big) = \big(I - \alpha D_k \bbX D_k\big) \wtbeta_{k - 1} + \alpha D_k
  X\tran \bfY.
\end{align} As alluded to in the beginning of Section \ref{sec::avg_drop}, the
gradient descent iterates should be related to the minimizer $\wtbeta$ of
\eqref{eq::averaged_loss}. It then seems natural to study the difference
$\wtbeta_k - \wtbeta$, with $\wtbeta$ as an ``anchoring point''. Comparing
$\wtbeta_k$ and $\wtbeta$ demands an explicit analysis, without marginalization
of the dropout noise.

To start, we rewrite the updating formula \eqref{eq::wt_drop_def_recall} in
terms of $\wtbeta_k - \wtbeta$. Using $D_k^2 = D_k$, $\Diag(\bbX) = \bbX_p - p
\ol{\bbX}$, and that diagonal matrices always commute, we obtain $D_k \bbX D_k =
D_k \ol{\bbX} D_k + D_k \Diag(\bbX) = D_k \ol{\bbX} D_k + D_k\bbX_p - p D_k
\ol\bbX$. As defined in \eqref{eq::mlls_drop}, $\bbX_p\wtbeta = X\tran \bfY$ and
thus \begin{align}
  \label{eq::wtbeta_rewrite}
  \begin{split}
    \wtbeta_k - \wtbeta &= \big(I - \alpha D_k \bbX D_k\big) \big(\wtbeta_{k -
    1} - \wtbeta\big) + \alpha D_k \ol{\bbX} \big(p I - D_k\big) \wtbeta\\
    &= \big(I - \alpha D_k \bbX_p\big) \big(\wtbeta_{k - 1} - \wtbeta\big) +
    \alpha D_k \ol{\bbX} \big(p I - D_k\big) \wtbeta_{k - 1}.
  \end{split}
\end{align} 

In both representations, the second term is centered and uncorrelated for
different values of $k$. Vanishing of the mean follows from the independence of
$D_k$ and $\big(\wtbeta, \wtbeta_{k - 1}\big)$, combined with $\E\big[D_k
\ol{\bbX} (p I - D_k)\big] = 0$, the latter being shown in
\eqref{eq::exp_ol_rewrite}. If $k > \ell$, independence of $D_k$ and
$\big(\wtbeta, \wtbeta_{k - 1}, \wtbeta_{\ell - 1}\big)$, as well as $\E\big[D_k
\ol{\bbX} (p I - D_k)\big] = 0$, imply $\Cov\big(D_k \ol{\bbX} (p I - D_k)
\wtbeta, D_\ell \ol{\bbX} (p I - D_\ell) \wtbeta\big) = 0$ and $\Cov\big(D_k
\ol{\bbX} (p I - D_k) \wtbeta_{k-1}, D_\ell \ol{\bbX} (p I - D_\ell)
\wtbeta_{\ell-1}\big) = 0$, which proves uncorrelatedness.

Defining $Z_k := \wtbeta_k - \wtbeta$, $G_k := I - \alpha D_k \bbX D_k$, and
$\bfxi_k := \alpha D_k \ol{\bbX} \big(p I - D_k\big) \wtbeta$, the first
representation in \eqref{eq::wtbeta_rewrite} may be identified with a lag one
vector autoregressive (VAR) process \begin{align}
  Z_k = G_k Z_{k - 1} + \bfxi_k
  \label{eq::RAR}
\end{align} with i.i.d.\ random coefficients $G_k$ and noise/innovation process
$\bfxi_k$. As just shown, $\E[\bfxi_k] = 0$ and $\Cov(\bfxi_k, \bfxi_\ell) = 0$
whenever $k\neq \ell$, so the noise process is centered and serially
uncorrelated. The random coefficients $G_k$ and $\bfxi_k$ are, however,
dependent. While most authors do not allow for random coefficients $G_k$ in VAR
processes, such processes are special cases of a random autoregressive process
(RAR) \citep{regis_et_al_2022}.

In the VAR literature, identifiability and estimation of the random coefficients
$G_k$ is considered \citep{nicholls_quinn_1982, regis_et_al_2022}. In contrast,
we aim to obtain bounds for the convergence of $\E\big[\wtbeta_k - \wtbeta\big]$
and $\Cov\big(\wtbeta_k - \wtbeta\big)$. Difficulties arise from the involved
structure and coupled randomness of $G_k$ and $\bfxi_k$. Estimation of
coefficients under dependence of $G_k$ and $\bfxi_k$ is treated in
\cite{hill_peng_2014}.

For a sufficiently small learning rate $\alpha$, the random matrices $I - \alpha
D_k \bbX D_k$ and $I - \alpha D_k \bbX_p$ in both representations in
\eqref{eq::wtbeta_rewrite} are contractive maps in expectation. By Lemma
\ref{lem::drop_moments_A}, their expected values coincide since \begin{align*}
  \E\big[I - \alpha D_k \bbX D_k\big] = \E\big[I - \alpha D_k \bbX_p\big] = I -
  \alpha p \bbX_p.
\end{align*} For the subsequent analysis, we impose the following mild
conditions that, among other things, establish contractivity of $I - \alpha p
\bbX_p$ as a linear map.

\begin{assumption}
  \label{ass::sec_wtbeta_all}
  The learning rate $\alpha$ and the dropout probability $p$ are chosen such
  that $\alpha p \norm*{\bbX} < 1$, the initialization $\wtbeta_0$ is a square
  integrable random vector that is independent of the data $\bfY$ and the model
  is in reduced form, meaning that $X$ does not have zero columns. 
\end{assumption}

For gradient descent without dropout and fixed learning rate, as defined in
\eqref{eq::gd_no_drop}, $\alpha \norm{\bbX} < 1$ guarantees converge of the
scheme in expectation. We will see shortly that dropout essentially replaces the
expected learning rate with $\alpha p$, which motivates the condition $\alpha p
\norm{\bbX} < 1$.

As a straightforward consequence of the definitions, we are now able to show
that $\wtbeta_k - \wtbeta$ vanishes in expectation at a geometric rate. For a
proof of this as well as subsequent results, see Appendix
\ref{sec::gd_results_proof}.

\begin{lemma}[Convergence of Expectation]
  \label{lem::wt_exp_conv}
  Given Assumption \ref{ass::sec_wtbeta_all}, $\norm[\big]{I - \alpha p \bbX_p}
  \leq 1 - \alpha p (1 - p) \min_{i} \bbX_{ii} < 1$ and for any $k = 0, 1,
  \ldots$, \begin{align*}
    \norm[\Big]{\E\big[\wtbeta_k - \wtbeta\big]}_2 \leq \norm[\big]{I -
    \alpha p \bbX_p}^k \norm[\Big]{\E\big[\wtbeta_0 - \wtbeta\big]}_2.
  \end{align*}
\end{lemma}

Before turning to the analysis of the covariance structure, we highlight a
property of the sequence $\E\big[\wtbeta_k \mid \bfY\big]$. As mentioned, these
conditional expectations may be viewed as gradient descent iterates generated by
the marginalized objective $\tfrac{1}{2} \E\big[\norm{\bfY - X D \bfbeta}_2^2
\mid \bfY\big]$ that gives rise to $\wtbeta$. Indeed, combining
\eqref{eq::wt_drop_def_recall} with $\E[D_{k + 1}] = p I_d$, $\E\big[D_{k + 1}
\bbX D_{k + 1}\big] = p\bbX_p$ from Lemma \ref{lem::drop_moments_A}, and
\eqref{eq::tikh_func_exp} yields \begin{align*}
  \E\big[\wtbeta_{k + 1} \mid \bfY\big] &= \E\big[\wtbeta_k \mid \bfY\big] +
  \alpha p X\tran \bfY - \alpha p \bbX_p \E\big[\wtbeta_k \mid \bfY\big]\\
  &= \E\big[\wtbeta_k \mid \bfY\big] - \dfrac{\alpha}{2} \nabla_{\E[\wtbeta_k
  \mid \bfY]} \E\Big[\norm[\big]{\bfY - X D \E[\wtbeta_k \mid \bfY]}_2^2 \mid
  \bfY\Big].
\end{align*} This establishes a connection between the dropout iterates and the
averaged analysis of the previous section. However, the relationship between the
(unconditional) covariance matrices $\Cov(\wtbeta_k)$ and the added noise
remains unclear. A new dropout matrix is sampled for each iteration, whereas
$\wtbeta$ results from minimization only after applying the conditional
expectation $\E[\ \cdot \mid \bfY]$ to the randomized objective function. Hence,
we may expect that $\wtbeta$ features smaller variance than the iterates as the
latter also depend on the noise added via dropout.

As a first result for the covariance analysis, we establish an extension of the
Gauss-Markov theorem stating that the covariance matrix of a linear estimator
lower-bounds the covariance matrix of an affine estimator, provided that both
estimators have the same asymptotic mean. Moreover, the covariance matrix of
their difference characterizes the gap. We believe that a similar result may
already be known, but we are not aware of any reference, so a full proof is
provided in Appendix \ref{sec::gd_results_proof} for completeness.

\begin{theorem}
  \label{thm::gm}
   In the linear regression model \eqref{eq::lin_reg}, consider estimators
   $\wtbeta_A = A X\tran \bfY$ and $\wtbetaaff = B \bfY + \bfa$, with $B \in
   \R^{d \times n}$ and $\bfa \in \R^d$ (possibly) random, but independent of
   $\bfY$, and $A \in \R^{d \times d}$ deterministic. Then, \begin{align*}
    \norm[\Big]{\Cov\big(\wtbetaaff\big) - \Cov\big(\wtbeta_A\big) -
    \Cov\big(\wtbetaaff - \wtbeta_A\big)}
    \leq 4 \norm{A} \sup_{\starbeta : \norm{\starbeta}_2 \leq 1}
    \norm[\big]{\E_{\starbeta}\big[\wtbetaaff - \wtbeta_A\big]}_2,
  \end{align*} where $\E_{\starbeta}$ denotes the expectation with respect to
  $\starbeta$ being the true regression vector in the linear regression model
  \eqref{eq::lin_reg}.
\end{theorem}

Since $\Cov\big(\wtbetaaff\big) - \Cov\big(\wtbeta_A\big) - \Cov\big(\wtbetaaff
- \wtbeta_A\big) = \Cov\big(\wtbetaaff - \wtbeta_A, \wtbeta_A\big) +
\Cov\big(\wtbeta_A, \wtbetaaff - \wtbeta_A\big)$, Theorem \ref{thm::gm} may be
interpreted as follows: if the estimators $\wtbetaaff$ and $\wtbeta_A$ are
nearly the same in expectation, then $\wtbetaaff - \wtbeta_A$ and $\wtbeta_A$
must be nearly uncorrelated. In turn, $\wtbeta_k$ may be decomposed into
$\wtbeta_A$ and (nearly) orthogonal noise $\wtbetaaff - \wtbeta_A$, so that
$\Cov\big(\wtbetaaff\big) \approx \Cov\big(\wtbeta_A\big) + \Cov\big(\wtbetaaff
- \wtbeta_A\big)$ is lower bounded by $\Cov\big(\wtbeta_A\big)$. Therefore, the
covariance matrix $\Cov\big(\wtbetaaff - \wtbeta_A\big)$ quantifies the gap in
the bound.

Taking $A := \bbX\inv$ and considering linear estimators with $\bfa = \bfzero$
recovers the usual Gauss-Markov theorem, stating that $\bbX\inv X\tran \bfY$ is
the best linear unbiased estimator (BLUE) for the linear model. Applying the
generalized Gauss-Markov theorem with $A = (\bbX + \Gamma)\inv$, where $\Gamma$
is a positive definite matrix, we obtain the following statement about
$\ell_2$-penalized estimators.

\begin{corollary}
  \label{cor::GM_l2}
  The minimizer $\wtbeta_\Gamma := (\bbX + \Gamma)\inv X\tran \bfY$ of the
  $\ell_2$-penalized functional $\norm[\big]{\bfY - X \bfbeta}_2^2 +
  \bfbeta\tran \Gamma \bfbeta$ has the smallest covariance matrix among all
  affine estimators with the same expectation as $\wtbeta_\Gamma$.
\end{corollary}

We now return to our analysis of the covariance structure induced by dropout. If
$A := \bbX_p\inv$, then $\wtbeta = \bbX_p\inv X\tran \bfY = A X\tran \bfY =
\wtbeta_A$ in Theorem \ref{thm::gm}. Further, the dropout iterates may be
rewritten as affine estimators $\wtbeta_k = B_k \bfY + \bfa_k$ with
\begin{equation*}
  \begin{split}
    B_k &:= \sum_{j = 1}^{k - 1} \left(\prod_{\ell = 0}^{k - j - 1} \Big(I -
    \alpha D_{k - \ell} \bbX D_{k - \ell}\Big)\right) \alpha D_j X\tran + \alpha
    D_k X\tran\\
    \bfa_k &:= \left(\prod_{\ell = 0}^{k - 1} \Big(I - \alpha D_{k - \ell} \bbX
    D_{k - \ell}\Big)\right) \wtbeta_0.
  \end{split}
\end{equation*} By construction, $(B_k, \bfa_k)$ and $\bfY$ are independent, so
Theorem \ref{thm::gm} applies. As shown in Lemma \ref{lem::wt_exp_conv},
$\E\big[\wtbeta_k - \wtbeta\big]$ vanishes exponentially fast, so we conclude
that $\Cov\big(\wtbeta_k\big)$ is asymptotically lower-bounded by
$\Cov\big(\wtbeta\big)$. Further, the covariance structure of $\wtbeta$ is
optimal in the sense of Corollary \ref{cor::GM_l2}.

We proceed by studying $\Cov(\wtbeta_k - \wtbeta)$, with the aim of quantifying
the gap between the covariance matrices. To this end, we exhibit a particular
recurrence for the second moments $\E\big[(\wtbeta_k - \wtbeta) (\wtbeta_k -
\wtbeta)\tran\big]$. Recall that $\odot$ denotes the Hadamard product, $B_p = p
B + (1 - p) \Diag(B)$, and $\ol{B} = B - \Diag(B)$.

\begin{lemma}[Second Moment - Recursive Formula]
  \label{lem::wt_var_rec}
  Under Assumption \ref{ass::sec_wtbeta_all}, for all positive integers $k$
  \begin{align*}
    \begin{split}
      &\norm[\Bigg]{\E\Big[\big(\wtbeta_k - \wtbeta\big) \big(\wtbeta_k -
      \wtbeta\big)\tran\Big]- S\bigg(\E\Big[\big(\wtbeta_{k - 1} - \wtbeta\big)
      \big(\wtbeta_{k - 1} - \wtbeta\big)\tran\Big]\bigg)}\\
      &\leq 6\norm[\big]{I - \alpha p \bbX_p}^{k - 1}
      \norm[\Big]{\E\big[(\wtbeta_0 - \wtbeta\big) \wtbeta\tran\big]},
    \end{split}
  \end{align*} where $S : \R^{d \times d} \to \R^{d \times d}$ denotes the
  affine operator \begin{align*}
    S(A) &= \big(I - \alpha p \bbX_p\big) A \big(I - \alpha p \bbX_p\big)\\
    &\quad + \alpha^2 p (1 - p) \Diag\big(\bbX_p A \bbX_p\big) + \alpha^2 p^2 (1
    - p)^2 \bbX \odot \ol{A + \E\big[\wtbeta \wtbeta\tran\big]} \odot \bbX\\
    &\quad + \alpha^2 p^2 (1 - p) \bigg(\Big(\ol{\bbX} \Diag\Big(A +
    \E\big[\wtbeta \wtbeta\tran\big]\Big) \ol{\bbX}\Big)_p + \ol{\bbX}
    \Diag\big(\bbX_p A\big) + \Diag\big(\bbX_p A\big) \ol{\bbX}\bigg).
  \end{align*}
\end{lemma}

Intuitively, the lemma states that the second moment of $\wtbeta_k - \wtbeta$
evolves as an affine dynamical system, up to some exponentially decaying
remainder. This may be associated with the implicit regularization of the
dropout noise, as illustrated empirically in \cite{wei_et_al_2020}.

Mathematically, the result may be motivated via the representation of the
dropout iterates as a random autoregressive process $Z_k = G_k Z_{k - 1} +
\bfxi_k$ in \eqref{eq::RAR}. Writing out $Z_k Z_k\tran = G_k Z_{k - 1} Z_{k -
1}\tran G_k + \bfxi_k \bfxi_k\tran + G_k Z_{k - 1} \bfxi_k\tran + \bfxi_k Z_{k -
1}\tran G_k$ and comparing with the proof of the lemma, we see that the
remainder term, denoted by $\rho_{k - 1}$ in the proof, coincides with the
expected value of the cross terms $G_k Z_{k - 1} \bfxi_k\tran + \bfxi_k Z_{k -
1}\tran G_k$. Moreover, the operator $S$ is obtained by computing \begin{align*}
  S(A) = \E\Big[G_k A G_k + \bfxi_k \bfxi_k\tran\Big].
\end{align*} As $(G_k, \bfxi_k)$ are i.i.d., $S$ does not depend on $k$.
Moreover, independence of $G_k$ and $Z_{k - 1}$ implies \begin{align*}
  \E\Big[Z_k Z_k\tran\Big] &= \E\Big[G_k Z_{k - 1} Z_{k - 1}\tran G_k + \bfxi_k
  \bfxi_k\tran\Big] + \rho_{k - 1}\\
  &= \E\Big[G_k \E\big[Z_{k - 1} Z_{k - 1}\tran\big] G_k + \bfxi_k
  \bfxi_k\tran\Big] + \rho_{k - 1}\\
  &= S\Big(\E\big[Z_{k - 1} Z_{k - 1}\tran\big]\Big) + \rho_{k - 1}.
\end{align*} Inserting the definition $Z_k = \wtbeta_k - \wtbeta$ results in the
statement of the lemma. The random vector $\bfxi_k$ depends on $\wtbeta$ and by
Theorem \ref{thm::gm} the correlation between $Z_k = \wtbeta_k - \wtbeta$ and
$\wtbeta$ decreases as $k\to \infty$. This leads to the exponentially decaying
bound for the remainder term $\rho_{k - 1}$.

The previous lemma entails equality between $\E\big[Z_k Z_k\tran\big]$ and
$S^k\big(\E[Z_0 Z_0\tran]\big)$, up to the remainder terms. The latter may be
computed further by decomposing the affine operator $S$ into its intercept and
linear part \begin{align}
  \label{eq::S0_Slin_def}
  S_0 := S(0) = \E\big[\bfxi_k \bfxi_k\tran\big] \qquad \mbox{and} \qquad
  S\lin(A) := S(A) - S_0 = \E\big[G_k  A G_k\big].
\end{align} If $S\lin$ were to have operator norm less than one, then the
Neumann series for $(\id - S\lin)\inv$ (see Lemma \ref{lem::neumann}) gives
\begin{align*}
  S^k(A) = \sum_{j = 0}^{k - 1} S\lin^j(S_0) + S\lin^k(A) \to \sum_{j =
  0}^\infty S\lin^j(S_0) = \big(\id - S\lin\big)\inv S_0,
\end{align*} with $\id$ the identity operator on $d \times d$ matrices.
Surprisingly, the operator ``forgets'' $A$ in the sense that the limit does not
depend on $A$ anymore. The argument shows that $\E\big[Z_k Z_k\tran\big] =
\E\big[(\wtbeta_k - \wtbeta) (\wtbeta_k - \wtbeta)\tran\big]$ should behave like
$\big(\id - S\lin\big)\inv S_0$ in the first order. The next result makes this
precise, taking into account the remainder terms and approximation errors.

\begin{theorem}[Second Moment - Limit Formula]
  \label{thm::var_lim}
  In addition to Assumption \ref{ass::sec_wtbeta_all} suppose $\alpha <
  \tfrac{\lambdamin(\bbX_p)}{3\norm{\bbX}^2}$, then, for any $k = 1, 2, \ldots$
  \begin{align*}
    \norm[\bigg]{\E\Big[\big(\wtbeta_k - \wtbeta\big) \big(\wtbeta_k -
    \wtbeta\big)\tran\Big] - \big(\id - S\lin\big)\inv S_0} \leq C k
    \norm[\big]{I - \alpha p \bbX_p}^{k - 1}
  \end{align*} and \begin{align*}
    \norm[\Big]{\Cov\big(\wtbeta_k-\wtbeta\big) - \big(\id - S\lin\big)\inv S_0}
    \leq C k \norm[\big]{I - \alpha p \bbX_p}^{k - 1}
  \end{align*} with constant $C$ given by
  \begin{align*}
    C := \norm[\bigg]{\E\Big[\big(\wtbeta_0 - \wtbeta\big) \big(\wtbeta_0 -
    \wtbeta\big)\tran\Big] - \big(\id - S\lin\big)\inv S_0} +
    6\norm[\Big]{\E\big[(\wtbeta_0 - \wtbeta) \wtbeta\tran\big]} +
    \norm[\Big]{\E\big[\wtbeta_0 - \wtbeta\big]}_2^2.
  \end{align*}
\end{theorem}

In short, $\Cov\big(\wtbeta_k - \wtbeta)$ converges exponentially fast to the
limit $(\id - S\lin)\inv S_0$. Combining the generalized Gauss-Markov Theorem
\ref{thm::gm} with Theorem \ref{thm::var_lim} also establishes \begin{align*}
  \Cov\big(\wtbeta_k\big) \to \Cov\big(\wtbeta\big) + (\id - S\lin)\inv S_0,
  \qquad \text{as} \ k \to \infty,
\end{align*} with exponential rate of convergence. Recall the intuition gained
from Theorem \ref{thm::gm} that $\wtbeta_k$ may be decomposed into a sum of
$\wtbeta$ and (approximately) orthogonal centered noise. We now conclude that up
to exponentially decaying terms, the covariance matrix of this orthogonal noise
must be given by $(\id - S\lin)\inv S_0$, which fully describes the (asymptotic)
gap between the covariance matrices of $\wtbeta$ and $\wtbeta_k$.

Taking the trace and noting that $\abs[\big]{\Tr(A)} \leq d \norm{A}$, we obtain
a bound for the convergence of $\wtbeta_k$ with respect to the squared Euclidean
loss, \begin{align}
  \abs[\bigg]{\E\Big[\norm[\big]{\wtbeta_k - \wtbeta}_2^2\Big] -
  \Tr\Big(\big(\id - S\lin\big)\inv S_0\Big)} \leq C d k \norm[\big]{I - \alpha
  p \bbX_p}^{k - 1}.
  \label{eq::MSE_bd}
\end{align} 

Since $(\id - S\lin)\inv S_0$ is a $d\times d$ matrix, the term $\Tr\big((\id -
S\lin)\inv S_0\big)$ describing the asymptotic discrepancy between $\wtbeta_k$
and $\wtbeta$ can be large in high dimensions $d$, even if the spectral norm of
$(\id - S\lin)\inv S_0$ is small. Since $\id - S\lin$ is a positive definite
operator, the matrix $(\id - S\lin)\inv S_0$ is zero if, and only if, $S_0$ is
zero. By \eqref{eq::S0_Slin_def}, $S_0 = \E[\bfxi_k \bfxi_k\tran]$. The explicit
form $\bfxi_k = \alpha D_k \ol{\bbX} \big(p I - D_k\big) \wtbeta$, shows that
$\bfxi_k = 0$ and $S_0 = 0$ provided that $\ol{\bbX} = 0$, meaning whenever
$\bbX$ is diagonal. To give a more precise quantification, we show that the
operator norm of $(\id - S\lin)\inv S_0$ is of order $\alpha p / (1 - p)^2$.

\begin{lemma}
  \label{lem::small_alpha_p}
  In addition to Assumption \ref{ass::sec_wtbeta_all} suppose $\alpha <
  \tfrac{\lambdamin(\bbX_p)}{3\norm{\bbX}^2}$, then, for any $k = 1, 2,\ldots$
  \begin{align*}
    \norm[\Big]{\Cov\big(\wtbeta_k\big) - \Cov\big(\wtbeta\big)} \leq \dfrac{k
    \norm{I - \alpha p \bbX_p}^{k - 1} C' + \alpha p C''}{(1 - p)^2}
  \end{align*} and \begin{align*}
    \norm[\Big]{\Cov\big(\wtbeta_k\big) - \Diag(\bbX)\inv \bbX
    \Diag(\bbX)\inv} \leq \dfrac{k \norm{I - \alpha p \bbX_p}^{k - 1} C' +
    p (1 + \alpha) C''}{(1 - p)^2},
  \end{align*} where $C'$ and $C''$ are constants that are independent of
  $(\alpha, p, k)$. 
\end{lemma}

The first bound describes the interplay between $\alpha p$ and $k$. Making
$\alpha p$ small will decrease the second term in the bound, but conversely
requires a larger number of iterations $k$ for the first term to decay.

In the second bound, $\Diag(\bbX)\inv \bbX \Diag(\bbX)\inv$ matches the
covariance matrix $\Cov\big(\wtbeta\big)$ of the marginalized loss minimizer
$\wtbeta$ up to a term of order $p$. Consequently, the covariance structures
induced by dropout and $\ell_2$-regularization approximately coincide for
sufficiently small $p$. However, in this regime we have $\bbX_p = p \bbX + (1 -
p) \Diag(\bbX) \approx \Diag(\bbX)$, and $\wtbeta = \bbX_p\inv X\tran \bfY
\approx \Diag(\bbX)\inv X\tran \bfY$ becomes extremely biased whenever the Gram
matrix $\bbX$ is not diagonal.

Theorem \ref{thm::gm} already establishes $\Cov(\wtbeta)$ as the optimal
covariance among all affine estimators that are asymptotically unbiased for
$\wtbeta$. To conclude our study of the gap between $\Cov(\wtbeta_k)$ and
$\Cov(\wtbeta)$, we provide a lower-bound.

\begin{theorem}[Sub-Optimality of Variance]
  \label{thm::wtbeta_var_sub}
  In addition to the assumptions of Theorem \ref{thm::var_lim}, suppose for
  every $\ell = 1, \ldots, d$ there exists $m \neq \ell$ such that $\bbX_{\ell
  m} \neq 0$, then \begin{align*}
    \lim_{k \to \infty} \Cov\big(\wtbeta_k\big) - \Cov\big(\wtbeta\big) \geq
    \dfrac{\alpha p (1 - p)^2 \lambdamin(\bbX)}{2 \norm{\bbX}^3} \min_{i \neq j
    : \bbX_{ij} \neq 0} \bbX_{ij}^2 \cdot I_d.
  \end{align*}
\end{theorem}

The lower-bound is positive whenever $\bbX$ is invertible. In general, Theorem
\ref{thm::wtbeta_var_sub} entails asymptotic statistical sub-optimality of the
gradient descent iterates $\wtbeta_k$ for a large class of design matrices.
Moreover, the result does not require any further assumptions on the tuning
parameters $\alpha$ and $p$, other than $\alpha$ being sufficiently small.

To summarize, compared with the marginalized loss minimizer $\wtbeta$, the
covariance matrix of the gradient descent iterates with dropout may be larger.
The difference may be significant, especially if the data dimension $d$ is
large. Proving the results requires a refined second moment analysis, based on
explicit computation of the dynamics of $\wtbeta_k - \wtbeta$. Simple heuristics
such as \eqref{eq::heuristic_av_GD} do not fully reveal the properties of the
underlying dynamics.

\subsubsection{Ruppert-Polyak Averaging}

To reduce the gradient noise induced by dropout, one may consider the running
average over the gradient descent iterates. This technique is also known as
Ruppert-Polyak averaging \citep{ruppert_1988, polyak_1990} The
$k$\textsuperscript{th} Ruppert-Polyak average of the gradient descent iterates
is given by \begin{align*}
  \rpbeta_k := \dfrac{1}{k} \sum_{\ell = 1}^k \wtbeta_\ell.
\end{align*} Averages of this type are well-studied in the stochastic
approximation literature, see \cite{polyak_juditsky_1992, gyorfi_walk_1996} for
results on linear regression and \cite{zhu_et_al_2021, dereich_kassing_2023} for
stochastic gradient descent. The next theorem illustrates convergence of
$\rpbeta_k$ towards $\wtbeta$.

\begin{theorem}
  \label{thm::rp_l2_conv}
  In addition to Assumption \ref{ass::sec_wtbeta_all}, suppose $\alpha <
  \tfrac{\lambdamin(\bbX_p)}{3\norm{\bbX}^2}$, then, for any $k = 1, 2, \ldots$
  \begin{align*}
    \norm[\bigg]{\E\Big[\big(\rpbeta_k - \wtbeta\big) \big(\rpbeta_k -
    \wtbeta\big)\tran\Big]} \leq \dfrac{2 \norm{\bbX}^2 \cdot
    \norm[\big]{\E[X\tran \bfY \bfY\tran X]}}{k (1 - p) (\min_i \bbX_{ii})^4} +
    \dfrac{2 C}{k^2 (\alpha p (1 - p) \min_i\bbX_{ii})^3},
  \end{align*} where $C$ is the constant from Theorem \ref{thm::var_lim}.
\end{theorem}

The first term in the upper bound is independent of $\alpha$ and decays at the
rate $k\inv$, whereas the second term scales with $(\alpha p)^{- 3} k^{- 2}$.
Accordingly, for small $\alpha p$, the second term will dominate initially,
until $k$ grows sufficiently large.

Since the right hand side eventually tends to zero, the theorem implies
convergence of the Ruppert-Polyak averaged iterates to the marginalized loss
minimizer $\wtbeta$, so the link between dropout and $\ell_2$-regularization
persists at the variance level. The averaging comes at the price of a slower
convergence rate $k\inv$ of the remainder terms, as opposed to the exponentially
fast convergence in Theorem \ref{thm::var_lim}. As in \eqref{eq::MSE_bd}, the
bound can be converted into a convergence rate for $\E\big[\norm{\rpbeta_k -
\wtbeta}_2^2\big]$ by taking the trace, \begin{align*}
  \E\Big[\norm[\big]{\rpbeta_k - \wtbeta}_2^2\Big] \leq d \left(\dfrac{2
  \norm{\bbX}^2 \cdot \norm[\big]{\E[X\tran \bfY \bfY\tran X]}}{k (1 - p)
  (\min_i \bbX_{ii})^4} + \dfrac{2 C}{k^2 (\alpha p (1 - p)
  \min_i\bbX_{ii})^3}\right).
\end{align*}

\subsection{Convergence of Simplified Dropout}
\label{sec::conv_whbeta}

To further illustrate how dropout and gradient descent are coupled, we will now
study the simplified dropout iterates \begin{align}
  \label{eq::whbeta_def_main}
  \whbeta_k = \whbeta_{k - 1} + \alpha D_k X\tran \big(\bfY - X \whbeta_{k -
  1}\big),
\end{align} as defined in \eqref{eq::wh_drop_def}. While the original dropout
reduces the model before taking the gradient, this version takes the gradient
first and applies dropout afterwards. As shown in Section
\ref{sec::drop_variations}, both versions coincide if the Gram matrix $\bbX$ is
diagonal. Recall from the discussion preceding Lemma \ref{lem::small_alpha_p}
that for diagonal $\bbX$, $\Cov\big(\wtbeta_k\big)$ converges to the optimal
covariance matrix. This suggests that for the simplified dropout, no additional
randomness in the limit $k\to \infty$ occurs.

The least squares objective $\bfbeta \mapsto \norm{\bfY - X\bfbeta}_2^2$ always
admits a minimizer, with any minimizer $\whbeta$ necessarily solving the
so-called normal equations $X\tran \bfY = \bbX \whbeta$. Provided $\bbX$ is
invertible, the least-squares estimator $\whbeta = \bbX\inv X\tran \bfY$ gives
the unique solution. We will not assume invertibility for all results below, so
we let $\whbeta$ denote any solution of the normal equations, unless specified
otherwise. In turn, \eqref{eq::whbeta_def_main} may be rewritten as
\begin{align}
  \label{eq::wh_rec_id}
  \whbeta_k - \whbeta = \big(I - \alpha D_k \bbX\big) \big(\whbeta_{k - 1} -
  \whbeta\big),
\end{align} which is simpler than the analogous representation of $\wtbeta_k$ as
a VAR process in \eqref{eq::RAR}.

As a first result, we will show that the expectation of the simplified dropout
iterates $\whbeta_k$ converges to the mean of the unregularized least squares
estimator $\whbeta$, provided that $\bbX$ is invertible. Indeed, using
\eqref{eq::wh_rec_id}, independence of $D_k$ and $\big(\whbeta_{k - 1} -
\whbeta\big)$, and $\E\big[D_k\big] = p I$, observe that \begin{align}
  \label{eq::whbeta_ind_step}
  \E\big[\whbeta_k - \whbeta\big] = \E\Big[\big(I - \alpha D_k \bbX\big)
  \big(\whbeta_{k - 1} - \whbeta\big)\Big] = \big(I - \alpha p \bbX\big)
  \E\big[\whbeta_{k - 1} - \whbeta\big].
\end{align} Induction on $k$ now shows $\E\big[\whbeta_k - \whbeta\big] = \big(I
- \alpha p \bbX\big)^k \E\big[\whbeta_0 - \whbeta\big]$ and so \begin{align*}
  \norm[\Big]{\E\big[\whbeta_k - \whbeta\big]}_2 = \norm[\Big]{\big(I -
  \alpha p \bbX\big)^k \E\big[\whbeta_0 - \whbeta\big]}_2
  \leq \norm[\big]{I - \alpha p \bbX}^k \norm[\Big]{\E\big[\whbeta_0 -
  \whbeta\big]}_2.
\end{align*} Assuming $\alpha p \norm{\bbX} < 1$, invertibility of $\bbX$
implies $\norm[\big]{I - \alpha p \bbX} < 1$. Consequently, the convergence is
exponential in the number of iterations.

Invertibility of $\bbX$ may be avoided if the initialization $\whbeta_0$ lies in
the orthogonal complement of the kernel of $\bbX$ and $\whbeta$ is the $\norm{\
\cdot\ }_2$-minimal solution to the normal equations. We can then argue that $(I
- \alpha p \bbX)^{k - 1} \E[\whbeta_0 - \whbeta]$ always stays in a linear
subspace on which $(I - \alpha p \bbX)$ still acts as a contraction.

We continue with our study of $\whbeta_k - \whbeta$ by employing the same
techniques as in the previous section to analyze the second moment. The linear
operator \begin{align}
  \label{eq::T_def}
  T(A) := \big(I - \alpha p \bbX\big) A \big(I - \alpha p \bbX\big) +
  \alpha^2 p (1 - p) \Diag\big(\bbX A \bbX\big),
\end{align} defined on $d\times d$ matrices, takes over the role of the affine
operator $S$ encountered in Lemma \ref{lem::wh_cov_rec}. In particular, the
second moments $A_k := \E\big[(\whbeta_k - \whbeta) (\whbeta_k -
\whbeta)\tran\big]$ evolve as the linear dynamical system \begin{align}
  \label{eq::T_rec_def}
  A_k = T\big(A_{k - 1}\big), \qquad k = 1, 2, \ldots
\end{align} without remainder terms. To see this, observe via
\eqref{eq::wh_rec_id} the identity $(\whbeta_k - \whbeta) (\whbeta_k -
\whbeta)\tran = (I - \alpha D_k \bbX) (\whbeta_{k - 1} - \whbeta) (\whbeta_{k -
1} - \whbeta)\tran (I - \alpha \bbX D_k)$. Taking the expectation on both sides,
conditioning on $D_k$, and recalling that $D_k$ is independent of $(\whbeta_k,
\whbeta)$ gives $A_k = \E\big[(I - \alpha D_k \bbX) A_{k - 1} (I - \alpha \bbX
D_k)\big]$. We have $\E[D_k] = p I_d$ and by Lemma \ref{lem::drop_moments},
$\E\big[D_k \bbX A_{k - 1} \bbX D_k\big] = p (\bbX A_{k-1} \bbX)_p = p^2 \bbX
A_{k - 1} \bbX + p (1 - p) \Diag(\bbX A_{k - 1} \bbX)$. Together with the
definition of $T(A)$, this proves \eqref{eq::T_rec_def}.

Further results are based on analyzing the recursion in \eqref{eq::T_rec_def}.
It turns out that convergence of $\whbeta_k$ to $\whbeta$ in second mean
requires a non-singular Gram matrix $\bbX$.

\begin{lemma}
  \label{lem::wh_div_ker}
  Suppose the initialization $\whbeta_0$ is independent of all other sources of
  randomness and the number of parameters satisfies $d \geq 2$, then there
  exists a singular $\bbX$, such that for any positive integer $k$,
  $\Cov\big(\whbeta_k\big) \geq \Cov\big(\whbeta_k - \whbeta\big) +
  \Cov\big(\whbeta\big)$ and \begin{align*}
    \norm[\Big]{\Cov\big(\whbeta_k - \whbeta\big)} \geq \alpha^2 p (1 - p).
  \end{align*}
\end{lemma}

For invertible $\bbX$, we can apply Theorem \ref{thm::gm} to show that
$\Cov\big(\whbeta\big) = \bbX\inv$ is the optimal covariance matrix for the
sequence of affine estimators $\whbeta_k$. The simplified dropout iterates
actually achieve the optimal variance when $\bbX$ is invertible, which stands in
contrast with the situation in Theorem \ref{thm::wtbeta_var_sub}.

\begin{theorem}
  \label{thm::wh_l2_conv}
  Suppose $\bbX$ is invertible, $\alpha \leq \min \Big\{\tfrac{1}{p
  \norm{\bbX}}, \tfrac{\lambdamin(\bbX)}{\norm{\bbX}^2}\Big\}$, and let
  $\whbeta_0$ be square-integrable, then, for any $k = 1, 2, \ldots$
  \begin{align*}
    \norm[\bigg]{\E\Big[\big(\whbeta_k - \whbeta\big) \big(\whbeta_k -
    \whbeta\big)\tran\Big]} \leq \big(1 - \alpha p\lambdamin(\bbX)\big)^k
    \norm[\bigg]{\E\Big[\big(\whbeta_0 - \whbeta\big) \big(\whbeta_0 -
    \whbeta\big)\tran\Big]}.
  \end{align*}
\end{theorem}

Intuitively, the result holds due to the operator $T$ in \eqref{eq::T_def} being
linear, as opposed to affine like in the case of Lemma \ref{lem::wt_var_rec}.
Choosing $\alpha$ sufficiently small ensures that $T$ acts as a contraction,
meaning $T^k(A) \to 0$ for any matrix $A$. Hence, linearity of $T$ serves as an
algebraic expression of the simplified dynamics. As in \eqref{eq::MSE_bd}, we
may take the trace to obtain the bound \begin{align*}
  \E\Big[\norm[\big]{\whbeta_k - \whbeta}_2^2\Big] \leq d \big(1 - \alpha
  p\lambdamin(\bbX)\big)^k \norm[\bigg]{\E\Big[\big(\whbeta_0 - \whbeta\big)
  \big(\whbeta_0 - \whbeta\big)\tran\Big]}.
\end{align*}


\section{Discussion and Outlook}
\label{sec::disc}

Our main contributions may be summarized as follows: We studied dropout in the
linear regression model, but unlike previous results, we explicitly analyzed the
gradient descent dynamics with new dropout noise being sampled in each
iteration. This allows us to characterize the limiting variance of the gradient
descent iterates exactly (Theorem \ref{thm::var_lim}), which sheds light on the
covariance structure induced via dropout. Our main tool in the analysis is a
particular recursion (Lemma \ref{lem::wt_var_rec}), which may be exhibited by
``anchoring'' the gradient descent iterates around the marginalized loss
minimizer $\wtbeta$. To further understand the interaction between noise and
gradient descent dynamics, we analyze the running average of the process
(Theorem \ref{thm::rp_l2_conv}) and a simplified version of dropout (Theorem
\ref{thm::wh_l2_conv}).

We view our analysis of the linear model as a fundamental first step towards
understanding the dynamics of gradient descent with dropout. Analyzing the
linear model has been a fruitful approach to study other phenomena in deep
learning, such as overfitting \citep{tsigler_bartlett_2023}, sharpness of local
minima \citep{bartlett_et_al_2023}, and in-context learning
\citep{zhang_et_al_2024}. We conclude by proposing some natural directions for
future work.

\textbf{Random minibatch sampling:} For yet another way of incorporating
dropout, we may compute the gradient based on a random subset of the data (mini
batches). In this case, the updating formula satisfies \begin{align}
  \ol{\bfbeta}_{k + 1} &= \ol{\bfbeta}_k - \alpha \nabla_{\ol \bfbeta_k}
  \dfrac{1}{2} \norm[\Big]{D_{k + 1} \big(\bfY - X \ol{\bfbeta}_k\big)}^2_2 =
  \ol{\bfbeta}_k + \alpha X\tran D_{k + 1} \big(\bfY - X \ol{\bfbeta}_k\big).
\end{align} The dropout matrices are now of dimension $n \times n$ and select a
random subset of the data points in every iteration. This version of dropout
also relates to randomly weighted least squares and resampling methods
\citep{duembgen_et_al_2013}. The update formula may be written in the form
$\ol{\bfbeta}_{k + 1} - \whbeta = (I-\alpha X\tran D_{k+1} X)(\ol{\bfbeta}_k -
\whbeta) + \alpha X\tran D_{k + 1}(\bfY - X \whbeta)$ with $\whbeta$ solving the
normal equations $X\tran \bfY = \bbX \whbeta$. Similarly to the corresponding
reformulation \eqref{eq::wtbeta_rewrite} of the original dropout scheme, this
defines a vector autoregressive process with random coefficients and lag one.

\textbf{Learning rates:} The proof ideas may be generalized to a sequence of
iteration-dependent learning rates $\alpha_k$. We expect this to come at the
cost of more involved formulas. Specifically, the operator $S$ in Lemma
\ref{lem::wt_var_rec} will depend on the iteration number through $\alpha_k$, so
the limit in Theorem \ref{thm::var_lim} cannot be expressed as $(\id -
S\lin)\inv S_0$ anymore.

\textbf{Random design and stochastic gradients:} We considered a fixed design
matrix $X$ and (full) gradient descent, whereas in machine learning practice
inputs are typically assumed to be random and parameters are updated via
stochastic gradient descent (SGD). The recent works
\cite{bos_schmidt-hieber_2023, schmidt-hieber_koolen_2023} derive convergence
rates for SGD considering linear regression and another form of noisy gradient
descent. We believe that parts of theses analyses carry over to dropout.

\textbf{Generic dropout distributions:} As already mentioned,
\cite{srivastava_et_al_2014} carry out simulations with dropout where $D_{ii}
\simiid \calN(1, 1)$. Gaussian dropout distributions are currently supported, or
easily implemented, in major software libraries \citep{abadi_et_al_2016,
chollet_2015, jia_et_al_2014, paszke_et_al_2019}. Analyzing a generic dropout
distribution with mean $\mu$ and variance $\sigma$ may also paint a clearer
picture of how the dropout noise interacts with the gradient descent dynamics.
For the linear regression model, results that marginalize over the dropout noise
generalize to arbitrary dropout distribution. In particular,
\eqref{eq::mlls_drop} turns into \begin{align*}
  \wtbeta :=&\ \argmin_{\bfbeta} \E\Big[\norm[\big]{\bfY - X D
  \bfbeta}^2_2\ \big\vert\ \bfY\Big] \\
  =&\ \argmin_{\bfbeta} \norm[\big]{\bfY - \mu X \bfbeta}^2_2 + \sigma^2
  \bfbeta\tran \Diag(\bbX) \bfbeta \\
  =&\ \mu \big(\mu^2 \bbX + \sigma^2 \Diag(\bbX)\big)\inv X\tran \bfY.
\end{align*} If $D_{ii} \simiid \calN(1, 1)$, then $\wtbeta = \big(\bbX +
\Diag(\bbX)\big)\inv X\tran \bfY$.

In contrast, treatment of the corresponding iterative dropout scheme seems more
involved. The analysis of dropout with Bernoulli distributions relies in parts
on the projection property $D^2 = D$. Without it, additional terms occur in the
moments in Lemma \ref{lem::drop_moments}, which is required for the computation
of the covariance matrix. For example, the formula $\E\big[D A D B D\big] = p
A_p B_p + p^2 (1 - p) \Diag(\ol{A} B)$ turns into \begin{align*}
  \E\big[D A D B D\big] &= \dfrac{1}{\mu_1} \big(\mu_1^2 A + \sigma^2
  \Diag(A)\big) \big(\mu_1^2 B + \sigma^2 \Diag(B)\big) + \sigma^2 \mu_1
  \Diag(\ol{A} \, \ol{B})\\
  &\qquad + \bigg(\mu_3 - \frac{\mu_2^2}{\mu_1}\bigg) \Diag(A) \Diag(B),
\end{align*} where $\mu_r$ denotes the $r$\textsuperscript{th} moment of the
dropout distribution and $\sigma^2$ its variance. For the Bernoulli
distribution, all moments equal $p$, so $\mu_3 - \mu_2^2 / \mu_1 = 0$ and the
last term disappears. Similarly, more terms will appear in the fourth moment of
$D$, making the expression for the operator corresponding to $S$ in Lemma
\ref{lem::wt_var_rec} more complicated.

\textbf{Inducing robustness via dropout:} Among the possible ways of explaining
the data, dropout should, by design, favor an explanation that is robust against
setting a random subset of the parameters to zero. \cite{mianjy_et_al_2018}
indicate that dropout in two-layer linear networks tends to equalize the norms
of different weight vectors.

To study the robustness properties of dropout, one may suggest analysis of loss
functions measuring prediction of the response vector if each estimated
regression coefficients is deleted with probability $p$. Given an estimator
$\whbeta$, a natural choice would be the loss \begin{align*}
  L\big(\whbeta,\starbeta\big) := \E\Big[\norm[\big]{X (D \whbeta -
  \starbeta)}_2^2\ \big\vert\ \bfY\Big] = \big(p \whbeta - \starbeta\big)\tran
  \bbX \big(p \whbeta - \starbeta\big) + p (1 - p) \whbeta\tran
  \Diag\big(\bbX\big) \whbeta,
\end{align*} with $D$ a new draw of the dropout matrix, independent of all other
randomness. This loss depends on the unknown true regression vector $\starbeta$.
Since $\E[\bfY] = X \starbeta$, an empirical version of the loss may replace $X
\starbeta$ with $\bfY$, considering $\E\big[\norm{X D \whbeta - \bfY}_2^2 \mid
\bfY\big]$. As shown in \eqref{eq::mlls_drop}, $\wtbeta = \bbX_p\inv X\tran
\bfY$ minimizes this loss function. This suggests that $\wtbeta$ may possess
some optimality properties for the loss $L(\ \cdot\ , \starbeta)$ defined above.

\textbf{Shallow networks with linear activation function:} Multi-layer neural
networks do not admit unique minimizers. In comparison with the linear
regression model, this poses a major challenge for the analysis of dropout.
\cite{mianjy_et_al_2018} consider shallow linear networks of the form $f(\bfx) =
U V\tran \bfx$ with $U = (\bfu_1, \ldots, \bfu_m)$ an $n \times m$ matrix and $V
= (\bfv_1, \ldots, \bfv_m)$ a $d\times m$ matrix. Suppose $D$ is an $m \times m$
dropout matrix. Assuming the random design vector $\bfX$ satisfies $\E\big[\bfX
\bfX\tran\big] = I_d$, and marginalizing over dropout noise applied to the
columns of $U$ (or equivalently to the rows of $V$) leads to an $\ell_2$-penalty
via \begin{equation}
  \label{eq::shallow_drop} 
  \begin{split}
    &\E\Big[\norm[\big]{\bfY - p\inv U D V\tran \bfX}_2^2 \mid \bfY\Big]\\
    =\ &\E\Big[\norm[\big]{\bfY - U V\tran \bfX}_2^2 \mid \bfY\Big] + \frac{1 -
    p}{p} \sum_{i = 1}^m \norm{\bfu_i}_2^2 \norm{\bfv_i}_2^2\\
    =\ &\E\Big[\norm[\big]{\bfY - U V\tran \bfX}_2^2 \mid \bfY\Big] + \frac{1 -
    p}{p} \Tr\Big(\Diag(U\tran U)\Diag(V\tran V)\Big).
  \end{split}
\end{equation} As an extension of our approach, it seems natural to investigate
whether gradient descent with dropout will converge to the same minimizer or
involve additional terms in the variance. In contrast with linear regression,
the marginalized loss \eqref{eq::shallow_drop} is non-convex and does not admit
a unique minimizer. Hence, we cannot simply center the gradient descent iterates
around a specific closed-form estimator, as in Section \ref{sec::gd_results}. To
extend our techniques we may expect to replace the centering estimator with the
gradient descent iterates for the marginalized loss function, demanding a more
careful analysis.

\cite{senen-cerda_sanders_2020} study the gradient flow associated with
\eqref{eq::shallow_drop} and exhibit exponential convergence of $U$ and $V$
towards a minimizer. Extending the existing result on gradient flows to gradient
descent is, however, non-trivial, see for example the gradient descent version
of Theorem 3.1 in \cite{bah_et_al_2022} provided in Theorem 2.4 of
\cite{nguegnang_et_al_2021}.

To be more precise, suppose $\bfY = W_\star \bfX + \bfeps$, where $\bfX$ and
$\bfeps$ are independent random vectors, so the task reduces to learning a
factorization $W_{\star} \approx U V\tran$ based on noisy evaluation of
$W_{\star} \bfX$. Consider the randomized loss \begin{align*}
  L\big(U, V\big) \mapsto \dfrac{1}{2} \norm[\Big]{\bfY - p\inv U D V\tran
  \bfX}_2^2,
\end{align*} with respective gradients \begin{align*}
  \nabla_{U} L\big(U, V\big) &= - \big(\bfY - p\inv U D V\tran \bfX\big)
  \bfX\tran V p\inv D\\
  \nabla_{V\tran} L\big(U, V\big) &= - p\inv D U\tran \big(\bfY - U p\inv D
  V\tran \bfX\big) \bfX\tran.
\end{align*} Given observations $(\bfY_k, \bfX_k) \simiid (\bfY, \bfX)$ and
independent dropout matrices $D_k \simiid D$, the factorized structure leads to
two coupled dynamical systems $U_{k + 1} = U_k - \alpha \nabla_{U_k} L(U_k,
V_k)$ and $V_{k + 1}\tran = V_k\tran - \alpha \nabla_{V_k\tran} L(U_k, V_k)$,
which are linked through the appearance of $V_k$ in $\nabla_{U_k} L(U_k, V_k)$
and $U_k$ in $\nabla_{V_k} L(U_k, V_k)$. Due to non-convexity of the underlying
marginalized objective \eqref{eq::shallow_drop}, the resulting dynamics should
be sensitive to initialization. Suppose $U_k = P_k(U_0, V_0)$ and $V_k =
Q_k(U_0, V_0)$ are given as random matrix polynomials $(P_k, Q_k)$ in $(U_0,
V_0)$, meaning finite sums of expressions like $A_1 X_{i_1} A_2 X_{i_2} \cdots
A_n X_{i_n} A_{n + 1}$, where $X_{i_j} \in \{U_0, U_0\tran, V_0, V_0\tran\}$ and
the $A_j$ are random coefficient matrices. Now, the gradient descent recursions
lead to \begin{equation}
  \label{eq::poly_rec}
  \begin{split}
    U_{k + 1} &= P_k(U_0, V_0) + \alpha \Big(\bfY_k  - P_k(U_0, V_0) p\inv D_k
    Q_k\tran(U_0, V_0) \bfX_k\Big) \bfX_k\tran Q_k(U_0, V_0) p\inv D_k\\
    &=: P_{k + 1}(U_0, V_0),\\
    V_{k + 1}\tran &= Q_k\tran(U_0, V_0) + \alpha p\inv D_k P_k\tran(U_0, V_0)
    \Big(\bfY_k - P_k(U_0, V_0) p\inv D_k Q_k\tran(U_0, V_0) \bfX_k\Big)
    \bfX_k\tran\\
    &=: Q_{k + 1}\tran(U_0, V_0),
  \end{split}
\end{equation} so $(U_{k + 1}, V_{k + 1})$ is also a polynomial in $(U_0, V_0)$
with random coefficients. A difficulty in analyzing this recursion is that the
degree of the polynomial increases exponentially fast. Indeed, since $P_{k + 1}$
includes the term $p^{- 2} P_k D_k Q_k\tran \bfX_k \bfX_k\tran Q_k D_k$, the
degree of $P_{k+1}$ is the degree of $P_k$ plus twice the degree of $Q_k$.
During each gradient descent step, additional randomness is introduced via the
newly sampled dropout matrix $D_k$ and the training data $(\bfY_k, \bfX_k)$.
Accordingly, the coefficients of $P_{k + 1}$ and $Q_{k + 1}$ fluctuate around
the coefficients of $\E[P_{k + 1} \mid U_k, V_k]$ and $\E[Q_{k + 1} \mid U_k,
V_k]$. A principled analysis of the resulting dynamics requires careful
accounting of how these fluctuations propagate through the iterations.
\cite{senen-cerda_sanders_2020} show that the gradient flow trajectories and
minimizers of \eqref{eq::shallow_drop} satisfy specific symmetries, so one
should hope to reduce the algebraic complexity of the problem by finding
analogous symmetries in the stochastic recursions \eqref{eq::poly_rec}.

Alternatively, one may consider layer-wise training of the weight matrices to
break the dependence between $U_k$ and $V_k$. Given $K_1 > 0$, suppose we keep
$U_0$ fixed while taking $K_1$ gradient steps \begin{align*}
  V_{k + 1}\tran = V_k\tran + \alpha p\inv D_k U_0\tran \big(\bfY - U_0 p\inv
  D_k V_k\tran \bfX_k\big) \bfX_k\tran
\end{align*} followed by $K_2 > 0$ gradient steps of the form \begin{align*}
  U_{k + 1} &= U_k + \alpha \big(\bfY - U_k p\inv D_k V_{K_1}\tran \bfX_k\big)
  \bfX_k\tran V_{K_1} p\inv D_k.
\end{align*} In each phase, the gradient descent recursion solves a linear
regression problem similar to to our analysis of the linear model. We leave the
details for future work.


\clearpage
\appendix

\section{Proofs for Section \ref{sec::avg_drop}}
\label{sec::avg_drop_proof}

\subsection{\texorpdfstring{Proof of Equation \eqref{eq::avg_obj_split}}{}}

Recall the definition $\wtbeta = \bbX_p\inv X\tran \bfY$. By Lemma
\ref{lem::drop_moments_A}, $\E[D\bbX D]=p \bbX_p$ and so \begin{align}
  \label{eq::del_loss_exp}
  \E\Big[D X\tran \big(\bfY - X D \wtbeta\big)\ \big\vert\ \bfY\Big] = p X\tran
  \bfY - p \bbX_p \wtbeta = 0.
\end{align} Note the identity $\bfY = X D \wtbeta + (\bfY - X D \wtbeta)$.
Hence, if $\whbeta$ denotes any estimator, $X D \whbeta - \bfY = X D (\whbeta -
\wtbeta) - (\bfY - X D \wtbeta)$. By \eqref{eq::del_loss_exp}, $\E\big[(\whbeta
- \wtbeta)\tran D X\tran (\bfY - X D \wtbeta) \mid \bfY\big] = 0$ and
\begin{align*}
  \E\Big[\norm[\big]{X D \whbeta - \bfY}_2^2\ \big\vert\ \bfY\Big] &=
  \E\Big[\norm[\big]{X D \wtbeta - \bfY}_2^2\ \big\vert\ \bfY\Big] +
  \E\Big[\norm[\big]{X D \big(\wtbeta - \whbeta\big)}_2^2\ \big\vert\
  \bfY\Big]\\
  &\qquad - 2 \E\Big[\big(\whbeta - \wtbeta\big)\tran D X\tran \big(\bfY - X D
  \wtbeta\big)\ \big\vert\ \bfY\Big]\\
  &= \E\Big[\norm[\big]{X D \wtbeta - \bfY}_2^2\ \big\vert\ \bfY\Big] +
  \E\Big[\norm[\big]{X D \big(\wtbeta - \whbeta\big)}_2^2\ \big\vert\ \bfY\Big],
\end{align*} which is the claimed expression. \CustomQED

\subsection{\texorpdfstring{Proof of Equation \eqref{eq::wtbeta_ones}}{}}

Let $r := \mathrm{rank}(X)$ and consider a singular value decomposition $X =
\sum_{\ell = 1}^r \sigma_\ell \bfv_\ell \bfw_\ell\tran$. The left-singular
vectors $\bfw_\ell$, $\ell = 1, \ldots, d$ are orthonormal, meaning $\bbX =
\sum_{\ell = 1}^r \sigma_\ell^2 \bfw_\ell \bfw_\ell\tran$ and since $\Diag(\bbX)
= \bbX_{11} \cdot I_d$, \begin{align*}
  \bbX_p = p \left(\sum_{\ell = 1}^{r} \sigma_\ell^2 \bfw_\ell
  \bfw_\ell\tran\right) + (1 - p) \bbX_{11} I_d.
\end{align*} Each left-singular vector $\bfw_\ell$ is an eigenvector of
$\bbX_p$, with associated eigenvalue $p \sigma_\ell^2 + (1 - p) \bbX_{11}$. If
$r < d$, suppose $\bfw_m$, $m = r + 1, \ldots, d$ complete the orthonormal
basis, then $\bbX_p \bfw_m = (1 - p) \bbX_{11} \bfw_m$. Consequently,
$\bbX_p\inv$ admits $\bfw_\ell$ as an eigenvector for every $\ell = 1, \ldots,
d$ and \begin{align*}
  \bbX_p\inv = \sum_{\ell = 1}^r \dfrac{1}{p \sigma_\ell^2 + (1 - p) \bbX_{11}}
  \bfw_\ell \bfw_\ell\tran + \sum_{m = r + 1}^d \dfrac{1}{(1 - p) \bbX_{11}}
  \bfw_m \bfw_m\tran.
\end{align*} By definition, $X\tran \bfY = \sum_{\ell = 1}^r \sigma_\ell
\bfw_\ell \bfv_\ell\tran \bfY$ and $\bfY = \sum_{\ell = 1}^r \sigma_\ell
\bfv_\ell \bfw_\ell\tran \starbeta + \bfeps$. Combining these facts now leads to
\begin{align*}
  \wtbeta = \sum_{\ell = 1}^r \dfrac{\sigma_\ell}{p
  \sigma_\ell^2 + (1 - p) \bbX_{11}} \bfw_\ell \bfv_\ell\tran \bfY
\end{align*} and \begin{align*}
  \E\big[X \wtbeta\big] &= \sum_{\ell = 1}^r \dfrac{\sigma^2_\ell}{p
  \sigma_\ell^2 + (1 - p) \bbX_{11}} \bfv_\ell \bfv_\ell\tran \E\big[\bfY\big]\\
  &= \sum_{\ell = 1}^r \dfrac{1}{p + (1 - p) \bbX_{11}/\sigma_\ell^2} \bfv_\ell
  \bfv_\ell\tran X \starbeta\\
  &= \sum_{\ell = 1}^r \dfrac{1}{p + (1 - p) \bbX_{11}/\sigma_\ell^2}
  \sigma_\ell \bfv_\ell \bfw_\ell\tran \starbeta.
\end{align*} 
\CustomQED

\subsection{\texorpdfstring{Proof of Equation \eqref{eq::data_miss_norm}}{}}

Set $A := X \big(\bbX + (p\inv - 1) \Diag(\bbX)\big)\inv X\tran$ and consider
$\bfw = \bfu + \bfv$ with $\bfu\tran \bfv = 0$ and $X\tran \bfv = \bfzero$.
Observe that \begin{align*}
  A^2 = A - X \big(\bbX + (p\inv - 1) \Diag(\bbX)\big)\inv (p\inv - 1)
  \Diag(\bbX) \big(\bbX + (p\inv - 1) \Diag(\bbX)\big)\inv X\tran
\end{align*} and thus $\bfw\tran A^2 \bfw = \bfu\tran A^2 \bfu \leq \bfu A \bfu
= \bfw\tran A\bfw$. If $\bfu$ is the zero vector, $0 = \bfw\tran A^2 \bfw =
\bfw\tran A\bfw$. Otherwise, $\Diag(\bbX) > 0$ implies $\big(\bbX + (p\inv - 1)
\Diag(\bbX)\big)\inv (p\inv - 1) \Diag(\bbX) \big(\bbX + (p\inv - 1)
\Diag(\bbX)\big)$ being positive definite, so we have the strict inequality
$\bfw\tran A^2 \bfw < \bfw\tran A \bfw$ whenever $\bfu \neq 0$.

Now, suppose that $A$ has an eigenvector $\bfw$ with corresponding eigenvalue
$\lambda \geq 1$, which implies $\bfu \neq \bfzero$. Then, we have $\bfw\tran
A^2 \bfw = \lambda^2 \geq \lambda = \bfw\tran A\bfw$. Equality only holds if
$\bfw\tran A \bfw = 1$. This contradicts the strict inequality $\bfw\tran A^2
\bfw < \bfw\tran A \bfw$ so all eigenvalues have to be strictly smaller than
one. \CustomQED


\section{Proofs for Section \ref{sec::gd_results}}
\label{sec::gd_results_proof}

\subsection{Proof of Lemma \ref{lem::wt_exp_conv}}

To show that $\norm[\big]{I - \alpha p \bbX_p} \leq 1 - \alpha p (1 - p) \min_i
\bbX_{ii} < 1$, note that $\norm{\bbX_p} \leq \norm{\bbX}$ by Lemma
\ref{lem::ol_spectral} and recall $\alpha p \norm{\bbX} < 1$ from
Assumption \ref{ass::sec_wtbeta_all}. Hence, \begin{align}
  \norm[\big]{I - \alpha p \bbX_p} = 1 - \alpha p \lambda_{\min}(\bbX_p).
  \label{eq::gd_mat_norm}
\end{align} For any vector $\bfv$ with $\norm{\bfv}_2 = 1$ and any two $d\times
d$ positive semi-definite matrices $A$ and $B$, $\bfv\tran (A + B)\bfv =
\bfv\tran A \bfv + \bfv\tran B \bfv \geq \lambdamin(A) + \lambdamin(B)$. Hence,
$\lambdamin(A + B) \geq \lambdamin(A) + \lambdamin(B)$ and
$\lambdamin(\bbX_p) \geq (1 - p) \lambdamin\big(\Diag(\bbX)\big) = (1 - p)
\min_i \bbX_{ii}$. By Assumption \ref{ass::sec_wtbeta_all}, the design matrix
$X$ has no zero columns, guaranteeing $\min_i \bbX_{ii} > 0$. Combined
with \eqref{eq::gd_mat_norm}, we now obtain $\norm[\big]{I - \alpha p
\bbX_p} \leq 1 - \alpha p (1 - p) \min_i \bbX_{ii} < 1$.

To prove the bound on the expectation, recall from \eqref{eq::wtbeta_rewrite}
that $\wtbeta_k - \wtbeta$ equals $\big(I - \alpha D_k \bbX_p\big)
\big(\wtbeta_{k - 1} - \wtbeta\big) + \alpha D_k \ol{\bbX} \big(p I - D_k\big)
\wtbeta_{k - 1}$. Lemma \ref{lem::drop_moments_A} shows that $\E[D A D] = p A_p$
and Lemma \ref{lem::ol_prop_B} gives $\ol{A}_p = p \ol{A}$. In turn, we have
$\E\big[D_k \ol{\bbX} D_k\big] = p \ol{\bbX}_p = p^2 \ol{\bbX}$ and
\begin{align}
  \E\big[D_k \ol{\bbX} (p I - D_k)\big] = p^2 \ol{\bbX} - p^2 \ol{\bbX} = 0.
  \label{eq::exp_ol_rewrite}
\end{align} Conditioning on all randomness except $D_k$ now implies
\begin{align}
  \E\big[\wtbeta_k - \wtbeta \mid \wtbeta, \wtbeta_{k - 1} \big] = \big(I -
  \alpha p \bbX_p\big) \big(\wtbeta_{k - 1} - \wtbeta\big).  
  \label{eq::exp_wtbeta_diff}
\end{align} By the tower rule $\E\big[\wtbeta_k - \wtbeta\big] = \big(I -
\alpha p \bbX_p\big) \E\big[\wtbeta_{k - 1} - \wtbeta\big]$, so induction
on $k$ gives \begin{align*}
  \E\big[\wtbeta_k - \wtbeta\big] = \big(I - \alpha p\bbX_p\big)^k
  \E\big[\wtbeta_0 - \wtbeta\big].
\end{align*} Sub-multiplicativity of the spectral norm implies
$\norm[\big]{(I - \alpha p \bbX_p)^k} \leq \norm*{I - \alpha p \bbX_p}^k$,
proving that $\norm[\big]{\E[\wtbeta_k - \wtbeta]}_2 \leq \norm[\big]{I -
\alpha p \bbX_p}^k \norm[\big]{\E[\wtbeta_0 - \wtbeta]}_2$. \CustomQED

\subsection{Proof of Theorem \ref{thm::gm}}

We have $\Cov\big(\wtbetaaff\big) = \Cov\big(\wtbeta_A + (\wtbetaaff -
\wtbeta_A)\big) = \Cov\big(\wtbeta_A\big) + \Cov\big(\wtbetaaff - \wtbeta_A\big)
+ \Cov\big(\wtbeta_A, \wtbetaaff - \wtbeta_A\big) + \Cov\big(\wtbetaaff -
\wtbeta_A, \wtbeta_A\big)$, so the triangle inequality implies \begin{align}
  \label{eq::aff_cov_upper}
  \norm[\Big]{\Cov\big(\wtbetaaff\big) - \Cov\big(\wtbeta_A\big) -
  \Cov\big(\wtbetaaff - \wtbeta_A\big)} \leq 2 \norm[\Big]{\Cov\big(\wtbetaaff -
  \wtbeta_A, \wtbeta_A\big)}.
\end{align} Write $B' := B - A X\tran$, then $\wtbetaaff - \wtbeta_A = B' \bfY +
\bfa$. When conditioned on $\bfY$, the estimator $\wtbeta_A = A X\tran
\bfY$ is deterministic. Hence, the law of total covariance yields
\begin{align*}
  \Cov\big(\wtbetaaff - \wtbeta_A, \wtbeta_A\big) &= \E\Big[\Cov\big(\wtbetaaff
  - \wtbeta_A, \wtbeta_A \mid \bfY\big)\Big] + \Cov\Big(\E\big[\wtbetaaff -
  \wtbeta_A \mid \bfY\big], \E\big[\wtbeta_A \mid \bfY\big]\Big)\\
  &= 0 + \Cov\big(\E[B'] \bfY + \E[\bfa], A X\tran \bfY\big) = \Cov\big(\E[B']
  \bfY, A X\tran \bfY\big).
\end{align*} Further, $\Cov(\bfY) = I$ implies $\Cov\big(\E[B'] \bfY, A X\tran
\bfY\big) = \E[B'] \Cov(\bfY) X A\tran = \E[B'] X A\tran$. Using $\wtbetaaff -
\wtbeta_A = B' \bfY + \bfa$, note that $\E[B'] X \starbeta =
\E_{\starbeta}\big[\wtbetaaff - \wtbeta_A\big] - \E_{\bfzero}\big[\wtbetaaff -
\wtbeta_A\big]$ with $\bfzero = (0, \ldots, 0)\tran$. Combining these
identities, sub-multiplicativity of the spectral norm, and the triangle
inequality leads to \begin{align*}
  \norm[\Big]{\Cov\big(\wtbetaaff - \wtbeta_A, \wtbeta_A\big)} &\leq \norm{A} \cdot
  \norm[\big]{\E[B'] X} \\
  &= \norm{A} \sup_{\starbeta : \norm{\starbeta}_2 \leq 1} \norm[\big]{\E[B'] X
  \starbeta}_2 \\
  &\leq \norm{A} \sup_{\starbeta : \norm{\starbeta}_2 \leq 1} \bigg(
  \norm[\Big]{\E_{\starbeta}\big[\wtbetaaff-\wtbeta_A\big]}_2 +
  \norm[\Big]{\E_{\bfzero}\big[\wtbetaaff-\wtbeta_A\big]}_2\bigg) \\
  &\leq 2 \norm{A} \sup_{\starbeta : \norm{\starbeta}_2 \leq 1}
  \norm[\Big]{\E_{\starbeta}\big[\wtbetaaff - \wtbeta_A\big]}_2.
\end{align*} Together with \eqref{eq::aff_cov_upper}, this proves the
result. \CustomQED

\subsection{Proof of Lemma \ref{lem::wt_var_rec}}

Given a random vector $U$ and a random element $V$, observe that $\E\big[\Cov(U
\mid V)\big] = \E\big[U U\tran\big] - \E\big[\E[U \mid V] \E[U \mid
V]\tran\big]$. Inserting $U = \wtbeta_k - \wtbeta$ and $V =
\big(\wtbeta_{k - 1}, \wtbeta\big)$, as well as defining $A_k :=
\E\big[(\wtbeta_k - \wtbeta) (\wtbeta_k - \wtbeta)\tran\big]$, leads to
\begin{align}
  \label{eq::wtbeta_rec_A_k}
  A_k = \E\Big[\E\big[\wtbeta_k - \wtbeta \mid \wtbeta_{k - 1}, \wtbeta\big]
  \E\big[\wtbeta_k - \wtbeta \mid \wtbeta_{k - 1}, \wtbeta\big]\tran\Big] +
  \E\Big[\Cov\big(\wtbeta_k - \wtbeta \mid \wtbeta_{k - 1}, \wtbeta\big)\Big]
\end{align} 
Recall $\E\big[\wtbeta_k - \wtbeta \mid \wtbeta, \wtbeta_{k-1}
\big] = (I - \alpha p\bbX_p) \big(\wtbeta_{k - 1} - \wtbeta\big)$ from
\eqref{eq::exp_wtbeta_diff}, and so 
\begin{align}
  \label{eq::wtbeta_rec_exp_cond}
  \E\Big[\E\big[\wtbeta_k - \wtbeta \mid \wtbeta_{k - 1}, \wtbeta\big]
  \E\big[\wtbeta_k - \wtbeta \mid \wtbeta_{k - 1}, \wtbeta\big]\tran \Big] =
  \big(I - \alpha p \bbX_p\big) A_{k - 1} \big(I - \alpha p \bbX_p\big),
\end{align} where $A_{k - 1} := \E\big[(\wtbeta_{k - 1} - \wtbeta) (\wtbeta_{k - 1}
- \wtbeta)\tran\big]$.

Evaluating the conditional covariance $\Cov\big(\wtbeta_k - \wtbeta \mid
\wtbeta_{k - 1}, \wtbeta\big)$ is the more challenging part, requiring
moments up to fourth order in $D_k$, see Lemma \ref{lem::drop_cov}. Recall that
\begin{align*}
  S(A) &= \big(I - \alpha p \bbX_p\big) A \big(I - \alpha p \bbX_p\big)\\
  &\quad + \alpha^2 p (1 - p) \Diag\big(\bbX_p A \bbX_p\big) + \alpha^2 p^2 (1 -
  p)^2 \bbX \odot \ol{A + \E\big[\wtbeta \wtbeta\tran\big]} \odot \bbX\\
  &\quad + \alpha^2 p^2 (1 - p) \bigg(\Big(\ol{\bbX} \Diag\Big(A +
  \E\big[\wtbeta \wtbeta\tran\big]\Big) \ol{\bbX}\Big)_p + \ol{\bbX}
  \Diag\big(\bbX_p A\big) + \Diag\big(\bbX_p A\big) \ol{\bbX}\bigg).
\end{align*}

\begin{lemma}
  \label{lem::exp_var_cond}
  For every positive integer $k$, \begin{align*}
    \E\Big[\Cov\big(\wtbeta_k - \wtbeta \mid \wtbeta_{k - 1}, \wtbeta\big)\Big]
    &= S(A_{k - 1}) - \big(I - \alpha p \bbX_p\big) A_{k - 1} \big(I - \alpha p
    \bbX_p\big) + \rho_{k - 1},
  \end{align*} with remainder $\rho_{k - 1}$ vanishing at the rate
  \begin{align*}
    \norm[\big]{\rho_{k - 1}} \leq 6 \norm[\big]{I - \alpha p \bbX_p}^{k - 1}
    \norm[\Big]{\E\big[(\wtbeta_0 - \wtbeta) \wtbeta\tran\big]}.
  \end{align*}
\end{lemma}

\begin{CustomProof}
  Recall from \eqref{eq::wtbeta_rewrite} that $\wtbeta_k - \wtbeta = \big(I-
  \alpha D_k \bbX_p\big) \big(\wtbeta_{k - 1} - \wtbeta\big) + \alpha D_k
  \ol{\bbX} \big(p I - D_k\big) \wtbeta_{k - 1}$. The covariance is invariant
  under deterministic shifts and sign flips, so \begin{align*}
    \Cov\Big(\wtbeta_k - \wtbeta\ \big\vert\ \wtbeta_{k - 1}, \wtbeta\Big) =
    \alpha^2 \Cov\Big(D_k \bbX_p \big(\wtbeta_{k - 1} - \wtbeta\big) + D_k
    \ol{\bbX} \big(D_k - p I\big) \wtbeta_{k - 1}\ \big\vert\ \wtbeta_{k -
    1},\wtbeta\Big).
  \end{align*} Applying Lemma \ref{lem::drop_cov} with $\bfu := \bbX_p
  \big(\wtbeta_{k - 1} - \wtbeta\big)$, $\ol{A} := \ol{\bbX}$, and $\bfv :=
  \wtbeta_{k - 1}$, we find \begin{equation}
    \label{eq::coveq}
    \begin{split}
      \dfrac{\Cov(\wtbeta_k - \wtbeta\ \vert\ \wtbeta_{k - 1},
      \wtbeta)}{\alpha^2 p (1 - p)} &= \Diag\Big(\bbX_p \big(\wtbeta_{k - 1} -
      \wtbeta\big) \big(\wtbeta_{k - 1} - \wtbeta\big)\tran
      \bbX_p\Big)\\
      &\quad + p \ol{\bbX} \Diag\Big(\bbX_p
      \big(\wtbeta_{k - 1} - \wtbeta\big) \wtbeta\tran_{k - 1}\Big)\\
      &\quad + p \Diag\Big(\wtbeta_{k - 1} \big(\wtbeta_{k - 1} -
      \wtbeta\big)\tran \bbX_p\Big) \ol{\bbX}\\
      &\quad + p \Big(\ol{\bbX} \Diag\big(\wtbeta_{k - 1}
      \wtbeta_{k - 1}\tran\big) \ol{\bbX}\Big)_p\\
      &\quad + p (1 - p) \bbX \odot \ol{\wtbeta_{k - 1}
      \wtbeta\tran_{k - 1}} \odot \bbX.
    \end{split}
  \end{equation} Set $B_{k - 1} := \E\big[(\wtbeta_{k - 1} - \wtbeta)
  \wtbeta\tran\big]$ and recall $A_{k - 1} = \E\big[(\wtbeta_{k - 1} -
  \wtbeta) (\wtbeta_{k - 1} - \wtbeta)\tran\big]$. Note the identities
  $\E\big[(\wtbeta_{k - 1} - \wtbeta) \wtbeta_{k - 1}\tran\big] = A_{k -
  1} + B_{k - 1}$ and $\E\big[\wtbeta_{k - 1} \wtbeta_{k - 1}\tran\big] =
  A_{k - 1} + \E\big[\wtbeta\wtbeta\tran\big] + B_{k - 1} + B_{k -
  1}\tran$. Taking the expectation of \eqref{eq::coveq}, multiplying
  both sides with $\alpha^2 p (1 - p)$, and using the definition of $S(A)$
  proves the claimed expression for $\E\big[\Cov(\wtbeta_k - \wtbeta \mid
  \wtbeta_{k - 1}, \wtbeta)\big]$ with remainder term \begin{equation}
    \label{eq::rho_eq}
    \begin{split}
      \rho_{k - 1} = \alpha^2 p (1 - p) \bigg( &p \ol{\bbX} \Diag\big(\bbX_p
      B_{k - 1}\big) + p \Diag\big( B_{k - 1}\tran \bbX_p\big) \ol{\bbX} + p
      \Big(\ol{\bbX} \Diag\big(B_{k - 1} + B_{k - 1}\tran\big)
      \ol{\bbX}\Big)_p\\
      &+ p (1 - p) \bbX \odot \ol{\big(B_{k - 1}+B_{k - 1}\tran\big)} \odot
      \bbX\bigg).
    \end{split}
  \end{equation}

  For any $d \times d$ matrices $A$ and $B$, Lemma \ref{lem::ol_spectral}
  provides the inequalities $\norm[\big]{\Diag(A)} \leq \norm{A}$, $\norm{A_p}
  \leq \norm{A}$, and $\norm{A \odot B} \leq \norm{A} \cdot \norm{B}$. If $A$ is
  moreover positive semi-definite, then also $\norm[\big]{\ol{A}} \leq
  \norm*{A}$. Combined with the sub-multiplicativity of the spectral norm, this
  implies \begin{align*}
    \norm[\big]{\ol{\bbX} \Diag(\bbX_p B_{k - 1})} &\leq \norm[\big]{\ol{\bbX}}
    \cdot \norm[\big]{\Diag(\bbX_p B_{k - 1})}\\
    &\leq \norm{\bbX} \cdot \norm{\bbX_p} \cdot \norm[\big]{B_{k - 1}} \leq
    \norm{\bbX}^2 \cdot \norm[\big]{B_{k-1}}\\ \norm[\Big]{\Big(\ol{\bbX}
    \Diag\big(B_{k - 1} + B_{k - 1}\tran\big) \ol{\bbX}\Big)_p} &\leq
    \norm[\big]{\ol{\bbX} \Diag\big(B_{k - 1} + B_{k - 1}\tran\big)
    \ol{\bbX}}\\
    &\leq \norm[\big]{\ol{\bbX}}^2 \cdot \norm[\big]{B_{k - 1} + B_{k - 1}\tran}
    \leq 2 \norm{\bbX}^2 \cdot \norm[\big]{B_{k - 1}}\\
    \norm[\Big]{\bbX \odot \ol{\big(B_{k - 1} + B_{k - 1}\tran\big)}
    \odot \bbX} &\leq 2 \norm{\bbX}^2 \cdot \norm[\big]{B_{k - 1}}.
  \end{align*} By Assumption \ref{ass::sec_wtbeta_all} also $\alpha p
  \norm{\bbX} < 1$, so combining the upper-bounds with \eqref{eq::rho_eq} leads
  to \begin{align}
    \label{eq::rho_upper}
    \norm[\big]{\rho_{k - 1}} \leq 6 (\alpha p)^2 \norm{\bbX}^2 \cdot
    \norm[\big]{B_{k - 1}} \leq 6 \norm[\big]{B_{k - 1}}.
  \end{align}

  The argument is to be completed by bounding $\norm[\big]{B_{k - 1}}$. Using
  \eqref{eq::wtbeta_rewrite}, we have $\big(\wtbeta_{k - 1} - \wtbeta\big)
  \wtbeta\tran = \big(I - \alpha D_{k - 1} \bbX_p\big) \big(\wtbeta_{k - 2} -
  \wtbeta\big) \wtbeta\tran + \alpha D_{k - 1} \ol{\bbX} \big(p I - D_{k -
  1}\big) \wtbeta_{k - 1} \wtbeta\tran$. In \eqref{eq::exp_ol_rewrite}, it
  was shown that $\E\big[D_{k - 1} \ol{\bbX} (p I - D_{k-1})\big] = 0$.
  Recalling that $D_{k - 1}$ is independent of $\big(\wtbeta, \wtbeta_{k -
  2}\big)$ and $\E[D_{k - 1}] = p I_d$, we obtain $B_{k - 1} = \big(I - \alpha p
  \bbX_p\big) B_{k - 2}$. By induction on $k$, $B_{k - 1} = \big(I - \alpha p
  \bbX_p\big)^{k - 1} \E\big[(\wtbeta_0 - \wtbeta\big) \wtbeta\tran\big]$. Using
  sub-multiplicativity of the spectral norm, \begin{align*}
    \norm[\big]{B_{k - 1}} \leq \norm{I - \alpha p \bbX_p}^{k - 1}
    \norm[\Big]{\E\big[(\wtbeta_0-\wtbeta) \wtbeta\tran\big]}.
  \end{align*} Together with \eqref{eq::rho_upper} this finishes the proof. 
\end{CustomProof}

Combining Lemma \ref{lem::exp_var_cond} with \eqref{eq::wtbeta_rec_A_k} and
\eqref{eq::wtbeta_rec_exp_cond} leads to $\norm[\big]{A_k - S(A_{k - 1})} =
\norm[\big]{\rho_{k - 1}}$ with remainder $\rho_{k - 1}$ as above. This
completes the proof of Lemma \ref{lem::wt_var_rec}. \CustomQED

\subsection{Proof of Theorem \ref{thm::var_lim}}

Let $S : \R^{d \times d} \to \R^{d \times d}$ be the affine operator introduced
in Lemma \ref{lem::wt_var_rec} and recall the definitions $S_0 := S(0)$ and
$S\lin(A) := S(A) - S_0$. First, the operator norm of $S\lin$ will be analyzed.

\begin{lemma}
  \label{lem::S_op_bound}
  The linear operator $S\lin$ satisfies $\norm*{S\lin}_{\op} \leq
  \norm[\big]{I - \alpha p\bbX_p} < 1$, provided that \begin{align*}
    \alpha < \min \left\{\dfrac{1}{p \norm{\bbX}}, \dfrac{\lambdamin(\bbX_p)}{3
    \norm{\bbX}^2}\right\},
  \end{align*} where $\lambdamin(\bbX_p)$ denotes the smallest eigenvalue of
  $\bbX_p$.
\end{lemma}

\begin{CustomProof}
  Let $A$ be a $d \times d$ matrix. Applying the triangle inequality,
  Lemma \ref{lem::ol_spectral}, and sub-multiplicativity of the spectral norm,
  \begin{align*}
    \norm[\big]{S\lin(A)} &\leq \norm[\big]{I - \alpha p \bbX_p}^2 \norm*{A} +
    \Big(\alpha^2 p (1 - p) + 3 \alpha^2 p^2 (1 - p) + \alpha^2 p^2 (1 -
    p)^2\Big) \norm*{\bbX}^2 \norm*{A} \nonumber\\
    &\leq \Big(\norm[\big]{I - \alpha p \bbX_p}^2 + 2 \alpha^2 p
    \norm*{\bbX}^2\Big) \norm*{A},
  \end{align*} where the second inequality follows from $p (1 - p) \leq 1/4$.

  As shown in \eqref{eq::gd_mat_norm}, $\norm[\big]{I - \alpha p \bbX_p} =
  1 - \alpha p \lambdamin(\bbX_p)$. Lemma \ref{lem::ol_spectral} now
  implies $\big(1 - \alpha p \lambdamin(\bbX_p)\big)^2 = 1 - 2 \alpha p
  \lambdamin(\bbX_p) + \alpha^2 p^2 \lambdamin(\bbX_p) \leq 1 - 2 \alpha p
  \lambdamin(\bbX_p)^2 + \alpha^2 p \norm{\bbX}^2$, so that
  \begin{align*}
    \norm[\big]{S\lin(A)} \leq \big(1 - 2\alpha p \lambdamin(\bbX_p) + 3
    \alpha^2 p \norm{\bbX}^2\big) \norm*{A}.
  \end{align*} If $\alpha < \lambdamin(\bbX_p) / \big(3 \norm{\bbX}^2\big)$,
  then also $3 \alpha^2 p \norm{\bbX}^2 \leq 3 \alpha p \lambdamin(\bbX_p)$, so
  that in turn $\norm{S\lin}_{\op} \leq \norm[\big]{I - \alpha p\bbX_p}$. The
  constraint $\alpha < 1 / \big(p \norm{\bbX}\big)$ now enforces $\alpha p
  \norm{\bbX} < 1$, which implies $\norm[\big]{I - \alpha p \bbX_p} < 1$.
\end{CustomProof}

As before, set $A_k := \E\big[(\wtbeta_k - \wtbeta) (\wtbeta_k -
\wtbeta)\tran\big]$ for each $k \geq 0$ and let $\rho_k := A_{k + 1} - S(A_k)$.
Using induction on $k$, we now prove \begin{align}
  \label{eq::S_no_neumann}
  A_k &= S\lin^k(A_0) + \sum_{\ell = 0}^{k - 1} S\lin^{\ell} \big(S_0 +
  \rho_{k - 1 - \ell}\big).
\end{align} Taking $k = 1$, $A_1 = S(A_0) + \rho_0 = S\lin(A_0) + S_0 + \rho_0$,
so the claimed identity holds. Assuming the identity is true for $k - 1$, the
recursion $A_k = S(A_{k - 1}) + \rho_{k - 1}$ leads to \begin{align*}
  A_k &= S\left(S\lin^{k - 1}(A_0) + \sum_{\ell = 0}^{k - 2} S\lin^{\ell}
  \Big(S_0 + \rho_{k - 2 - \ell}\Big)\right) + \rho_{k - 1}\\
  &= S\lin^k (A_0) + S\lin \left(\sum_{\ell = 0}^{k - 2} S\lin^{\ell} \Big(S_0 +
  \rho_{k - 2 - \ell}\Big)\right) + S_0 + \rho_{k - 1}\\
  &= S\lin^k (A_0) + \sum_{\ell = 0}^{k - 1} S\lin^{\ell} \Big(S_0 + \rho_{k - 1
  - \ell}\Big), 
\end{align*} thereby establishing the induction step and proving
\eqref{eq::S_no_neumann}.
  
Assuming $\norm*{S\lin}_{\op} \leq \norm{I-\alpha p \bbX_p} < 1$, we move on to
show the bound \begin{align*}
  \norm[\Big]{A_k - \big(\id - S\lin\big)\inv S_0} \leq \norm[\big]{I - \alpha p
  \bbX_p}^k \norm[\Big]{A_0 - \big(\id - S\lin\big)\inv S_0} + C_0 k \norm*{I - \alpha
  p\bbX_p}^{k - 1}
\end{align*} with $\id$ the identity operator on $\R^{d \times d}$, and $C_0 :=
6\norm[\big]{\E\big[(\wtbeta_0 - \wtbeta) \wtbeta\tran\big]}$. By linearity,
$\sum_{\ell = 0}^{k - 1} S\lin^{\ell} \big(S_0 + \rho_{k - 1 - \ell}\big) =
\sum_{\ell = 0}^{k - 1} S\lin^{\ell}(S_0) + \sum_{m = 0}^{k - 1} S\lin^m(\rho_{k
- 1 - m})$. Since $\norm*{S\lin}_{\op} < 1$, Lemma \ref{lem::neumann}
asserts that $\big(\id - S\lin\big)\inv = \sum_{\ell = 0}^\infty
S\lin^\ell$ and \begin{align}
  \label{eq::S_op_neumann}
  \sum_{\ell = 0}^{k - 1} S\lin^{\ell}\big(S_0\big) &= \big(\id -
  S\lin\big)\inv S_0 + \left(\sum_{\ell = 0}^{k - 1} S\lin^{\ell} - \big(\id -
  S\lin\big)\inv\right) S_0 \nonumber\\
  &= \big(\id - S\lin\big)\inv S_0 - \sum_{\ell = k}^\infty
  S\lin^\ell\big(S_0\big) \nonumber\\
  &= \big(\id - S\lin\big)\inv S_0 - S\lin^k \Big(\big(\id - S\lin\big)\inv
  S_0\Big).
\end{align}

Lemma \ref{lem::wt_var_rec} ensures $\norm{\rho_{k - 1 - m}} \leq C_0
\norm[\big]{I - \alpha p \bbX_p}^{k - 1 - m}$ for all $m \leq k - 1$. Moreover,
$\norm{S\lin}_{\op} \leq \norm[\big]{I - \alpha p \bbX_p}$ and hence
$\norm[\big]{S\lin^m(\rho_{k - 1 - m})} \leq C_0 \norm[\big]{I - \alpha p
\bbX_p}^{k - 1}$ for every $m = 0, 1, \ldots$, so the triangle inequality
implies \begin{align}
  \label{eq::S_rho_bound}
  \norm*{\sum_{m = 0}^{k - 1} S\lin^m(\rho_{k - 1 - m})} \leq C_0 k
  \norm[\big]{I-\alpha p \bbX_p}^{k - 1}.
\end{align} 

Combining \eqref{eq::S_no_neumann} and \eqref{eq::S_op_neumann}, as well as
applying the triangle inequality and the bound \eqref{eq::S_rho_bound},
leads to the first bound asserted in Theorem \ref{thm::var_lim}, \begin{align}
  \norm[\Big]{A_k - \big(\id - S\lin\big)\inv S_0} &\leq \norm[\bigg]{S\lin^k
  \Big(A_0 - \big(\id - S\lin\big)\inv S_0\Big)} + \norm*{\sum_{m = 0}^{k - 1}
  S\lin^m(\rho_{k - 1 - m})} \nonumber \\
  &\leq \norm[\big]{I - \alpha p \bbX_p}^k \norm[\Big]{A_0 - \big(\id -
  S\lin\big)\inv S_0} + C_0 k \norm[\big]{I-\alpha p \bbX_p}^{k - 1}.
  \label{eq::sec_mom_bound}
\end{align}

To show the corresponding bound for the variance, observe that
$\Cov\big(\wtbeta_k - \wtbeta\big) = A_k - \E\big[\wtbeta_k - \wtbeta\big]
\E\big[\wtbeta_k - \wtbeta\big]\tran$. Lemma \ref{lem::wt_exp_conv} and
\eqref{eq::tensor_cont} imply \begin{align*}
  \norm[\Big]{\Cov\big(\wtbeta_k - \wtbeta\big) - A_k} &=
  \norm[\Big]{\E\big[\wtbeta_k - \wtbeta\big] \E\big[\wtbeta_k -
  \wtbeta\big]\tran}\\
  &\leq \norm[\Big]{\E\big[\wtbeta_k - \wtbeta\big]}_2^2\\
  &\leq \norm[\big]{I - \alpha p \bbX_p}^{k - 1} \norm[\Big]{\E\big[\wtbeta_0 -
  \wtbeta\big]}_2^2.
\end{align*} Together with \eqref{eq::sec_mom_bound} and the triangle
inequality, this proves the second bound asserted in Theorem \ref{thm::var_lim}.
\CustomQED

\subsection{Proof of Lemma \ref{lem::small_alpha_p}}

Applying Theorem \ref{thm::gm}, Lemma \ref{lem::wt_exp_conv}, and
the triangle inequality, \begin{align}
  \label{eq::cov_conv_rate}
  \norm[\Big]{\Cov\big(\wtbeta_k\big) - \Cov\big(\wtbeta\big)} \leq
  \norm[\Big]{\Cov\big(\wtbeta_k - \wtbeta\big)} + 4 \norm[\big]{\bbX_p\inv}
  \norm[\big]{I - \alpha p \bbX_p}^k \sup_{\starbeta : \norm{\starbeta}_2
  \leq 1} \norm[\Big]{\E_{\starbeta}\big[\wtbeta_0 - \wtbeta\big]}_2.
\end{align} Lemma \ref{lem::ol_spectral} implies $\bbX_p \geq (1 - p)
\Diag(\bbX)$, so that $\norm[\big]{\bbX_p\inv} = \lambdamin(\bbX_p)\inv \leq
\big((1 - p) \min_i \bbX_{ii}\big)\inv$. Next, $\E_{\starbeta}\big[\wtbeta\big]
= \bbX_p\inv \bbX \starbeta$ entails equality between $\sup_{\starbeta :
\norm{\starbeta}_2 \leq 1} \norm[\big]{\E_{\starbeta}[\wtbeta]}_2$ and
$\norm[\big]{\bbX_p\inv \bbX} \leq \big((1 - p) \min_i \bbX_{ii}\big)\inv
\norm{\bbX}$. The second term on the right-hand side of
\eqref{eq::cov_conv_rate} is then bounded by $C_1 \norm[\big]{I - \alpha p
\bbX_p}^k / (1 - p)^2$, for some constant $C_1$ independent of $(\alpha, p, k)$.

To prove the first claim of the lemma, it now suffices to show
\begin{align}
  \label{eq::joint_cov_bound}
  \norm[\Big]{\Cov\big(\wtbeta_k - \wtbeta\big)} \leq \dfrac{1}{(1 - p)^2}
  \Big(k \norm[\big]{I - \alpha p \bbX_p}^{k - 1}C_2 + \alpha p C_3\Big),
\end{align} where $C_2$ and $C_3$ are constants independent of $(\alpha, p, k)$.
As $\norm[\big]{\bbX_p\inv} \leq \big((1 - p) \min_i \bbX_{ii}\big)\inv$, the
constant $C$ in Theorem \ref{thm::var_lim} satisfies $C \leq C_4 / (1 - p)^2 +
\norm[\big]{(\id - S\lin)\inv S_0}$, with $C_4$ depending only on the
distribution of $\big(\bfY, \wtbeta_0, X\big)$. Consequently, Theorem
\ref{thm::var_lim} and the triangle inequality imply \begin{equation}
  \label{eq::joint_cov_bound_init}
  \begin{split}
    \norm[\Big]{\Cov\big(\wtbeta_k - \wtbeta\big)} &\leq
    \norm[\Big]{\Cov\big(\wtbeta_k - \wtbeta) - \big(\id - S\lin\big)\inv S_0} +
    \norm[\Big]{\big(\id - S\lin\big)\inv S_0}\\
    &\leq k \norm[\big]{I - \alpha p \bbX_p}^{k - 1} \left(\dfrac{C_4}{(1 -
    p)^2} + \norm[\Big]{\big(\id - S\lin\big)\inv S_0}\right) +
    \norm[\Big]{\big(\id - S\lin\big)\inv S_0}.
  \end{split}
\end{equation}

Consider a bounded linear operator $G$ on $\R^{d\times d}$ satisfying
$\norm{G}_{\op} < 1$. For an arbitrary $d \times d$ matrix $A$, Lemma
\ref{lem::neumann} asserts $(\id - G)\inv A = \sum_{\ell = 0}^\infty
G^{\ell}(A)$ and therefore $\norm[\big]{(\id - G)\inv A} \leq \sum_{\ell =
0}^\infty \norm{G}_{\op}^\ell \cdot \norm{A} = \big(1 - \norm{G}_{\op}\big)\inv
\norm{A}$. Theorem \ref{thm::gm} states that $\norm*{S\lin}_{\op} \leq
\norm[\big]{I - \alpha p \bbX_p}$. As shown following \eqref{eq::gd_mat_norm},
$\norm[\big]{I - \alpha p \bbX_p} \leq 1 - \alpha p (1 - p) \min_i \bbX_{ii}$.
Therefore, \begin{align*}
  \norm[\Big]{\big(\id - S\lin\big)\inv S_0} \leq \big(1 -
  \norm{S\lin}_{\op}\big)\inv \norm{S_0} \leq \big(\alpha p (1 - p) \min_i
  \bbX_{ii}\big)\inv \norm{S_0}.
\end{align*} Taking $A = 0$ in Lemma \ref{lem::wt_var_rec}, $S_0 = \alpha^2 p^2
(1 - p) \big(\ol{\bbX} \Diag\big(\E\big[\wtbeta \wtbeta\tran\big]\big)
\ol{\bbX}\big)_p + \alpha^2 p^2 (1 - p)^2 \bbX \odot \ol{\E[\wtbeta
\wtbeta\tran]} \odot \bbX$. Using Lemma \ref{lem::ol_spectral} and
$\norm[\big]{\bbX_p\inv} \leq \big((1 - p) \min_i \bbX_{ii}\big)\inv$,
\begin{align*}
  \norm{S_0} &\leq \alpha^2 p^2 (1 - p) \Big(\norm[\Big]{\ol{\bbX}
  \Diag\big(\E\big[\wtbeta \wtbeta\tran\big]\big) \ol{\bbX}} + \norm[\Big]{\bbX
  \odot \ol{\E\big[\wtbeta \wtbeta\tran\big]} \odot \bbX}\Big)\\
  &\leq \dfrac{(\alpha p \norm{\bbX})^2 \norm[\big]{\E[X\tran \bfY \bfY\tran
  X]}}{(1 - p) (\min_i \bbX_{ii})^2}
\end{align*} proving that \begin{align}
  \label{eq::extra_var_bound}
  \norm[\big]{(\id - S\lin)\inv S_0} \leq \alpha p (1 - p)^{-2}  (\min_i
  \bbX_{ii})^{-3} \norm{\bbX}^2 \cdot \norm[\big]{\E[X\tran \bfY \bfY\tran X]}.
\end{align} Note that $\alpha p \norm{\bbX}^2 \leq \norm{\bbX}$ by Assumption
\ref{ass::sec_wtbeta_all}. Applying these bounds in
\eqref{eq::joint_cov_bound_init} leads to \begin{align*}
  \norm*{\Cov\big(\wtbeta_k - \wtbeta\big)} &\leq \dfrac{k \norm{I - \alpha p
  \bbX_p}^{k - 1}}{(1 - p)^2} \left(C_4 + \dfrac{\norm[\big]{\bbX}
  \norm[\big]{\E[X\tran \bfY \bfY\tran X]}}{(\min_i \bbX_{ii})^3}\right)\\
  &\qquad + \dfrac{\alpha p \norm[\big]{\bbX}^2 \norm[\big]{\E[X\tran \bfY
  \bfY\tran X]}}{(1 - p)^2 (\min_i \bbX_{ii})^3},
\end{align*} which proves \eqref{eq::joint_cov_bound}. Combined with
\eqref{eq::cov_conv_rate}, this proves the first claim of the lemma since
\begin{align}
  \label{eq::cov_bound_small_ap}
  \norm[\Big]{\Cov\big(\wtbeta_k\big) - \Cov\big(\wtbeta\big)} \leq \dfrac{k
  \norm{I - \alpha p \bbX_p}^{k - 1} (C_1 + C_2) + \alpha p C_3}{(1 - p)^2}.
\end{align}
  
To start proving the second claim, recall that $\Cov\big(\wtbeta\big) =
\bbX_p\inv \bbX \bbX_p\inv$. Hence, the triangle inequality leads to
\begin{align*}
  &\norm[\Big]{\Cov\big(\wtbeta_k\big) - \Diag(\bbX)\inv \bbX \Diag(\bbX)\inv}\\
  &\leq \norm[\Big]{\Cov\big(\wtbeta_k\big) - \Cov\big(\wtbeta\big)} +
  \norm[\Big]{\Diag(\bbX)\inv \bbX \Diag(\bbX)\inv - \bbX_p\inv \bbX
  \bbX_p\inv}.
\end{align*} Let $A$, $B$, and $C$ be square matrices of the same dimension,
with $A$ and $B$ invertible. Observe the identity $A\inv C A\inv - B\inv C B\inv
= A\inv (B - A) B\inv C A\inv + B\inv C A\inv (B - A) B\inv$, so
sub-multiplicativity implies $\norm[\big]{A\inv C A\inv - B\inv C B\inv} \leq 2
\max \big\{\norm{A\inv}, \norm{B\inv}\big\} \norm{A\inv} \cdot \norm{B\inv}
\cdot \norm{A - B} \cdot \norm{C}$. Using $\norm[\big]{\bbX_p\inv} \leq \big((1
- p) \min_i \bbX_{ii}\big)\inv$, and inserting $A = \Diag(\bbX)$, $B = \bbX_p =A
+ p \ol{\bbX}$, and $C = \bbX$ results in \begin{align*}
  \norm[\Big]{\Diag(\bbX)\inv \bbX \Diag(\bbX)\inv - \bbX_p\inv \bbX \bbX_p\inv}
  \leq \dfrac{p C_5}{(1 - p)^2}
\end{align*} with $C_5$ independent of $(\alpha, p, k)$. Combined with
\eqref{eq::cov_bound_small_ap}, this results in \begin{align*}
  \norm[\Big]{\Cov\big(\wtbeta_k\big) - \Diag(\bbX)\inv \bbX \Diag(\bbX)\inv}
  \leq \dfrac{k \norm{I - \alpha p \bbX_p}^{k - 1} (C_1 + C_2) + \alpha p
  C_3 + p C_5}{(1 - p)^2},
\end{align*} which proves the second claim of the lemma by enlarging $C''$, if
necessary. \CustomQED

\subsection{Proof of Theorem \ref{thm::wtbeta_var_sub}}

Start by noting that $\lambdamin(A) = \inf_{\bfv : \norm{\bfv} = 1} \bfv\tran A
\bfv$ for symmetric matrices, see \cite{horn_johnson_2013}, Theorem 4.2.6. Using
super-additivity of infima, observe the lower bound \begin{align}
  &\liminf_{k \to \infty} \inf_{\bfv : \norm{\bfv} = 1} \bfv\tran
  \Big(\Cov\big(\wtbeta_k\big) - \Cov\big(\wtbeta\big)\Big) \bfv
  \nonumber\\
  \geq &\liminf_{k \to \infty} \bigg(\lambdamin \Big(\Cov\big(\wtbeta_k -
  \wtbeta\big)\Big) - \sup_{\bfv : \norm{\bfv} = 1} \abs[\Big]{\bfv\tran
  \Big(\Cov\big(\wtbeta_k\big) - \Cov\big(\wtbeta\big) - \Cov\big(\wtbeta_k -
  \wtbeta\big)\Big)\bfv}\bigg) \nonumber\\
  \geq &\liminf_{k \to \infty} \lambdamin \Big(\Cov\big(\wtbeta_k -
  \wtbeta\big)\Big) - \limsup_{k \to \infty}
  \norm[\Big]{\Cov\big(\wtbeta_k\big) - \Cov\big(\wtbeta\big) -
  \Cov\big(\wtbeta_k - \wtbeta\big)}. \label{eq::cov_diff_lower}
\end{align} Combining Lemma \ref{lem::wt_exp_conv} and Theorem \ref{thm::gm},
the limit superior in \eqref{eq::cov_diff_lower} vanishes. Further,
$\Cov\big(\wtbeta_k - \wtbeta\big)$ converges to $\big(\id - S\lin\big)\inv S_0$
by Theorem \ref{thm::var_lim}, so it suffices to analyze the latter matrix.

For the next step, the matrix $S_0 := S(0)$ in Theorem \ref{thm::var_lim} will
be lower-bounded. Taking $A = 0$ in Lemma \ref{lem::wt_var_rec} and exchanging
the expectation with the $\Diag$ operator results in \begin{align}
  S_0 &= \alpha^2 p^2 (1 - p) \bigg(\ol{\bbX} \Diag\Big(\E\big[\wtbeta
  \wtbeta\tran\big]\Big) \ol{\bbX}\bigg)_p + \alpha^2 p^2 (1 - p)^2 \bbX \odot
  \ol{\E\big[\wtbeta \wtbeta\tran\big]} \odot \bbX \nonumber \\
  &= \alpha^2 p^2 (1 - p) \E\Bigg[p\ol{\bbX} \Diag\big(\wtbeta
  \wtbeta\tran\big) \ol{\bbX} + (1 - p) \bigg(\Diag\Big(\ol{\bbX}
  \Diag\big(\wtbeta \wtbeta\tran\big) \ol{\bbX}\Big) + \bbX \odot \ol{\wtbeta
  \wtbeta\tran} \odot \bbX\bigg)\Bigg].
  \label{eq::lambda_zero_expd}
\end{align}

The first matrix in \eqref{eq::lambda_zero_expd} is always positive
semi-definite and we will now lower bound the matrix $B:= \Diag\big(\ol{\bbX}
\Diag(\wtbeta \wtbeta\tran) \ol{\bbX}\big) + \bbX \odot \ol{\wtbeta
\wtbeta\tran} \odot \bbX $. Given distinct $i, j = 1, \ldots, d$, symmetry of
$\bbX$ implies \begin{align*}
  \Big(\ol{\mathbb{X}} \Diag\big(\wtbeta \wtbeta\tran\big)
  \ol{\mathbb{X}}\Big)_{ii} &= \sum_{k = 1}^d \ol{\bbX}_{ik} \Diag\big(\wtbeta
  \wtbeta\tran\big)_{kk} \ol{\bbX}_{ki} = \sum_{k = 1}^d \indic_{\{k \neq i\}}
  \bbX_{ik}^2 \wt{\beta}_k^2,\\
  \Big(\bbX \odot \ol{\wtbeta \wtbeta\tran} \odot \bbX\Big)_{ij} &=
  \bbX_{ij}^2 \wt{\beta}_i \wt{\beta}_j.
\end{align*} In turn, for any unit-length vector $\bfv$, \begin{align}
  \bfv\tran \E[B] \bfv &= \E\left[\sum_{i = 1}^d \sum_{k = 1}^d \indic_{\{k
  \neq i\}} v_i^2 \bbX_{ik}^2 \wt{\beta}_k^2 + \sum_{\ell = 1}^d \sum_{m =
  1}^d \Big(\indic_{\{\ell \neq m\}} v_\ell \bbX_{\ell m} \wt{\beta}_\ell\Big)
  \Big(\indic_{\{\ell \neq m\}} \wt{\beta}_m \bbX_{\ell m} v_{m}\Big)\right]
  \nonumber\\
  &= \sum_{\ell = 1}^d \sum_{m = 1}^d \indic_{\{\ell \neq m\}} \bbX_{\ell m}^2
  \E\Big[v_\ell^2 \wt{\beta}_m^2 + v_\ell \wt{\beta}_\ell v_m
  \wt{\beta}_m\Big] \nonumber\\
  &= \dfrac{1}{2} \sum_{\ell = 1}^d \sum_{m = 1}^d \indic_{\{\ell \neq m\}}
  \bbX_{\ell m}^2 \E\Big[\big(v_\ell \wt{\beta}_m + v_m
  \wt{\beta}_\ell\big)^2\Big],
  \label{eq::S_null_id}
\end{align} where the last equality follows by noting that each square $(v_\ell
\wt{\beta}_m + v_m \wt{\beta}_\ell)^2$ appears twice in
\eqref{eq::S_null_id} since the expression is symmetric in $(\ell, m)$. Every
summand in \eqref{eq::S_null_id} is non-negative. If $v_\ell \neq 0$, then there
exists $m(\ell) \neq \ell$ such that $\bbX_{\ell m} \neq 0$. Write $\bfw(\ell)$
for the vector with entries \begin{align*}
  w_i(\ell) = \begin{cases}
    v_{\ell} &\mbox{if } i = m(\ell),\\
    v_{m(\ell)} &\mbox{if } i = \ell,\\
    0 &\mbox{otherwise.}
  \end{cases}
\end{align*} By construction, $\E\big[(v_\ell \wt{\beta}_{m(\ell)} + v_{m(\ell)}
\wt{\beta}_{\ell})^2\big] = \bfw(\ell)\tran \E\big[\wtbeta \wtbeta\tran\big]
\bfw(\ell) \geq \bfw(\ell)\tran \Cov\big(\wtbeta\big) \bfw(\ell)$. Recall that
$\Cov\big(\wtbeta\big) = \bbX_p\inv \bbX \bbX_p\inv$ and note that
$\lambdamin\big(\bbX_p\inv \bbX \bbX_p\inv\big) \geq \lambdamin(\bbX) /
\norm{\bbX_p}^2$. Together with $\norm[\big]{\bfw(\ell)}_2^2 \geq v_\ell^2$
and $\sum_{\ell = 1}^d v_\ell^2 = \norm{\bfv}_2^2 = 1$, \eqref{eq::S_null_id}
now satisfies \begin{align*}
  &\dfrac{1}{2} \sum_{\ell = 1}^d \sum_{m = 1}^d \indic_{\{\ell \neq m\}}
  \bbX_{\ell m}^2 \E\Big[\big(v_\ell \wt{\beta}_m + v_m
  \wt{\beta}_\ell\big)^2\Big]\\
  \geq\ &\dfrac{1}{2} \sum_{\ell = 1}^d \indic_{\{v_\ell \neq 0\}} \bbX_{\ell
  m(\ell)}^2 \bfw(\ell)\tran \Cov\big(\wtbeta\big) \bfw(\ell)\\
  \geq\ &\dfrac{\lambdamin(\bbX)}{2 \norm{\bbX_p}^2} \sum_{\ell
  = 1}^{d} \indic_{\{v_\ell \neq 0\}} \norm[\big]{\bfw(\ell)}_2^2 \min_{m :
  \bbX_{\ell m} \neq 0} \bbX_{\ell m}^2\\
  \geq\ &\dfrac{\lambdamin(\bbX)}{2 \norm{\bbX_p}^2} \min_{i \neq j : \bbX_{ij}
  \neq 0} \bbX_{ij}^2.
\end{align*} As $\lambdamin(S_0) \geq \alpha^2 p^2 (1 - p)^2 \lambdamin(B)$,
this proves the matrix inequality \begin{align}
  \label{eq::S_null_bound}
  S_0 \geq \dfrac{\alpha^2 p^2 (1 - p)^2 \lambdamin(\bbX)}{2 \norm{\bbX_p}^2}
  \min_{i \neq j : \bbX_{ij} \neq 0} \bbX_{ij}^2 \cdot I_d.
\end{align}

Next, let $\bfxi$ be a centered random vector with covariance matrix $M$ and
suppose $D$ is a $d \times d$ dropout matrix, independent of $\bfxi$.
Conditioning on $\bfxi$, the law of total variance states \begin{align}
  &\qquad \Cov\Big(\big(I - \alpha D \bbX_p\big) \bfxi + \alpha D \ol{\bbX} \big(p
  I - D\big) \bfxi\Big) \nonumber\\
  &= \Cov\Big(\E\big[(I - \alpha D \bbX_p) \bfxi + D \ol{\bbX} (p I - D) \bfxi
  \mid \bfxi\big]\Big)\\
  &\qquad + \E\Big[\Cov\big((I - \alpha D \bbX_p) \bfxi + \alpha D \ol{\bbX} (p
  I - D) \bfxi \mid \bfxi\big)\Big] \nonumber\\
  &= \big(I - \alpha p \bbX_p\big) M \big(I - \alpha p \bbX_p\big) + \alpha^2
  \E\bigg[\Cov\Big(D \bbX_p \bfxi + D \ol{\bbX} \big(D - p I\big) \bfxi \mid
  \bfxi\Big)\bigg]. \label{eq::S_lin_mono}
\end{align} Applying Lemma \ref{lem::drop_cov} with $A := \bbX$, $\bfu := \bbX_p
\bfxi$, and $\bfv := \bfxi$ now shows that $S\lin(M) = \Cov\big((I - \alpha D
\bbX_p\big) \bfxi + \alpha D \ol{\bbX} (p I - D) \bfxi\big)$. The second term in
\eqref{eq::S_lin_mono} is always positive semi-definite, proving that $S\lin(M)
\geq \lambdamin(I - \alpha p \bbX_p)^2 \lambdamin(M) \cdot I_d$. As $(\id -
S\lin)\inv = \sum_{\ell = 0}^\infty S\lin^{\ell}$ and $\lambdamin\big(I - \alpha
p \bbX_p\big) = 1 - \alpha p \norm{\bbX_p}$, this implies \begin{align*}
  \big(\id - S\lin\big)\inv M &\geq \left(\lambdamin(M) \sum_{\ell = 0}^\infty
  \lambdamin\big(I - \alpha p \bbX_p\big)^{2 \ell}\right) \cdot I_d\\
  &= \dfrac{\lambdamin(M)}{2 \alpha p \norm{\bbX_p} - (\alpha p)^2
  \norm{\bbX_p}^2} \cdot I_d \geq \dfrac{\lambdamin(M)}{\alpha p \norm{\bbX_p}}
  \cdot I_d.
\end{align*} Lemma \ref{lem::ol_spectral} moreover gives $\norm{\bbX_p} \leq
\norm{\bbX}$. Together with the lower-bound \eqref{eq::S_null_bound} for
$\lambdamin(S_0)$, this proves the result. \CustomQED

\subsection{Proof of Theorem \ref{thm::rp_l2_conv}}

As in Section \ref{sec::conv_wtbeta}, write $\rpbeta_k := k\inv \sum_{j = 1}^k
\wtbeta_j$ for the running average of the iterates and define \begin{align}
  \rpA_k := \E\Big[\big(\rpbeta_k - \wtbeta\big) \big(\rpbeta_k
  -\wtbeta\big)\tran\Big] = \dfrac{1}{k^2} \sum_{j, \ell = 1}^k
  \E\Big[\big(\wtbeta_j -\wtbeta\big)\big(\wtbeta_\ell -\wtbeta\big)\tran\Big].
  \label{eq::rp_ol_A}
\end{align} Suppose $j > \ell$ and take $r = 0, \ldots, j - \ell$. Using
induction on $r$, we now prove that $\E\big[(\wtbeta_j -\wtbeta) (\wtbeta_\ell -
\wtbeta)\tran\big] = (I - \alpha p \bbX_p)^r \E\big[(\wtbeta_{j - r} - \wtbeta)
(\wtbeta_\ell - \wtbeta)\tran\big]$. The identity always holds when $r = 0$.
Next, suppose the identity holds for some $r - 1 < j - \ell$. Taking $k = j + 1
- r$ in \eqref{eq::wtbeta_rewrite}, $\wtbeta_{j + 1 - r} - \wtbeta = \big(I -
\alpha D_{j + 1 - r} \bbX_p\big) \big(\wtbeta_{j - r} - \wtbeta\big) + \alpha
D_{j + 1 - r} \ol{\bbX} \big(p I - D_{j + 1 - r}\big) \wtbeta_{j - r}$. Since $j
- r \geq \ell$, $D_{j + 1 - r}$ is by assumption independent of $\big(\wtbeta,
\wtbeta_{j - r}, \wtbeta_\ell\big)$. Recall from \eqref{eq::exp_ol_rewrite} that
$\E\big[D_{j + 1 - r} \ol{\bbX} (p I - D_{j + 1 - r})\big]=0$. Conditioning on
$\big(\wtbeta, \wtbeta_{j - r}, \wtbeta_\ell\big)$ and applying tower rule now
gives \begin{align*}
  \E\Big[\big(\wtbeta_{j + 1 - r} - \wtbeta\big) \big(\wtbeta_\ell -
  \wtbeta\big)\tran\Big] &= \E\Big[I - \alpha D_{j + 1 - r} \bbX_p\Big]
  \E\Big[\big(\wtbeta_{j - r} - \wtbeta\big) \big(\wtbeta_\ell -
  \wtbeta\big)\tran\Big]\\
  &\qquad + \alpha \E\Big[D_{j + 1 - r} \ol{\bbX} \big(p I - D_{j + 1 -
  r}\big)\Big] \E\Big[\wtbeta_{j - r} \big(\wtbeta_\ell -
  \wtbeta\big)\tran\Big]\\
  &= \big(I - \alpha p \bbX_p\big) \E\Big[\big(\wtbeta_{j - r} -\wtbeta\big)
  \big(\wtbeta_\ell -\wtbeta\big)\tran\Big].
\end{align*} Together with the induction hypothesis, this proves the desired
equality \begin{align*}
  \E\Big[\big(\wtbeta_j - \wtbeta\big) \big(\wtbeta_\ell -\wtbeta\big)\tran\Big]
  &= \big(I - \alpha p\bbX_p\big)^{r - 1} \E\Big[\big(\wtbeta_{j + 1 - r} -
  \wtbeta\big) \big(\wtbeta_\ell - \wtbeta\big)\tran\Big]\\
  &= \big(I - \alpha p \bbX_p\big)^r \E\Big[\big(\wtbeta_{j - r} - \wtbeta\big)
  \big(\wtbeta_\ell -\wtbeta\big)\tran\Big].
\end{align*} For $j < \ell$, transposing and flipping the roles of $j$ and
$\ell$ also shows that $\E\big[(\wtbeta_j - \wtbeta) (\wtbeta_\ell
-\wtbeta)\tran\big] = \E\big[(\wtbeta_j - \wtbeta) (\wtbeta_{\ell - r} -
\wtbeta)\tran\big] (I - \alpha p \bbX_p)^r$ with $r = 0, \ldots, \ell - j$.

Defining $A_\ell := \E\big[(\wtbeta_\ell -\wtbeta) (\wtbeta_\ell
-\wtbeta)\tran\big]$ and taking $r = \abs{j - \ell}$, \eqref{eq::rp_ol_A}
may now be rewritten as \begin{align}
  \rpA_k = \dfrac{1}{k^2} \sum_{j = 1}^k \left(\sum_{\ell = 0}^j \big(I - \alpha
  p \bbX_p\big)^{j - \ell} A_\ell + \sum_{\ell = j + 1}^\infty A_j \big(I -
  \alpha p \bbX_p)^{\ell - j}\right).
  \label{eq::rp_ol_A_exp}
\end{align} Set $\gamma := \norm{I - \alpha p\bbX_p}$, then $\gamma \leq  1 -
\alpha p (1 - p) \min_i \bbX_{ii} < 1$ by Lemma \ref{lem::wt_exp_conv}. Note
also that $\sum_{r = 0}^j \gamma^r \leq \sum_{r = 0}^\infty \gamma^r = (1 -
\gamma)\inv$. Using the triangle inequality and sub-multiplicativity of the
spectral norm, \eqref{eq::rp_ol_A_exp} then satisfies \begin{align*}
  \norm[\big]{\rpA_k} \leq \dfrac{1}{k^2} \sum_{j = 1}^k \left(\sum_{\ell = 0}^j \gamma^{j -
  \ell} \norm[\big]{A_\ell} + \sum_{\ell = j + 1}^\infty \norm[\big]{A_j}
  \gamma^{\ell - j}\right) &\leq \dfrac{2}{k^2} \sum_{\ell = 1}^k
  \norm[\big]{A_\ell} \sum_{r = 0}^\infty \gamma^r\\
  &= \dfrac{2}{k^2 (1 - \gamma)}
  \sum_{\ell = 1}^k \norm[\big]{A_\ell}.
\end{align*} As shown in Theorem \ref{thm::var_lim}, $\norm{A_\ell} \leq
\norm[\big]{(\id - S\lin)\inv S_0} + C \ell \gamma^{\ell - 1}$ for some constant
$C$. Observing that $\sum_{\ell = 1}^\infty \ell \gamma^{\ell - 1} =
\partial_{\gamma} \sum_{\ell = 1}^\infty \gamma^\ell = \partial_{\gamma} \big((1
- \gamma)\inv - 1\big) = (1 - \gamma)^{-2}$, this implies \begin{align*}
  \norm[\big]{\rpA_k} \leq \dfrac{2}{k (1 - \gamma)} \norm[\big]{(\id - S\lin)\inv
  S_0} + \dfrac{2 C}{k^2 (1 - \gamma)^3}.
\end{align*} To complete the proof note that $\gamma \leq  1 - \alpha p (1 - p)
\min_i \bbX_{ii}$ may be rewritten as $(1 - \gamma)\inv \leq \big(\alpha p (1 -
p) \min_{i} \bbX_{ii}\big)\inv$ and $\norm[\big]{(\id - S\lin)\inv S_0} \leq
\alpha p (1 - p)^{-2}  (\min_i \bbX_{ii})^{-3} \norm{\bbX}^2 \cdot
\norm[\big]{\E[X\tran \bfY \bfY\tran X]}$ by \eqref{eq::extra_var_bound}.
\CustomQED

\subsection{Proof of Lemma \ref{lem::wh_div_ker}}

Recall the definition $T(A) := \big(I - \alpha p \bbX\big) A \big(I - \alpha p
\bbX\big) + \alpha^2 p (1 - p) \Diag\big(\bbX A \bbX\big)$.

\begin{lemma}
  \label{lem::wh_cov_rec}
  For every $k = 1, 2, \ldots$ \begin{align*}
    \Cov\big(\whbeta_k - \whbeta\big) &= \big(I - \alpha p \bbX\big)
    \Cov\big(\whbeta_{k - 1} - \whbeta\big) \big(I - \alpha p \bbX\big)\\
    &\qquad + \alpha^2 p (1 - p) \Diag\bigg(\bbX \E\Big[\big(\whbeta_{k - 1} -
    \whbeta\big) \big(\whbeta_{k - 1} - \whbeta\big)\tran\Big] \bbX\bigg)\\
    &\geq T\big(\Cov\big(\whbeta_{k-1} - \whbeta\big)\big),
  \end{align*} with equality if $\whbeta_0 = \E\big[\whbeta\big]$ almost surely.
\end{lemma}

\begin{CustomProof}
  Recall the definition $A_k := \E\big[(\whbeta_k - \whbeta) (\whbeta_k -
  \whbeta)\tran\big]$, so that $\Cov\big(\whbeta_k - \whbeta\big) = A_k -
  \E\big[\whbeta_k - \whbeta\big] \E\big[\whbeta_k - \whbeta\big]\tran$. As
  shown in \eqref{eq::T_rec_def}, $A_k = T(A_{k - 1})$ and hence \begin{align*}
    \Cov\big(\wtbeta_k - \wtbeta\big) &= T\big(A_{k - 1}\big) - \E\big[\whbeta_k
    - \whbeta\big] \E\big[\whbeta_k - \whbeta\big]\tran
  \end{align*} By definition, $T(A_{k - 1}) = \big(I - \alpha p \bbX\big) A_{k -
  1} \big(I - \alpha p \bbX) + \alpha^2 p (1 - p) \Diag\big(\bbX A_{k - 1}
  \bbX\big)$. Recall from \eqref{eq::whbeta_ind_step} that $\E\big[\whbeta_k -
  \whbeta\big] = \big(I - \alpha p \bbX\big) \E\big[\whbeta_{k - 1} -
  \whbeta\big]$, so $\Cov\big(\whbeta_{k - 1} - \whbeta\big) = A_{k - 1} -
  \E\big[\whbeta_{k - 1} - \whbeta\big] \E\big[\whbeta_{k - 1} -
  \whbeta\big]\tran$ implies \begin{align*}
    \big(I - \alpha p \bbX\big)A_{k-1} \big(I - \alpha p \bbX) =\big(I - \alpha
    p \bbX\big)\Cov\big(\whbeta_{k-1} - \whbeta\big) \big(I - \alpha p \bbX)+
    \E[\whbeta_k-\whbeta]\E[\whbeta_k-\whbeta]\tran.
  \end{align*} Together, these identities prove the first claim.

  The lower bound follows from $\bbX \E\big[\whbeta_{k - 1} - \whbeta\big]
  \E\big[\whbeta_{k - 1} - \whbeta\big]\tran \bbX$ being positive semi-definite.
  A positive semi-definite matrix has non-negative diagonal entries, meaning
  $\Diag\big(\bbX \E[\whbeta_{k - 1} - \whbeta] \E[\whbeta_{k - 1}-\whbeta]\tran
  \bbX\big)$ is also positive semi-definite. Next, note that $A_{k - 1} =
  \Cov\big(\whbeta_{k - 1} - \whbeta\big) + \E\big[\whbeta_{k - 1} -
  \whbeta\big] \E\big[\whbeta_{k - 1} - \whbeta\big]\tran$ and in turn
  \begin{equation}
    \label{eq::T_diag_ineq}
    \begin{split}
      \Diag\big(\bbX A_{k - 1} \bbX\big) &= \Diag\big(\bbX \Cov(\whbeta_{k - 1}
      - \whbeta) \bbX\big) + \Diag\big(\bbX \E[\whbeta_{k - 1} - \whbeta]
      \E[\whbeta_{k - 1} - \whbeta]\tran \bbX\big)\\
      &\geq \Diag\big(\bbX \Cov(\whbeta_{k - 1} - \whbeta) \bbX\big).
    \end{split}
  \end{equation} Together with the first part of the lemma and the definition of
  $T$, the lower-bound follows.

  Lastly, if $\whbeta_0 = \E[\whbeta]$ almost surely, then
  \eqref{eq::whbeta_ind_step} implies $\E\big[\whbeta_k - \whbeta\big] = \big(I
  - \alpha p \bbX\big)^k \E\big[\whbeta_0 - \whbeta\big] = 0$, so equality holds
  in \eqref{eq::T_diag_ineq}.
\end{CustomProof}

Consider arbitrary positive semi-definite matrices $A \geq B$, then $\bfw\tran
(A - B) \bfw \geq 0$ for all vectors $\bfw$. Given any vector $\bfv$, this
implies \begin{align*}
  \bfv\tran T(A - B) \bfv &= \bfv\tran \big(I - \alpha p \bbX\big) \big(A -
  B\big) \big(I - \alpha p \bbX\big) \bfv + \alpha^2 p (1 - p) \sum_{\ell = 1}^d
  v_\ell^2 \bfe_\ell\tran \bbX \big(A - B\big) \bbX
  \bfe_\ell\\
  &\geq 0
\end{align*} 
with $\bfe_\ell$ the $\ell$\textsuperscript{th} standard basis vector.
Accordingly, $T$ is operator monotone with respect to the ordering of positive
semi-definite matrices, in the sense that $T(A) \geq T(B)$ whenever $A \geq B$.
Using induction on $k$, Lemma \ref{lem::wh_cov_rec} may now be rewritten as
\begin{align}
  \Cov\big(\whbeta_k - \whbeta\big) \geq T^k\Big(\Cov\big(\whbeta_0 -
  \whbeta\big)\Big).
  \label{eq::T_cov_pow}
\end{align} To complete the proof, the right-hand side of \eqref{eq::T_cov_pow}
will be analyzed for a suitable choice of $\bbX$.

\begin{lemma}
  \label{lem::wh_div_ker_rec}
  Suppose $\whbeta_0$ is independent of all other sources of randomness.
  Consider the linear regression model with a single observation $n = 1$, number
  of parameters $d \geq 2$, and design matrix $X = \bfone\tran$. Then,
  $\Cov\big(\whbeta_k\big) \geq \Cov\big(\whbeta_k - \whbeta\big) +
  \Cov\big(\whbeta\big)$ and for any $d$-dimensional vector $\bfv$ satisfying
  $\bfv\tran \bfone = 0$ and every $k = 1, 2, \ldots$ \begin{align*}
    \bfv\tran \Cov\big(\whbeta_k - \whbeta\big) \bfv\tran \geq
    \alpha^2 p (1 - p) \norm*{\bfv}_2^2.
  \end{align*}
\end{lemma}

\begin{CustomProof}
  By definition, $\bbX = \bfone \bfone\tran$ is the $d \times d$-matrix with all
  entries equal to one. Consequently, $\bbX^k = d^{k - 1} \bbX$ for all $k \geq
  1$ and $\bbX X\tran = \bfone \bfone\tran \bfone = d \bfone = d X\tran$, so
  $\whbeta := d\inv X\tran \bfY$ satisfies the normal equations $X\tran \bfY =
  \bbX \whbeta$. 

  To prove $\Cov\big(\whbeta_k\big) \geq \Cov\big(\whbeta_k - \whbeta\big) +
  \Cov\big(\whbeta\big)$, note that $\Cov\big(\whbeta_k - \whbeta + \whbeta\big)
  \geq \Cov\big(\whbeta_k - \whbeta\big) + \Cov\big(\whbeta\big)$ whenever
  $\Cov\big(\whbeta_k - \whbeta, \whbeta\big) \geq 0$. By conditioning on
  $\big(\whbeta_k, \whbeta\big)$, the identity $\E\big[\whbeta_k - \whbeta\big]
  = \big(I - \alpha p \bbX\big) \E\big[\whbeta_{k - 1} - \whbeta\big]$ was
  shown in \eqref{eq::whbeta_ind_step}. The same argument also proves
  $\E\big[(\whbeta_k - \whbeta) \whbeta\tran\big]  = \big(I - \alpha p \bbX\big)
  \E\big[(\whbeta_{k - 1} - \whbeta) \whbeta\tran\big]$. Induction on $k$
  and the assumed independence between $\whbeta_0$ and $\whbeta$ now lead to
  \begin{align}
    \Cov\big(\whbeta_k - \whbeta, \whbeta\big) &= \E\big[(\whbeta_k - \whbeta)
    \whbeta\tran\big] - \E\big[\whbeta_k - \whbeta\big]
    \E\big[\whbeta\tran\big] \nonumber\\
    &= \big(I - \alpha p \bbX\big) \Cov\big(\whbeta_{k - 1} - \whbeta,
    \whbeta\big) \nonumber\\
    &= \big(I - \alpha p \bbX\big)^k \Cov\big(\whbeta_0 - \whbeta,
    \whbeta\big) \nonumber\\
    &= \big(I - \alpha p \bbX\big)^k \Cov\big(\whbeta\big).
    \label{eq::wh_cross_cov_singular}
  \end{align} Next, note that $\Cov\big(\wtbeta\big) = d^{-2} \bbX$ and
  $\big(I - \alpha p \bbX\big) \bbX = \big(1 - \alpha p d\big) \bbX$. As $\alpha
  p \norm{\bbX} = \alpha p d < 1$, \eqref{eq::wh_cross_cov_singular} satisfies
  \begin{align*}
    \big(I - \alpha p \bbX\big)^k \Cov\big(\whbeta\big) &= \dfrac{1}{d^2} \big(I
    - \alpha p \bbX\big)^k \bbX = \dfrac{1}{d^2} \big(1 - \alpha p d\big)^k
    \bbX \geq 0
  \end{align*} which proves the first claim. 

  To prove the second claim, we first show that there are real sequences
  $\{\nu_k\}_k$ and $\{\lambda_k\}_k$, not depending on the distribution of
  $\whbeta_0$, such that \begin{align}
    \label{eq::wh_cov_ind_lower}
    \Cov\big(\whbeta_k - \whbeta\big) \geq \nu_k I_d + \tfrac{\lambda_k}{d} \bbX
  \end{align} for all $k \geq 0$, with equality if $\whbeta_0 =
  \E\big[\whbeta\big]$ almost surely. When $k = 0$, independence between
  $\whbeta_0$ and $\whbeta$, as well as $\Cov\big(\whbeta\big) = d^{-2} \bbX$,
  imply $\Cov\big(\whbeta_0 - \whbeta\big) = \Cov\big(\whbeta_0\big) +
  \Cov\big(\whbeta\big) \geq d^{-2} \bbX$. Moreover, equality holds whenever
  $\whbeta_0$ is deterministic.

  For the sake of induction suppose the claim is true for some $k - 1$.
  Lemma \ref{lem::wh_cov_rec} and operator monotonicity of $T$ then imply
  \begin{align*}
    \Cov\big(\whbeta_k - \whbeta\big) \geq T\Big(\Cov\big(\whbeta_{k - 1} -
    \whbeta\big)\Big) \geq T\Big(\nu_{k - 1} I_d + \tfrac{\lambda_{k - 1}}{d}
    \bbX\Big).
  \end{align*} In case $\whbeta = \E\big[\whbeta\big]$ almost surely, Lemma
  \ref{lem::wh_cov_rec} and the induction hypothesis assert equality in the
  previous display. Recall $\bbX^\ell = d^{\ell - 1} \bbX$, so $\big(I - \alpha
  p \bbX\big)^2 = I + \big((\alpha p)^2 d - 2 \alpha p \big) \bbX$ as well as
  $\big(I - \alpha p \bbX\big) \bbX \big(I - \alpha p \bbX\big) = (1 - \alpha p
  d)^2 \bbX$. Note also that $\Diag(\bbX) = I_d$. Setting $c := \alpha p d$,
  expanding the definition $T(A) = \big(I - \alpha p \bbX\big) A \big(I -
  \alpha p \bbX\big) + \alpha^2 p (1 - p) \Diag\big(\bbX A \bbX\big)$ now
  results in \begin{align}
    &T\left(\nu_{k - 1} I_d + \dfrac{\lambda_{k - 1}}{d} \bbX\right) \nonumber\\
    =\ &\Big(\nu_{k - 1} + \alpha^2 p (1 - p) d (\nu_{k - 1} + \lambda_{k -
    1})\Big) \cdot I_d + \left(\nu_{k - 1} \Big((\alpha p)^2 d - 2 \alpha p
    \Big) + \dfrac{\lambda_{k - 1} (1 - \alpha p d)^2}{d}\right) \cdot \bbX
    \nonumber\\
    =\ &\Big(\nu_{k - 1} + \alpha^2 p (1 - p) d (\nu_{k - 1} + \lambda_{k -
    1})\Big) \cdot I_d - \dfrac{\nu_{k - 1} - (1 - \alpha p d)^2 (\nu_{k - 1} +
    \lambda_{k - 1})}{d} \cdot \bbX \label{eq::wh_cov_ind_step}
  \end{align} This establishes the induction step and thereby proves
  \eqref{eq::wh_cov_ind_lower} for all $k \geq 0$.

  As $\big(\nu_k, \lambda_k)$ do not depend on the distribution of
  $\whbeta_0$, taking $\whbeta_0 = \E\big[\whbeta\big]$ shows that $\nu_k +
  \lambda_k \geq 0$ for all $k \geq 0$ since \begin{align*}
    0 \leq \bfone\tran \Cov\big(\whbeta_k - \whbeta\big) \bfone = \bfone\tran
    \Big(\nu_k I_d + \tfrac{\lambda_k}{d} \bbX\Big) \bfone = d \big(\nu_k +
    \lambda_k\big).
  \end{align*} Consequently, \eqref{eq::wh_cov_ind_step} implies $\nu_k = \nu_{k
  - 1} + \alpha^2 p (1 - p) d (\nu_{k - 1} + \lambda_{k - 1}) \geq \nu_{k - 1}$,
  proving that $\{\nu_k\}_k$ is non-decreasing in $k$.

  Lastly, we show that $\nu_1 \geq \alpha^2 p (1 - p)$. To this end,
  recall that $\bbX^3 = d^2 \bbX$ and $\Diag(\bbX) = I$. As $T$ is operator
  monotone and $T(A) \geq \alpha^2 p (1 - p) \Diag\big(\bbX A \bbX\big)$,
  independence of $\whbeta_0$ and $\whbeta$ results in \begin{align*}
    \Cov\big(\whbeta_1 - \whbeta\big) &\geq T\Big(\Cov\big(\whbeta_0 -
    \whbeta\big)\Big) \geq T\Big(\Cov\big(\whbeta\big)\Big) \geq \alpha^2 p (1 -
    p) \Diag\big(\bbX d^{-2} \bbX \bbX\big)\\
    &= \dfrac{\alpha^2 p (1 - p)}{d^2} \Diag\big(\bbX^3\big) = \alpha^2 p (1 -
    p) \Diag(\bbX) = \alpha^2 p (1 - p) I_d,
  \end{align*} so $\nu_1 \geq \alpha^2 p (1 - p)$.

  To complete the proof, observe that $\bfv\tran \bfone = 0$ implies
  $\bbX \bfv = 0$. Accordingly, \begin{align*}
     \bfv\tran \Cov\big(\whbeta_k - \whbeta\big) \bfv\tran \geq \bfv\tran
     \Big(\nu_k I_d + \tfrac{\lambda_k}{d} \bbX\Big) \bfv\tran \geq \nu_k
     \norm*{\bfv}_2^2 \geq \nu_1 \norm*{\bfv}_2^2 \geq \alpha^2 p (1 - p)
     \norm*{\bfv}_2^2
  \end{align*} which yields the second claim of the lemma.
\end{CustomProof}

\subsection{Proof of Theorem \ref{thm::wh_l2_conv}}

Recall that $T(A) := \big(I - \alpha p \bbX\big) A \big(I - \alpha p \bbX\big) +
\alpha^2 p (1 - p) \Diag\big(\bbX A \bbX\big)$, as defined in \eqref{eq::T_def}.
If $A_k := \E\big[(\whbeta_k - \whbeta) (\whbeta_k - \whbeta)\tran\big]$, then
$A_k = T\big(A_{k - 1}\big)$ by \eqref{eq::T_rec_def}. 

For an arbitrary $d\times d$ matrix $A$, the triangle inequality, Lemma
\ref{lem::ol_spectral}, and submultiplicativity of the spectral norm imply
\begin{align}
  \label{eq::T_norm_bound}
  \norm[\big]{T(A)} &\leq \Big(\norm[\big]{I - \alpha p \bbX}^2 + \alpha^2 p (1
  - p) \norm*{\bbX}^2\Big) \norm*{A}.
\end{align}
  
As $\bbX$ is positive definite, $\alpha p \norm*{\bbX} < 1$ implies
$\norm[\big]{I - \alpha p \bbX} = 1 - \alpha p \lambdamin(\bbX)$. If $\alpha
\leq \tfrac{\lambdamin(\bbX)}{\norm{\bbX}^2}$, then
\begin{align*}
  \norm[\big]{I - \alpha p \bbX}^2 + \alpha^2 p (1 - p) \norm{\bbX}^2
  &= 1 - 2 \alpha p \lambdamin(\bbX) + \alpha^2 \big(p^2 \lambdamin(\bbX)^2 + p
  (1 - p) \norm{\bbX}^2\big)\\
  &\leq \big(1 -  2 \alpha p \lambdamin(\bbX)\big) + \alpha^2 p \norm{\bbX}^2 \\
  &\leq 1 - \alpha p \lambdamin(\bbX).
\end{align*} Together with \eqref{eq::T_norm_bound} this leads to $\norm{A_k}
\leq \big(1 - \alpha p \lambdamin(\bbX)\big) \norm{A_{k-1}}$. By induction on
$k$, $\norm{A_k} \leq \big(1 -  \alpha p \lambdamin(\bbX)\big)^k \norm{A_0}$,
completing the proof. \CustomQED


\section{Higher Moments of Dropout Matrices}
\label{sec::higher_moments}

Deriving concise closed-form expressions for third and fourth order expectations
of the dropout matrices presents one of the main technical challenges
encountered in Section \ref{sec::gd_results_proof}. 

All matrices in this section will be of dimension $d \times d$ and all vectors
of length $d$. Moreover, $D$ always denotes a random diagonal matrix such that
$D_{ii} \simiid \mathrm{Ber}(p)$ for all $i = 1, \ldots, d$. The diagonal
entries of $D$ are elements of $\{0, 1\}$, meaning $D = D^k$ for all positive
powers $k$. 

Given a matrix $A$ and $p \in (0, 1)$, recall the definitions \begin{align*}
  A_p &:= pA + (1-p) \Diag(A)\\
  \ol{A} &:= A - \Diag(A).
\end{align*} The first lemma contains some simple identities.

\begin{lemma}
  For arbitrary matrices $A$ and $B$, $p \in (0, 1)$, and a diagonal matrix $F$,
  \begin{enumerate}[label=(\alph*), ref=\thelemma(\alph*)]
    \item $\ol{AF} = \ol{A} F$ and $\ol{FA} = F \ol{A}$ \label{lem::ol_prop_A}
    \item $\ol{A}_p = p \ol{A} = \ol{A_p}$ \label{lem::ol_prop_B}
    \item $\Diag\big(\ol{A} B\big) = \Diag\big(A \ol{B}\big)$.
      \label{lem::ol_prop_C}
  \end{enumerate}    
\end{lemma}

\begin{CustomProof}
  \begin{enumerate}[label = (\alph*)]
    \item By definition, $(\ol{A} F)_{ij} = F_{jj} \indic_{\{i \neq j\}} A_{ij}$
      for all $i,j \in \{1, \dots,  d\}$, which equals $\ol{AF}$. The second
      equality then follows by transposition.
    \item Clearly, $\Diag(\ol{A}) = 0$ and in turn $\overline{A}_p = p \ol{A} +
      (1 - p) \Diag(\ol{A}) = p \ol{A}$. On the other hand, $\Diag(A_p)
      = \Diag(A)$ and hence $\ol{A_p} = p A + (1 - p) \Diag(A) - \Diag(A) = p A
      - p \Diag(A)$ equals $p \ol{A}$ as well.
    \item Observe that $\Diag\big(\ol{A} B\big)$ equals $\Diag\big(AB\big) -
      \Diag\big(\Diag(A) B\big)$. As $\Diag\big(\Diag(A) B\big) = \Diag(A)
      \Diag(B) = \Diag\big(A \Diag(B)\big)$, the claim follows.
  \end{enumerate}
\end{CustomProof}

With these basic properties at hand, higher moments involving the dropout matrix
$D$ may be computed by carefully accounting for equalities between the involved
indices.

\begin{lemma}
  \label{lem::drop_moments}
  Given arbitrary matrices $A$, $B$, and $C$, the following hold:
  \begin{enumerate}[label = (\alph*), ref = \thelemma(\alph*)]
    \item $\E\big[DAD\big] = p A_p$ \label{lem::drop_moments_A}
    \item $\E\big[DADBD\big] = p A_p B_p + p^2 (1-p) \Diag(\ol{A} B)$
      \label{lem::drop_moments_B}
    \item $\E\big[DADBDCD\big] = p A_p B_p C_p + p^2 (1 - p)
      \Big(\Diag\big(\ol{A} B_p \ol{C}\big) + A_p \Diag(\ol{B} C) + \Diag(A
      \ol{B}) C_p + (1 - p) A \odot \ol{B}^\top \odot C\Big)$
      \label{lem::drop_moments_C}
  \end{enumerate}
\end{lemma}

\begin{CustomProof}
  \begin{enumerate}[label=(\alph*)]
    \item Recall that $D = D^2$ and hence $\E[D] = \E[D^2] = p I$, meaning
      $\E\big[D \Diag(A) D\big] = \Diag(A) \E[D]=p \Diag(A)$. On the other hand,
      $(D \ol{A} D)_{ij} = D_{ii} D_{jj} A_{ij} \indic_{\{i \neq j\}}$ implies
      $\E\big[D \ol{A} D\big] = p^2 \ol{A}$ due to independence of $D_{ii}$ and
      $D_{jj}$. Combining both identities, \begin{align*}
        \E\big[DAD\big] &= \E\Big[D \ol{A} D + D \Diag(A) D\Big]\\
        &= p^2 \ol{A} + p \Diag(A)\\
        &= p \Big(pA - p \Diag(A) + \Diag(A)\Big) = p A_p.
      \end{align*}
    \item First, note that $D = D^2$ and commutativity of diagonal matrices
      imply \begin{equation}
        \label{eq::D_cube_expanded}
        \begin{split}
          DADBD &= D \ol{A} DBD + D \Diag(A) D BD\\
          &= D \ol{\ol{A} DB} D + \Diag(A)DBD + D\Diag\big(\ol{A} DB\big)D\\
          &= D \ol{\ol{A} DB} D + \Diag(A)DBD + \Diag(\ol{A} DBD).
        \end{split}
      \end{equation} Applying Lemma~\ref{lem::ol_prop_A} twice, $D \ol{\ol{A}
      DB} D = \ol{D \ol{A} DBD}$ has no non-zero diagonal entries. Moreover,
      taking $i,j \in \{1, \ldots, d\}$ distinct, \begin{align*}
        \Big(D \ol{\ol{A} DB} D\Big)_{ij} &= D_{ii} D_{jj} \Big(\ol{\ol{A}
        DB}\Big)_{ij}
        = D_{ii} D_{jj} \sum_{k = 1}^d A_{ik} \indic_{\{i \neq k\}} D_{kk}
        B_{kj},
      \end{align*} so that both $i \neq j$ and $i \neq k$. Therefore, 
      \begin{align*}
        \E\Big[D \ol{\ol{A} DB} D\Big] &= \E[D] \E\Big[\ol{\ol{A} DB} D\Big]\\
        &= p \E\big[\ol{A} DBD\big] - p \Diag\big(\ol{A} \E[DBD]\big)\\
        &= p A \E[DBD] - p \Diag(A) \E[DBD] - p \Diag\big(\ol{A} \E[DBD]\big).
      \end{align*} Reinserting this expression into the expectation of
      \eqref{eq::D_cube_expanded} and applying Part~(i) of the lemma now results
      in the claimed identity \begin{align*}
        \E\big[DADBD\big] &= p A \E[DBD] + (1-p) \Diag(A) \E[DBD] + (1-p)
        \Diag\big(\ol{A} \E[DBD]\big)\\
        &= p \big(p A + (1-p) \Diag(A)\big) B_p + p (1-p) \Diag\big(\ol{A}
        B_p\big)\\
        &= p A_p B_p + p (1-p) \Diag\big(\ol{A} B_p\big)\\
        &= p A_p B_p + p^2 (1-p) \Diag\big(\ol{A}B\big),
      \end{align*}
      where $\Diag(\ol{A} \Diag(B))=0$ by Lemma~\ref{lem::ol_prop_A}.
    \item Following a similar strategy as in Part~(ii), observe that
      \begin{equation}
        \label{eq::D_quart_expanded}
        \begin{split}
          DADBDCD &= DAD \ol{B} DCD + DA D \Diag(B) D CD\\
          &= DAD \ol{B} D \ol{C} D + DAD \Diag(B) CD + DAD \ol{B} D \Diag(C) D\\
          &= D \ol{AD \ol{B}} D \ol{C} D + DAD \Diag(B) CD + DAD \ol{B} \Diag(C)
          D\\
          &\qquad + D \Diag(AD \ol{B}) \ol{C} D.
        \end{split}
      \end{equation} By construction of the latter matrix, \begin{align*}
        \big(D \ol{AD \ol{B}} D \ol{C} D\big)_{ij} &= D_{ii} D_{jj} \sum_{k =
        1}^d \big(\ol{AD \ol{B}}\big)_{ik} \big(D \ol{C}\big)_{kj}\\
        &= D_{ii} D_{jj} \sum_{k = 1}^d \indic_{\{k \neq j\}} D_{kk} C_{kj}
        \sum_{\ell = 1}^d \indic_{\{i \neq k\}} A_{i\ell} D_{\ell\ell}
        \indic_{\{\ell \neq k\}} B_{\ell k},
      \end{align*} meaning $k$ is always distinct from the other indices. The
      $k$ index corresponds to the third $D$-matrix from the left, so this
      proves $\E\big[D \ol{AD \ol{B}} D \ol{C} D\big] = p \E\big[D \ol{AD
      \ol{B}} \, \ol{C} D\big]$. Reversing the overlines of the latter
      expression in order, note that \begin{align*}
        D \ol{AD \ol{B}} \, \ol{C} D &= DADBCD - D \Diag(AD \ol{B}) \ol{C} D -
        DAD \ol{B} \Diag(C) D\\
        &\qquad - DAD \Diag(B) CD.
      \end{align*} Note that the subtracted terms match those added to $D \ol{AD
      \ol{B}} D \ol{C} D$ in \eqref{eq::D_quart_expanded} exactly. In turn,
      these identities prove \begin{align}
        &\E\big[DADBDCD\big] \nonumber\\
        =\ &\E\Big[D \ol{AD \ol{B}} D \ol{C} D + DAD \Diag(B) CD + DAD \ol{B}
        \Diag(C) D + D \Diag(AD \ol{B}) \ol{C} D\Big]\nonumber \\
        =\ &p \E\Big[D \ol{AD \ol{B}} \, \ol{C} D\Big] + \E\Big[DAD \Diag(B) CD
        + DAD \ol{B} \Diag(C) D + D \Diag(AD \ol{B}) \ol{C} D\Big] \nonumber \\
        =\ &p \E\Big[DADBCD\Big] + (1 - p) \E\Big[DAD \Diag(B) CD\Big] \nonumber \\
        &\quad + (1 - p) \E\Big[DAD \ol{B} \Diag(C) D\Big] + (1 - p) \E\Big[D
        \Diag(AD \ol{B}) \ol{C} D\Big]\nonumber \\
        =\ &\E\big[DAD B_p CD\big] + (1 - p) \Big(\E\big[DAD \ol{B} \Diag(C)
        D\big] + \E\big[D \Diag(AD \ol{B}) \ol{C} D\big]\Big).
        \label{eq::mom_sum}
      \end{align} The first and second term in the last equality may be computed
      via Part~(ii) of the lemma, whereas the third term remains to be treated.
     
      By definition, the diagonal entries of $\ol{B}$ and $\ol{C}$ are all zero,
      so $\big(D \Diag(AD \ol{B}) \ol{C} D\big)_{ii}$ equals $0$ for all $i \in
      \{1, \ldots, d\}$. Moreover, taking $i \neq j$ implies \begin{align*}
        \big(D \Diag(AD \ol{B}) \ol{C} D\big)_{ij} &= D_{ii} \big(AD
        \ol{B}\big)_{ii} C_{ij} D_{jj}\\ 
        &= D_{ii} C_{ij} D_{jj} \sum_{k = 1}^d A_{ik} D_{kk} B_{ki} \indic_{\{k
        \neq i\}}.
      \end{align*} On the set $\{i \neq j\} \cap \{i \neq k\}$, the entry
      $D_{ii}$ is independent of $D_{jj}$ and $D_{kk}$. Consequently,
      \begin{align*}
        \E\Big[\big(D \Diag(AD \ol{B}) \ol{C} D\big)_{ij}\Big] &= \E[D_{ii}]
        \sum_{k = 1}^d A_{ik} \ol{B}_{ki} C_{ij} \E[D_{jj} D_{kk}]\\
        &= \sum_{k = 1}^d A_{ik} \ol{B}_{ki} C_{ij} \big(p^3 + p^2 (1 - p)
        \indic_{\{k = j\}}\big)\\
        &= p^3 \sum_{k = 1}^d A_{ik} \ol{B}_{ki}
        C_{ij} + p^2 (1 - p) A_{ij}
        \ol{B}_{ji} C_{ij}.
      \end{align*} In matrix form, the previous equation reads \begin{align*}
        \E\Big[D \Diag\big(AD \ol{B}\big) \ol{C} D\Big] = p^3 \Diag(A \ol{B})
        \ol{C} + p^2 (1-p) A \odot \ol{B}^\top \odot C
      \end{align*} where $\odot$ denotes the Hadamard product.

      Reinserting the computed expressions into \eqref{eq::mom_sum}, as well as
      noting that $\big(\ol{B} \Diag(C)\big)_p = \ol{B}_p \Diag(C)$ and
      $\Diag\big(\ol{A} \, \ol{B} \Diag(C)\big) = \Diag(\ol{A} \, \ol{B})
      \Diag(C)$ yields \begin{equation}
        \label{eq::ful_exp_disp}
        \begin{split}
          \E\big[DADBDCD\big] &= p A_p \big(B_p C\big)_p + p^2 (1 - p)
          \Diag\big(\ol{A} B_p C\big)\\
          &\qquad + p (1 - p) A_p \ol{B}_p \Diag(C) + p^2 (1 - p)^2 \Diag(\ol{A}
          \, \ol{B}) \Diag(C)\\
          &\qquad + p^3 (1 - p) \Diag(A \ol{B}) \ol{C} + p^2 (1 - p)^2 A \odot
          \ol{B}^\top \odot C.
        \end{split}
      \end{equation} Next, using the identity $\ol{B}_p = B_p - \Diag(B)$ of
      Lemma~\ref{lem::ol_prop_B}, we combine the first and third terms of the
      latter display into \begin{align*}
        p A_p \big(B_p C\big)_p + p (1 - p) A_p \ol{B}_p \Diag(C) &= p A_p B_p
        \big(p C + (1 - p) \Diag(C)\big)\\
        &\qquad + p (1 - p) A_p \Diag\big(B_p C\big)\\
        &\qquad - p (1 - p) A_p \Diag(B) \Diag(C)\\
        &= p A_p B_p C_p + p^2 (1 - p) A_p \Diag\big(\ol{B} C\big).
      \end{align*} Regarding the second term of \eqref{eq::ful_exp_disp},
      observe that \begin{align*}
        \Diag\big(\ol{A} B_p C\big) &= \Diag\big(\ol{A} B_p \ol{C}\big) +
        \Diag\big(\ol{A} B_p\big) \Diag(C)\\
        &= \Diag\big(\ol{A} B_p \ol{C}\big) + p\Diag(\ol{A} B) \Diag(C),
      \end{align*} where the second equality follows from $\Diag\big(\ol{A}
      \Diag(B)\big) = 0$. Lastly, $\Diag(\ol{A} \, \ol{B}) = \Diag(\ol{A} B) =
      \Diag(A \ol{B})$ by Lemma \ref{lem::ol_prop_C}, so the fourth and fifth
      term of \eqref{eq::ful_exp_disp} combine into \begin{align*}
        &p \Diag(A \ol{B}) \ol{C} + (1 - p) \Diag(\ol{A} \, \ol{B}) \Diag(C) \\
        =\ &\Diag(A \ol{B}) \big(p \ol{C} + \Diag(C)\big)- p \Diag(A \ol{B})
        \Diag(C)\\
        =\ &\Diag(A \ol{B}) C_p - p \Diag(A \ol{B}) \Diag(C),
      \end{align*} where the common factor $p^2 (1 - p)$ is omitted in the
      display.
      
      Using these identities, Equation \eqref{eq::ful_exp_disp} now turns into
      \begin{align*}
        \E\big[DADBDCD\big] &= p A_p B_p C_p + p^2 (1 - p)^2 A \odot \ol{B}^\top
        \odot C\\
        &\qquad + p^2 (1 - p) \Big(\Diag\big(\ol{A} B_p C\big) - p \Diag(A
        \ol{B})\Diag(C)\Big)\\
        &\qquad + p^2 (1 - p) \Big(A_p \Diag\big(\ol{B} C\big) + \Diag(A \ol{B})
        C_p\Big).
      \end{align*} Noting that $\Diag\big(\ol{A} B_p C\big) - p \Diag(A
      \ol{B})\Diag(C)$ equals $\Diag\big(\ol{A} B_p \ol{C}\big)$ finishes the
      proof.
  \end{enumerate}
\end{CustomProof}

In principle, any computations involving higher moments of $D$ may be
accomplished with the proof strategy of Lemma \ref{lem::drop_moments}. A
particular covariance matrix is needed in Section \ref{sec::gd_results_proof},
which will be given in the next lemma.

\begin{lemma}
  \label{lem::drop_cov}
  Given a symmetric matrix $A$, as well as vectors $\bfu$ and $\bfv$,
  \begin{align*}
    \dfrac{\Cov\big(D \bfu + D \ol{A} (D - pI) \bfv\big)}{p (1 - p)} &=
    \Diag(\bfu \bfu\tran) + p \ol{A} \Diag(\bfv \bfu\tran) + p \Diag(\bfu
    \bfv\tran) \ol{A}\\
    &\qquad + p \big(\ol{A} \Diag(\bfv \bfv\tran) \ol{A}\big)_p + p (1 - p) A
    \odot \ol{\bfv \bfv\tran} \odot A.
  \end{align*}
\end{lemma}

\begin{CustomProof}
  The covariance of the sum is given by
  \begin{align}
    \label{eq::app_cov_sum}
    \Cov\big(D \bfu + D \ol{A} (D - pI) \bfv\big) &= \Cov\big(D
    \bfu\big) + \Cov\big(D \ol{A} (D - pI) \bfv\big) + \Cov\big(D
    \bfu, D \ol{A} (D - pI) \bfv\big)\nonumber\\
    &\qquad + \Cov\big(D \ol{A} (D - pI) \bfv, D\bfu\big)\nonumber\\
    &=: T_1 + T_2 + T_3 + T_4
  \end{align} where each of the latter terms will be treated separately. To this
  end, set $B_\bfu := \bfu \bfu\tran$, $B_\bfv := \bfv \bfv\tran$, $B_{\bfu,
  \bfv}:= \bfu \bfv\tran$ and $B_{\bfv, \bfu} := \bfv \bfu\tran$.

  First, recall that $\E\big[D B_\bfu D \big] = p (B_\bfu)_p$ by Lemma
  \ref{lem::drop_moments_A} and $\E[D]=p$, so the definition of covariance
  implies \begin{align}
    T_1 &= \E\big[DB_\bfu D\big] - \E[D] B_\bfu \E[D] = p^2 B_\bfu + p (1 - p)
    \Diag(B_\bfu) - p^2 B_\bfu \nonumber\\
    &= p (1 - p) \Diag(B_\bfu). \label{eq::app_t1}
  \end{align} Moving on to $T_4$, observe that \begin{align*}
    T_4 &= \Cov\big(D \ol{A} (D - p I) \bfv, D \bfu\big) = \E\big[D \ol{A} (D -
    p I) B_{\bfv, \bfu} D\big] - \E\big[D \ol{A} (D - p I) \bfv\big] \E\big[D
    \bfu\big]\tran.
  \end{align*} By the same argument as in \eqref{eq::exp_ol_rewrite},
  \begin{align}
    \label{eq::app_ix}
    \E\big[D \ol{A} (D - p I) \bfv\big] &= \E\big[D \ol{A} D \big] \bfv - p
    \E\big[D \ol{A}\big] \bfv = p^2 \ol{A} \bfv - p^2 \ol{A} \bfv = 0
  \end{align} so that in turn \begin{align*}
    T_4 = \E\big[D \ol{A} D B_{\bfv, \bfu} D\big] - p \E\big[D \ol{A} B_{\bfv,
    \bfu} D\big].
  \end{align*} Recall $\ol{A}_p = p \ol{A}$ from Lemma \ref{lem::ol_prop_B}.
  Applying Lemma \ref{lem::drop_moments_B} for the first term in the previous
  display and Lemma \ref{lem::drop_moments_A} for the second term now leads to
  \begin{align}
    \label{eq::app_t4}
    T_4 &= p^2 \ol{A} (B_{\bfv, \bfu})_p + p^2 (1 - p) \Diag\big(\ol{A} B_{\bfv,
    \bfu}\big) - p^2 \big(\ol{A} B_{\bfv, \bfu}\big)_p \nonumber\\
    &= p^3 \ol{A} B_{\bfv, \bfu} + p^2 (1 - p) \ol{A} \Diag(B_{\bfv, \bfu}) +
    p^2 (1 - p) \Diag\big(\ol{A} B_{\bfv, \bfu}\big) - p^3 \ol{A} B_{\bfv, \bfu}
    \nonumber\\
    &\qquad - p^2 (1 - p) \Diag\big(\ol{A} B_{\bfv, \bfu}\big) \nonumber\\
    &= p^2 (1 - p) \ol{A} \Diag(B_{\bfv, \bfu}).
  \end{align} Using a completely analogous argument, the reflected term $T_3$
  satisfies \begin{align}
    \label{eq::app_t3}
    T_3 = p^2 (1 - p) \Diag\big(B_{\bfu,\bfv}\big) \ol{A}.
  \end{align}

  The last term $T_2$ necessitates another decomposition into four sub-problems.
  First, recall from \eqref{eq::app_ix} that $\E\big[D \ol{A} (D - p I)
  \bfv\big]$ vanishes, which leads to \begin{align}
    \label{eq::app_t2}
    T_2 &= \Cov\big(D \ol{A} (D - p I) \bfv\big) \nonumber\\
    &= \E\big[D \ol{A} (D - p I) B_\bfv (D - p I) \ol{A} D\big] \nonumber\\
    &= \E\big[D \ol{A} D B_\bfv D \ol{A} D\big] - p \E\big[D \ol{A} B (D - p I)
    \ol{A} D\big] - p \E\big[D \ol{A} (D - p I) B_\bfv \ol{A} D\big] \nonumber\\
    &\qquad + p^2 \E\big[D \ol{A} B_\bfv \ol{A} D\big] \nonumber\\
    &=: T_{2,1} - T_{2,2} - T_{2,3} - T_{2,4}
  \end{align} where the last term is negative since $- p B_\bfv (D - p I) - p (D
  - p I) B_\bfv = - p B_\bfv D - p D B_\bfv + 2 p^2 B_\bfv$. Recall once more
  the identity $\ol{A}_p = p\ol{A}$ from Lemma \ref{lem::ol_prop_B} and apply
  Lemma \ref{lem::drop_moments_C}, to rewrite $T_{2,1}$ as \begin{align*}
    T_{2,1} &= \E\big[D \ol{A} D B_\bfv D \ol{A} D\big]\\
    &= p^3 \ol{A} (B_\bfv)_p \ol{A} + p^2 (1 - p) \Big(\Diag\big(\ol{A}
    (B_\bfv)_p \ol{A}\big) + p \ol{A} \Diag\big(\ol{B_\bfv} A\big) + p
    \Diag\big(A \ol{B_\bfv}\big) \ol{A}\Big)\\
    &\qquad + p^2 (1 - p)^2 A \odot \ol{B_\bfv} \odot A.
  \end{align*} As for $T_{2, 2}$, start by noting that \begin{align*}
    T_{2,2} &= p \E\big[D \ol{A} B_\bfv (D - p I) \ol{A} D\big]\\
    &= p \E\big[D \ol{A} B_\bfv D\ol{A} D\big] - p^2 \E\big[D \ol{A} B_\bfv \ol{A}
    D\big]\\
    &= p^3 \big(\ol{A} B_\bfv\big)_p \ol{A} + p^3 (1 - p) \Diag\big(\ol{\ol{A}
    B_\bfv} \ol{A}\big) - p^3 \big(\ol{A} B_\bfv \ol{A}\big)_p,
  \end{align*} where Lemma~\ref{lem::drop_moments_A} computes the second
  expectation and Lemma~\ref{lem::drop_moments_B} the first expectation. To
  progress, note first the identities \begin{align*}
    \big(\ol{A} B_\bfv\big)_p \ol{A} - \big(\ol{A} B_\bfv \ol{A}\big)_p &=
    (1 - p) \Big(\Diag\big(\ol{A} B_\bfv\big) \ol{A} - \Diag\big(\ol{A} B_\bfv
    \ol{A}\big)\Big)\\
    \Diag\big(\ol{\ol{A} B_\bfv} \ol{A}\big) &= \Diag\big(\ol{A} B_\bfv
    \ol{A}\big)
  \end{align*} so that \begin{align*}
    T_{2,2} &= p^3 (1 - p) \Big(\Diag\big(\ol{A} B_\bfv\big) \ol{A} -
    \Diag\big(\ol{A} B_\bfv \ol{A}\big) + \Diag\big(\ol{A} B_\bfv
    \ol{A}\big)\Big)
    = p^3 (1 - p) \Diag\big(\ol{A} B_\bfv\big) \ol{A}.
  \end{align*} By symmetry, the reflected term $T_{2, 3}$ then also satisfies
  $T_{2,3} = p^3 (1 - p) \ol{A} \Diag\big(B_\bfv \ol{A}\big)$. Lastly, applying
  Lemma \ref{lem::drop_moments_A} to $T_{2,4}$ results in $T_{2,4} = p^3
  \big(\ol{A} B_\bfv \ol{A}\big)_p$. To finish the treatment of $T_2$, inserting
  the computed expressions for $T_{2,1}, T_{2,2}, T_{2,3}$, and $T_{2,4}$ into
  \eqref{eq::app_t2} and combining like terms now leads to \begin{align*}
    T_2 &= p^3 \ol{A} (B_\bfv)_p \ol{A} + p^2 (1 - p)^2 A \odot \ol{B_\bfv}
    \odot A\\
    &\qquad + p^2 (1 - p) \Big(\Diag\big(\ol{A} (B_\bfv)_p \ol{A}\big) + p
    \ol{A} \Diag\big(\ol{B_\bfv} A\big) + \Diag\big(A \ol{B_\bfv}\big) p
    \ol{A}\Big)\\
    &\qquad - p^3 (1 - p) \Diag\big(\ol{A}B_\bfv\big) \ol{A} - p^3 (1 - p) \ol{A}
    \Diag\big(B_\bfv \ol{A}\big) + p^3 \big(\ol{A} B_\bfv \ol{A}\big)_p\\
    &= p^3 \Big(\ol{A} (B_\bfv)_p \ol{A} - \big(\ol{A} B_\bfv \ol{A}\big)_p\Big) +
    p^2 (1 - p) \Diag\big(\ol{A} (B_\bfv)_p \ol{A}\big) + p^2 (1 - p)^2 A \odot
    \ol{B_\bfv} \odot A\\
    &= p^3 (1 - p) \ol{A} \Diag(B_\bfv) \ol{A} + p^2 (1 - p)^2 \Diag\big(\ol{A}
    \Diag(B_\bfv) \ol{A}\big) p^2 (1 - p)^2 A \odot \ol{B_\bfv} \odot A\\
    &= p^2 (1 - p) \big(\ol{A} \Diag(B_\bfv) \ol{A}\big)_p+ p^2 (1 - p)^2 A
    \odot \ol{B_\bfv} \odot A.
  \end{align*} To conclude the proof, insert this expression for $T_{2}$ into
  \eqref{eq::app_cov_sum}, together with $T_1$ as in \eqref{eq::app_t1}, $T_3$
  as in \eqref{eq::app_t3}, and $T_4$ as in \eqref{eq::app_t4} to obtain the
  desired identity \begin{align*}
    \Cov\big(D \bfu + D \ol{A} (D - pI) \bfv\big) &= p (1 - p) \Diag(B_\bfu) +
    p^2 (1 - p) \big(\ol{A} \Diag(B_\bfv) \ol{A}\big)_p\\
    &\qquad + p^2 (1 - p)^2 A \odot \ol{B}_\bfv \odot A\\
    &\qquad + p^2 (1 - p) \big(\Diag(B_{\bfu, \bfv}) \ol{A} + \ol{A}
    \Diag(B_{\bfv,\bfu})\big).
  \end{align*}
\end{CustomProof}


\section{Auxiliary Results}
\label{sec::aux}

Below we collect identities and definitions referenced in other sections. 

\medskip

\noindent \textbf{Neumann series:} Let $V$ denote a real or complex Banach space
with norm $\norm{\ \cdot\ }$. Recall that the operator norm of a linear operator
$\lambda$ on $V$ is given by $\norm*{\lambda}_{\op} = \sup_{\bfv \in V :
\norm{\bfv} \leq 1} \norm*{\lambda(\bfv)}$.

\begin{lemma}[\cite{helemskii_2006}, Proposition~5.3.4]
  \label{lem::neumann}
  Suppose $\lambda : V \to V$ is bounded, linear, and satisfies $\norm*{\id -
  \lambda}_{\op} < 1$, then $\lambda$ is invertible and $\lambda\inv = \sum_{i =
  0}^\infty \big(\id - \lambda\big)^i$. 
\end{lemma}

\medskip

\noindent \textbf{Bounds on singular values:} Recall that $\norm{\ \cdot\ }_2$
denotes the Euclidean norm on $\R^d$ and $\norm{\ \cdot\ }$ the spectral norm on
$\R^{d \times d}$, which is given by the largest singular value
$\sigma_{\max}(\cdot)$. The spectral norm is sub-multiplicative in the sense
that $\norm{AB} \leq \norm{A} \norm{B}$. The spectral norm of a vector $\bfv \in
\R^d$, viewed as a linear functional on $\R^d$, is given by $\norm*{\bfv} =
\norm*{\bfv}_2$, proving that \begin{align}
  \label{eq::tensor_cont}
  \norm[\big]{\bfv \bfw\tran} \leq \norm{\bfv}_2 \norm{\bfw}_2
\end{align} for any vectors $\bfv$ and $\bfw$ of the same length. 

Recall the definitions $\ol{A} = A - \Diag(A)$ and $A_p = p A + (1 - p)
\Diag(A)$ with $p \in (0, 1)$.

\begin{lemma}
  \label{lem::ol_spectral}
  Given $d\times d$ matrices $A$ and $B$, the inequalities
  $\norm[\big]{\Diag(A)} \leq \norm{A}$, $\norm{A_p} \leq \norm{A}$, and
  $\norm{A \odot B} \leq \norm{A} \cdot \norm{B}$ hold. Moreover, if $A$ is
  symmetric and positive semi-definite, then also $\norm[\big]{\ol{A}} \leq
  \norm{A}$.
\end{lemma}

\begin{CustomProof}
  For any matrix $A$, the maximal singular value $\sigma_{\max}(A)$ can be
  computed from the variational formulation $\sigma_{\max}(A) = \max_{\bfv \in
  \R^d \setminus \{\bfzero\}} \norm{A \bfv}_2/\norm{\bfv}_2,$ see
  \cite{horn_johnson_2013}, Theorem 4.2.6.

  Let $\bfe_i$ denote the $i$\textsuperscript{th} standard basis vector. The
  variational formulation of the maximal singular value implies
  $\norm{\Diag(A)}^2 = \max_i A_{ii}^2$ which is bounded by $\max_i \sum_{k =
  1}^d A_{ki}^2 = \max_i \big(A\tran A\big)_{ii} = \max_i \bfe_i\tran A\tran A
  \bfe_i$. The latter is further bounded by $\norm[\big]{A\tran A}$, proving the
  first statement. The second inequality follows from the first since
  \begin{equation*}
    \norm{A_p} \leq p \norm{A} + (1 - p) \norm[\big]{\Diag(A)} \leq  p \norm{A}
    + (1 - p) \norm{A} = \norm{A}.
  \end{equation*} For the inequality concerning the Hadamard product, see
  Theorem 5.5.7 of \cite{horn_johnson_1991}.

  For the last inequality, note that semi-definiteness entails $\min_i
  \Diag(A)_{ii} \geq 0$. Fixing $\bfv \in \R^d$, this ensures $\bfv\tran \ol{A}
  \bfv \leq \bfv\tran A \bfv$, which completes the proof.
\end{CustomProof}

\newpage


\printbibliography

@inproceedings{abadi_et_al_2016,
  title     = {{TensorFlow}: A System for {Large-Scale} Machine Learning},
  author    = {Mart{\'\i}n Abadi and
               Paul Barham and
               Jianmin Chen and
               Zhifeng Chen and
               Andy Davis and
               Jeffrey Dean and
               Matthieu Devin and
               Sanjay Ghemawat and
               Geoffrey Irving and
               Michael Isard and
               Manjunath Kudlur and
               Josh Levenberg and
               Rajat Monga and
               Sherry Moore and
               Derek G. Murray and
               Benoit Steiner and
               Paul Tucker and
               Vijay Vasudevan and
               Pete Warden and
               Martin Wicke and
               Yuan Yu and
               Xiaoqiang Zheng},
  booktitle = {12th USENIX Symposium on Operating Systems Design and
               Implementation},
  pages     = {265--283},
  year      = {2016},
  publisher = {USENIX Association},
  ISBN      = {978-1-931971-33-1}
}

@article{andrews_1988,
  title   = {Laws of Large Numbers for Dependent Non-Identically Distributed
             Random Variables},
  author  = {Andrews, Donald W. K.},
  journal = {Econometric Theory},
  pages   = {458--467},
  volume  = {4},
  number  = {3},
  year    = {1988},
  ISSN    = {0266-4666}
}

@inproceedings{arora_et_al_2021,
  title     = {Dropout: Explicit Forms and Capacity Control},
  author    = {Arora, Raman and
               Bartlett, Peter and
               Mianjy, Poorya and
               Srebro, Nathan},
  booktitle = {38th International Conference on Machine Learning},
  pages     = {351--361},
  year      = {2021},
  publisher = {Proceedings of Machine Learning Research},
}

@inproceedings{ba_frey_2013,
  title     = {Adaptive dropout for training deep neural networks},
  author    = {Ba, Jimmy and
               Frey, Brendan},
  booktitle = {Advances in Neural Information Processing Systems 26},
  pages     = {3084--3092},
  year      = {2013},
  publisher = {Curran Associates, Inc.},
  ISBN      = {978-1-632660-24-4}
}

@article{bah_et_al_2022,
  title   = {Learning deep linear neural networks: {Riemannian} gradient flows
             and convergence to global minimizers},
  author  = {Bah, Bubacarr and
             Rauhut, Holger and
             Terstiege, Ulrich and
             Westdickenberg, Michael},
  journal = {Information and Inference: A Journal of the IMA},
  pages   = {307--353},
  volume  = {11},
  number  = {1},
  year    = {2022},
  ISSN    = {2049-8772}
}

@inproceedings{baldi_sadowski_2013,
  title     = {Understanding Dropout},
  author    = {Baldi, Pierre and
               Sadowski, Peter J.},
  booktitle = {Advances in Neural Information Processing Systems 26},
  pages     = {2814--2822},
  year      = {2013},
  publisher = {Curran Associates, Inc.},
  ISBN      = {978-1-632660-24-4}
}

@article {bartlett_et_al_2023,
  title   = {The Dynamics of Sharpness-Aware Minimization: Bouncing Across
             Ravines and Drifting Towards Wide Minima},
  author  = {Bartlett, Peter L. and
             Long, Philip M. and
             Bousquet, Olivier},
  journal = {Journal of Machine Learning Research},
  volume  = {24},
  year    = {2023},
  number  = {316},
  pages   = {1--36},
  ISSN    = {1533-7928}
}

@misc{bos_schmidt-hieber_2023,
  title        = {Convergence guarantees for forward gradient descent in the
                  linear regression model},
  author       = {{Bos}, Thijs and {Schmidt-Hieber}, Johannes},
  howpublished = {arXiv:2309.15001 [math.ST]},
  year         = {2023},
}

@inproceedings{cavazza_et_al_2018,
  title     = {Dropout as a Low-Rank Regularizer for Matrix Factorization},
  author    = {Cavazza, Jacopo and
               Morerio, Pietro and
               Haeffele, Benjamin and
               Lane, Connor and
               Murino, Vittorio and
               Vidal, Rene},
  booktitle = {21st International Conference on Artificial Intelligence and
               Statistics},
  pages     = {435--444},
  year      = {2018},
  publisher = {Proceedings of Machine Learning Research},
}

@misc{chollet_2015,
  title        = {Keras},
  author       = {Chollet, Fran\c{c}ois and others},
  year         = {2015},
  howpublished = {\url{https://keras.io}}
}

@article{cybenko_1998,
  title   = {Approximation by superpositions of a sigmoidal function},
  author  = {Cybenko, George},
  journal = {Mathematics of Control, Signals, and Systems},
  pages   = {303--314},
  volume  = {2},
  number  = {4},
  year    = {1989},
  ISSN    = {0932-4194}
}

@article{dereich_kassing_2023,
  title   = {Central limit theorems for stochastic gradient descent with
             averaging for stable manifolds},
  author  = {Dereich, Steffen and
             Kassing, Sebastian},
  journal = {Electronic Journal of Probability},
  pages   = {1--48},
  volume  = {28},
  year    = {2023},
  ISSN    = {1083-6489}
}

@incollection{duembgen_et_al_2013,
  title     = {Stochastic search for semiparametric linear regression models},
  author    = {D{\"u}mbgen, Lutz and
               Samworth, Richard J. and
               Schuhmacher, Dominic},
  booktitle = {From Probability to Statistics and Back: High-Dimensional Models
               and Processes. A Festschrift in Honor of Jon A. Wellner},
  pages     = {78--90},
  year      = {2013},
  publisher = {Institute of Mathematical Statistics},
  ISBN      = {978-0-940600-83-6}
}

@book{efron_hastie_2016,
  title     = {Computer Age Statistical Inference. Algorithms, Evidence, and
               Data Science},
  author    = {Efron, Bradley and Hastie, Trevor},
  series    = {Institute of Mathematical Statistics Monographs},
  volume    = {5},
  year      = {2016},
  publisher = {Cambridge University Press},
  ISBN      = {978-1-107-14989-2; 978-1-316-57653-3}
}

@inproceedings{gal_ghahramani_2016a,
  title     = {Dropout as a {Bayesian} Approximation: Representing Model
               Uncertainty in Deep Learning},
  author    = {Gal, Yarin and
               Ghahramani, Zoubin},
  booktitle = {33rd International Conference on International Conference on
               Machine Learning},
  pages     = {1050--1059},
  year      = {2016},
  publisher = {Proceedings of Machine Learning Research}
}

@inproceedings{gal_ghahramani_2016b,
  title     = {A Theoretically Grounded Application of Dropout in Recurrent
               Neural Networks},
  author    = {Gal, Yarin and
               Ghahramani, Zoubin},
  booktitle = {Advances in Neural Information Processing Systems 29},
  pages     = {1027--1035},
  year      = {2016},
  publisher = {Curran Associates, Inc.},
  ISBN      = {978-1-510838-81-9}
}

@article{gao_zhou_2016,
  title   = {Dropout {Rademacher} complexity of deep neural networks},
  author  = {Gao, Wei and
             Zhou, Zhi-Hua},
  journal = {Science China Information Sciences},
  pages   = {2104:1--12},
  volume  = {59},
  number  = {7},
  year    = {2016},
  ISSN    = {1869-1919}
}

@book{goodfellow_et_al_2016,
  title     = {Deep learning},
  author    = {Goodfellow, Ian and
               Bengio, Yoshua and
               Courville, Aaron},
  series    = {Adaptive Computation and Machine Learning},
  year      = {2016},
  publisher = {MIT Press},
  ISBN      = {978-0-262-03561-3; 978-0-262-33743-4}
}

@article{gyorfi_walk_1996,
  title   = {On the Averaged Stochastic Approximation for Linear Regression},
  author  = {Gy\"{o}rfi, L\'{a}szl\'{o}
             and Walk, Harro},
  journal = {SIAM Journal on Control and Optimization},
  pages   = {31--61},
  volume  = {34},
  number  = {1},
  year    = {1996},
  ISSN    = {0363-0129}
}

@book{helemskii_2006,
  title     = {Lectures and Exercises on Functional Analysis},
  author    = {Helemskii, Aleksandr Ya.},
  series    = {Translations of Mathematical Monographs},
  volume    = {233},
  year      = {2006},
  publisher = {American Mathematical Society},
  ISBN      = {0-8218-4098-3}
}

@inproceedings{helmbold_long_2017,
  title     = {Surprising properties of dropout in deep networks},
  author    = {Helmbold, David P. and
               Long, Philip M.},
  booktitle = {Conference on Learning Theory},
  pages     = {1123--1146},
  year      = {2017},
  publisher = {Proceedings of Machine Learning Research}
}

@article{hill_peng_2014,
  title   = {Unified interval estimation for random coefficient autoregressive
             models},
  author  = {Hill, Jonathan and
             Peng, Liang},
  journal = {Journal of Time Series Analysis},
  pages   = {282--297},
  volume  = {35},
  number  = {3},
  year    = {2014},
  ISSN    = {0143-9782}
}

@book{horn_johnson_1991,
  title     = {Topics in Matrix Analysis},
  author    = {Horn, Roger A. and
               Johnson, Charles R.},
  year      = {1991},
  publisher = {Cambridge University Press},
  ISBN      = {0-521-30587-X}
}

@book{horn_johnson_2013,
  title     = {Matrix Analysis},
  author    = {Roger A. {Horn} and
               Charles R. {Johnson}},
  edition   = {2nd},
  year      = {2013},
  publisher = {Cambridge University Press},
  ISBN      = {978-0-521-54823-6; 978-0-521-83940-2}
}

@article{hornik_1991,
  title   = {Approximation capabilities of multilayer feedforward networks},
  author  = {Kurt Hornik},
  journal = {Neural Networks},
  pages   = {251--257},
  volume  = {4},
  number  = {2},
  year    = {1991},
  ISSN    = {0893-6080}
}

@inproceedings{jia_et_al_2014,
  title      = {Caffe: Convolutional Architecture for Fast Feature Embedding},
  author     = {Jia, Yangqing and
                Shelhamer, Evan and
                Donahue, Jeff and
                Karayev, Sergey and
                Long, Jonathan and
                Girshick, Ross and
                Guadarrama, Sergio and
                Darrell, Trevor},
  booktitle = {22nd ACM International Conference on Multimedia},
  pages     = {675--678},
  year      = {2014},
  publisher = {Association for Computing Machinery},
  ISBN      = {978-1-4503-3063-3}
}

@inproceedings{kingma_et_al_2015,
  title     = {Variational Dropout and the Local Reparameterization Trick},
  author    = {Kingma, Diederik P. and
               Salimans, Tim and
               Welling, Max},
  booktitle = {Advances in Neural Information Processing Systems 28}, 
  pages     = {2575--2583},
  year      = {2015},
  publisher = {Curran Associates, Inc.},
  ISBN      = {978-1-510825-02-4}
}

@inproceedings{krizhevsky_et_al_2012,
  title     = {ImageNet Classification with Deep Convolutional Neural Networks},
  author    = {Krizhevsky, Alex and
               Sutskever, Ilya and
               Hinton, Geoffrey E},
  booktitle = {Advances in Neural Information Processing Systems 25},
  pages     = {1097--1105},
  year      = {2012},
  publisher = {Curran Associates, Inc.},
  ISBN      = {978-1-627480-03-1}
}

@article{leshno_et_al_1993,
  title   = {Multilayer feedforward networks with a nonpolynomial activation
             function can approximate any function},
  author  = {Moshe Leshno and
             Vladimir Ya. Lin and
             Allan Pinkus and
             Shimon Schocken},
  journal = {Neural Networks},
  pages   = {861--867},
  volume  = {6},
  number  = {6},
  year    = {1993},
  ISSN    = {0893-6080}
}

@article{manita_et_al_2022,
  title   = {Universal Approximation in Dropout Neural Networks},
  author  = {Oxana A. Manita and
             Mark A. Peletier and
             Jacobus W. Portegies and
             Jaron Sanders and
             Albert Senen-Cerda},
  journal = {Journal of Machine Learning Research},
  pages   = {1--46},
  volume  = {23},
  number  = {19},
  year    = {2022},
  ISSN    = {1533-7928}
}

@misc{mcallester_2013,
  title  = {A {PAC-Bayesian} Tutorial with A Dropout Bound}, 
  author = {David McAllester},
  year = {2013},
  howpublished = {arXiv:1307.2118 [cs.LG]},
}

@inproceedings{mianjy_arora_2019,
  title     = {On Dropout and Nuclear Norm Regularization},
  author    = {Mianjy, Poorya and
               Arora, Raman},
  booktitle = {36th International Conference on Machine Learning},
  pages     = {4575--4584},
  year      = {2019},
  publisher = {Proceedings of Machine Learning Research}
}

@inproceedings{mianjy_arora_2020,
  title = {On Convergence and Generalization of Dropout Training},
  author = {Mianjy, Poorya and Arora, Raman},
  booktitle = {Advances in Neural Information Processing Systems 33},
  pages = {21151--21161},
  year = {2020},
  publisher = {Curran Associates, Inc.},
  ISBN      = {978-1-713829-54-6}
}

@inproceedings{mianjy_et_al_2018,
  title     = {On the Implicit Bias of Dropout},
  author    = {Mianjy, Poorya and
               Arora, Raman and
               Vidal, Rene},
  booktitle = {35th International Conference on Machine Learning},
  pages     = {3540--3548},
  year      = {2018},
  publisher = {Proceedings of Machine Learning Research}
}

@article{moradi_et_al_2020,
  title   = {A Survey of Regularization Strategies for Deep Models},
  author  = {Moradi, Reza and
             Berangi, Reza and
             Minaei, Behrouz},
  journal = {Artificial Intelligence Review},
  pages   = {3947–3986},
  volume  = {53},
  number  = {6},
  year    = {2020},
  ISSN    = {0269--2821}
}

@book{nicholls_quinn_1982,
  title     = {Random Coefficient Autoregressive Models: An Introduction},
  author    = {Nicholls, Des F. and Quinn, Barry G.},
  series    = {Lecture Notes in Statistics},
  volume    = {11},
  year      = {1982},
  publisher = {Springer New York},
  ISBN      = {978-0-387-90766-6}
}

@misc{nguegnang_et_al_2021,
  title        = {Convergence of gradient descent for learning linear neural
                  networks},
  author       = {Nguegnang, Gabin Maxime and
                  Rauhut, Holger and
                  Terstiege, Ulrich},
  howpublished = {arXiv:2108.02040 [cs.LG]},
  year         = {2021}
}

@inproceedings{paszke_et_al_2019,
  title     = {{PyTorch}: {An} Imperative Style, High-Performance Deep Learning
               Library}, 
  author    = {Adam Paszke and
               Sam Gross and
               Francisco Massa and
               Adam Lerer and
               James Bradbury and Gregory Chanan and
               Trevor Killeen and
               Zeming Lin and
               Natalia Gimelshein and
               Luca Antiga and
               Alban Desmaison and
               Andreas Köpf and
               Edward Yang and
               Zach DeVito and
               Martin Raison and
               Alykhan Tejani and
               Sasank Chilamkurthy and
               Benoit Steiner and
               Lu Fang and
               Junjie Bai and
               Soumith Chintala},
  booktitle = {Advances in Neural Information Processing Systems 32},
  year      = {2019},
  publisher = {Curran Associates, Inc.},
  ISBN      = {978-1-713807-93-3}
}

@article{polyak_1990,
  title = {New method of stochastic approximation type},
  author = {Polyak, Boris},
  journal = {Avtomatica i Telemekhanika},
  pages = {98-107},
  volume = {7},
  year = {1990},
  ISSN   = {0005--2310}
}

@article{polyak_juditsky_1992,
  title   = {Acceleration of stochastic approximation by averaging},
  author  = {Boris Polyak and
             Anatoli B. Juditsky},
  journal = {SIAM Journal on Control and Optimization},
  pages   = {838--855},
  volume  = {30},
  number  = {4},
  year    = {1992},
  ISSN    = {0363-0129}
}

@article{regis_et_al_2022,
  title   = {Random autoregressive models: a structured overview},
  author  = {Regis, Marta and
             Serra, Paulo and
             van den Heuvel, Edwin R.},
  journal = {Econometric Reviews},
  pages = {207--230},
  volume = {41},
  number = {2},
  year = {2022},
  ISSN = {0747-4938}
}

@misc{ruppert_1988,
  title        = {Efficient Estimations from a Slowly Convergent {Robbins-Monro}
                  Process},
  author       = {David Ruppert},
  howpublished = {Technical Report 781. Cornell University Operations Research
                  and Industrial Engineering},
  year         = {1988}
}

@article{santos_papa_2022,
  title     = {Avoiding Overfitting: A Survey on Regularization Methods for
               Convolutional Neural Networks},
  author    = {Santos, Claudio Filipi Gon\c{c}alves Dos and
               Papa, Jo\~{a}o Paulo},
  journal   = {ACM Computing Surveys},
  pages     = {213:1--25},
  volume    = {54},
  number    = {10s},
  year      = {2022},
  ISSN      = {0360-0300}
}

@misc{schmidt-hieber_koolen_2023,
  author = {{Schmidt-Hieber}, Johannes and {Koolen}, Wouter M.},
  title = "{Hebbian learning inspired estimation of the linear regression parameters from queries}",
  year = {2023},
  howpublished = {arXiv:2311.03483 [math.ST]}
}

@article{senen-cerda_sanders_2020,
  title   = {Asymptotic Convergence Rate of Dropout on Shallow Linear Neural
             Networks},
  author  = {Senen-Cerda, Albert and
             Sanders, Jaron},
  journal = {Proceedings of the ACM on Measurement and Analysis of Computing
             Systems},
  pages   = {32:1--53},
  volume  = {6},
  number  = {2},
  year    = {2022},
  ISSN    = {2476-1249}
}

@article{srivastava_et_al_2014,
  title   = {Dropout: A Simple Way to Prevent Neural Networks from Overfitting},
  author  = {Nitish Srivastava and
             Geoffrey Hinton and
             Alex Krizhevsky and
             Ilya Sutskever and
             Ruslan Salakhutdinov},
  journal = {Journal of Machine Learning Research},
  pages   = {1929--1958},
  year    = {2014},
  volume  = {15},
  number  = {56},
  ISSN    = {1533-7928}
}

@article{tsigler_bartlett_2023,
  author  = {Alexander Tsigler and
             Peter L. Bartlett},
  title   = {Benign Overfitting in Ridge Regression},
  journal = {Journal of Machine Learning Research},
  year    = {2023},
  volume  = {24},
  number  = {123},
  pages   = {1--76},
  ISSN    = {1533-7928}
}

@inproceedings{wager_et_al_2013,
  title     = {Dropout Training as Adaptive Regularization},
  author    = {Wager, Stefan and
               Wang, Sida and
               Liang, Percy S},
  booktitle = {Advances in Neural Information Processing Systems 26},
  pages     = {351--359},
  year      = {2013},
  publisher = {Curran Associates, Inc.},
  ISBN      = {978-1-632660-24-4}
}

@inproceedings{wan_et_al_2013,
  title     = {Regularization of Neural Networks using DropConnect},
  author    = {Wan, Li and
               Zeiler, Matthew and
               Zhang, Sixin and
               Le Cun, Yann and
               Fergus, Rob},
  booktitle = {30th International Conference on Machine Learning},
  pages     = {1058--1066},
  year      = {2013},
  publisher = {Proceedings of Machine Learning Research}
}

@inproceedings{wang_manning_2013,
  title     = {Fast dropout training},
  author    = {Wang, Sida and
               Manning, Christopher},
  booktitle = {30th International Conference on Machine Learning},
  pages     = {118--126},
  year      = {2013},
  publisher = {Proceedings of Machine Learning Research}
}

@inproceedings{wei_et_al_2020,
  title     = {The Implicit and Explicit Regularization Effects of Dropout},
  author    = {Wei, Colin and
               Kakade, Sham and
               Ma, Tengyu},
  booktitle = {37th International Conference on Machine Learning},
  pages     = {10181--10192},
  year      = {2020},
  publisher = {Proceedings of Machine Learning Research}
}

@article{wu_gu_2015,
  title   = {Towards dropout training for convolutional neural networks},
  author  = {Haibing Wu and Xiaodong Gu},
  journal = {Neural Networks},
  pages   = {1--10},
  volume  = {71},
  year    = {2015},
  ISSN    = {0893-6080}
}

@inproceedings{zhai_wang_2018,
  title     = {Adaptive Dropout with {Rademacher} Complexity Regularization},
  author    = {Zhai, Ke and
               Wang, Huan},
  booktitle = {6th International Conference on Learning Representations},
  year      = {2018}
}

@article{zhang_et_al_2024,
  author  = {Ruiqi Zhang and
             Spencer Frei and
             Peter L. Bartlett},
  title   = {Trained Transformers Learn Linear Models In-Context},
  journal = {Journal of Machine Learning Research},
  year    = {2024},
  volume  = {25},
  number  = {49},
  pages   = {1--55},
  ISSN    = {1533-7928}
}

@article{zhu_et_al_2021,
  title   = {Online Covariance Matrix Estimation in Stochastic Gradient
             Descent},
  author  = {Wanrong Zhu and
             Xi Chen and
             Wei Biao Wu},
  journal = {Journal of the American Statistical Association},
  pages   = {393--404},
  volume  = {118},
  number  = {541},
  year    = {2021},
  ISSN    = {0162-1459}
}


\end{document}